\documentclass[12pt]{article} 
\usepackage{amsmath} 
\usepackage{amssymb} 
\usepackage{amscd}
\usepackage{theorem} 
\usepackage{pb-diagram}  

\def\Onabla{\stackrel{0}{\nabla}}
\def\oplusinf{\mathop{\oplus}} 
\def\otimesinf{\mathop{\otimes}}
 
\def\im{{\mbox{Im}}} 
\def\ker{{\mbox{Ker}}} 
\def\Der{{\mbox{Der}}}
\def\ham{{\mbox{Ham}}} 
\def\aut{{\mbox{Aut}}}

\def\Talpha#1{\vbox{\ialign{##\crcr
    $\alpha$\crcr\noalign{\kern2pt\nointerlineskip}
	   $\hfil\displaystyle{#1}\hfil$\crcr}}}
\def\hom{{\mbox{Hom}}}
\def\ob{{\mbox{Ob}}} \def\fin{\mbox{End}} \def\card{{\mbox{card}}}
\def\dif{{\mathbf{Dif}}} \def\alg{{\mathbf{Alg}}} 
\def\algcom{{\mathbf{Algcom}}}
\def\jord{{\mathbf{Jord}}} \def\lie{{\mathbf{Lie}}} 
\def\ag{{\mathbf{A}}}
\def\cg{{\mathbf{C}}}
\def\slg{\mathfrak{sl}}
\def\sug{\mathfrak{su}}
\def\cali{{\cal I}} 
\def\calj{{\cal J}} 
\def\calp{{\cal P}} 
\def\cala{{\cal A}}
\def\calb{{\cal B}} 
 
\def\calm{{\cal M}} 
\def\caln{{\cal N}}
 
\def\cals{{\cal S}} 
\def\calc{{\cal C}} 
\def\cale{{\cal E}}
\def\cals{{\cal S}}
\def\calu{{\cal U}}
\def\calv{{\cal V}}
 
\def\calr{{\cal R}} 
 
 \def\fracF{\mathfrak F} 
\def\fracT{\mathfrak T}
  
\def\fraca{{\mathfrak A}} 
 
\def\fracg{{\mathfrak g}}
\def\fracc{{\mathfrak C}} 
\def\bbbone{\mbox{\rm 1\hspace {-.6em} l}}

\def\gr{{\mbox{gr}}} 
\def\Int{{\mbox{Int}}}
\def\Out{{\mbox{Out}}}
\def\lied{{\mbox{Lie}}}

\def\diag{{\mbox{\scriptsize Diag}}}
\def\der{{\mbox{\scriptsize Der}}}
\def\derth{\mathrm{\scriptsize Der}}
\def\diagth{\mathrm{\scriptsize Diag}}

\def\Derth{\mathrm{Der}}
\def\obth{\mathrm{Ob}}
\def\homth{\mathrm{Hom}}
\def\os{\underline{\Omega}}
\def\hamth{\mathrm{Ham}}

\newtheorem{lemma}{LEMMA} \newtheorem{proposition}{PROPOSITION}
\newtheorem{corol}{COROLLARY} \newtheorem{theo}{THEOREM}

\begin{document} 

\begin{center} 
{\Large\bf LECTURES  ON}

\vspace{3mm}

{\Large\bf GRADED DIFFERENTIAL ALGEBRAS}

\vspace{3mm}

{\Large\bf  AND NONCOMMUTATIVE GEOMETRY}

\end{center} 
\vspace{0.75cm}

\begin{center} Michel DUBOIS-VIOLETTE \\
\vspace{0.3cm} {\small Laboratoire de Physique Th\'eorique
\footnote{Unit\'e Mixte de Recherche du Centre National de la
Recherche Scientifique - UMR 8627}\\ Universit\'e Paris XI,
B\^atiment 210\\ F-91 405 Orsay Cedex, France\\
patricia$@$th.u-psud.fr}\\ 
\end{center} \vspace{1cm}

\begin{center} \today \end{center}

\vspace {1cm}

\begin{abstract}
These notes contain a survey of some aspects of the theory of graded differential algebras and of noncommutative differential calculi as well as of some applications connected with physics. They also give a description of several new developments.
\end{abstract}

\vspace {2cm}

\noindent LPT-ORSAY 99/100

\vspace {0,5cm}

\noindent {\sl To be published in the
Proceedings of the Workshop on Noncommutative Differential Geometry and its Application to Physics, Shonan-Kokusaimura, Japan, May 31 - June 4, 1999.}
 
\newpage

\tableofcontents

\newpage

\section{Introduction}

The correspondence between ``spaces" and ``commutative
algebras" is  by now familiar in mathematics and in theoretical
physics. This  correspondence allows an algebraic translation
of various geometrical concepts on spaces in  terms of the
appropriate algebras of functions on these spaces. Replacing
these  commutative algebras by noncommutative algebras, i.e.
forgetting commutativity,  leads then to noncommutative
generalizations of geometries  where notions of  ``spaces of
points" are not involved. Such a noncommutative generalization
of  geometry was a need in physics for the formulation of
quantum theory and the  understanding of its relations with
classical physics. In fact, the relation between  spectral
theory and geometry has been implicitely understood very early
in  physics.\\

Gel'fand's transformation associates to each compact
topological  space $X$ the algebra $C(X)$ of complex continuous
functions on $X$. Equipped with  the sup norm, $C(X)$ is a
commutative unital $C^\ast$-algebra. One of the  main points of
Gel'fand theory is that  {\sl the correspondence $X\mapsto
C(X)$  defines an equivalence between the category of compact
topological spaces and  the category of commutative unital
$C^\ast$-algebras.} The compact space $X$ is  then identified
to the spectrum of $C(X)$, (i.e. to the set of  homomorphisms
of unital $\ast$-algebras of $C(X)$ into $\mathbb C$ equipped
with the weak  topology). Let $X$ be a compact space and  let
$\cale (X)$ denote the category of  finite rank complex vector
bundles over $X$. To any vector bundle $E$ of $\cale  (X)$  one
can associate the $C(X)$-module $\Gamma(E)$ of all continuous
sections of  $E$. The module $\Gamma(E)$ is a finite projective
$C(X)$-module and the  Serre-Swan theorem asserts that {\it the
correspondence $E\mapsto \Gamma (E)$  defines an equivalence
between the category $\cale (X)$ and the category $\calp 
(C(X))$ of finite projective $C(X)$-modules}. Thus the compact
spaces and the  complex vector bundles over them can be
replaced by the commutative unital  $C^\ast$-algebras and the
finite projective modules over them. In this sense
noncommutative  unital $C^\ast$-algebras provide
``noncommutative generalizations" of  compact spaces whereas
the notion of finite projective right module  over them is a
corresponding generalization of the notion of complex vector
bundle.  It is worth noticing here that for the latter
generalization one can use as well  left modules but these are
not the only possibilities (see below) and that  something else
has to be used for the generalization of the notion of real
vector  bundle.\\

\noindent \underbar{Remark 1}. Let $X$ be an arbitrary
topological  space, then the algebra $C^{\mathbf b}(X)$ of
complex continuous bounded  functions on $X$ is a
$C^\ast$-algebra if one equips it with the sup norm. In view
of  Gel'fand theory one has $C^{\mathbf b}(X)=C(\hat X)$ (as
$C^\ast$-algebras),  where  $\hat X$ denotes the spectrum of
$C^{\mathbf b}(X)$. The spectrum $\hat X$ is  a compact space
and the evaluation defines a continuous mapping $e:X\mapsto 
\hat X$ with dense image $(\overline{e(X)}=\hat X $). The
compact space $\hat X$  is called {\it the Stone-\u{C}ech
compactification of} $X$ and the pair  $(e,\hat X$) is
characterized  (uniquely up to an isomorphism)  by the
following  universal property: {\it For any continuous mapping
$f:X\mapsto Y$ of $X$ into  a compact space $Y$ there is a
unique continuous mapping $\hat f:\hat X\mapsto  Y$ such that
$f=\hat f\circ  e$}. Notice that $e: X\mapsto \hat X$ is
generally  not injective and that it is an isomorphism, i.e.
$X=\hat X$, if and only if $X$ is  compact. The above universal
property means that $\hat X$ is the biggest  compactification
of $X$ . For instance if $X$ is locally compact then $e$ is 
injective, i.e. $X\subset \hat X$ canonically, but $\hat X$ is
generally much bigger  than the one point compactification
$X\cup\{\infty\}$ of $X$, (e.g. for $X=\mathbb  R$ the
canonical projection $\hat \mathbb R \rightarrow \mathbb R\cup 
\{\infty\}$ has a huge inverse image of $\infty$).\\

If instead of (compact) topological spaces one is interested in
the  geometry of measure spaces, what replaces algebras of
continuous functions are of  course algebras of measurable
functions. In this case the class of algebras  is the class of
commutative $W^\ast$-algebras (or von Neumann algebras). The 
noncommutative generalizations are therefore provided by
general (noncommutative) $W^\ast$-algebras. It has been shown
by A.~Connes that the  corresponding noncommutative measure
theory (i.e. the theory of von Neumann  algebras) has a very
rich structure with no classical (i.e. commutative)
counterpart  (e.g. the occurrence of a canonical dynamical
system) \cite{connes:00}.\\

In the case of differential geometry, it is more or less
obvious that  the appropriate class of commutative algebras are
algebras of smooth  functions. Indeed if $X$ is a smooth
manifold and if $\calc$ is the algebra of  complex smooth
function on $X$, ($\calc=C^\infty(X))$, one can  reconstruct 
$X$ with its smooth structure and the objects attached to $X$,
(differential  forms, etc.), by starting from $\calc$
considered as an abstract (commutative) unital $\ast$-algebra.
As a set $X$ can be identified with the set of  characters of
$\calc$, i.e. with the set of homomorphisms of unital
$\ast$-algebras  of $\calc$ into $\mathbb C$; its differential
structure is connected with the  abundance of derivations of
$\calc$ which identify with the smooth vector fields  on $X$ as
well known. In fact,  in \cite{kosz},  J.L. Koszul gave a
powerful  algebraic generalization of differential geometry in
terms of a commutative  (associative) algebra $\calc$, of
$\calc$-modules and connections (called  derivation laws there)
on these modules. For the applications to differential 
geometry, $\calc$ is of course the algebra of smooth functions
on a smooth manifold and  the $\calc$-modules are modules of
smooth sections of smooth vector  bundles over the manifold.\\

In this approach what generalizes the vector fields are the 
derivations of $\calc$ (into itself). The space $\Der(\calc)$
of all derivations of $\calc$  is a Lie algebra and a
$\calc$-module, both structures being connected by
$[X,fY]=f[X,Y]+X(f)Y$ for $X,Y\in \Der(\calc)$ and $f\in
\calc$.  Using the latter property one can extract, (by
$\calc$-multilinearity), a graded  differential algebra
generalizing the algebra of differential forms, from the 
graded differential algebra $C_\wedge(\Der(\calc),\calc)$ of
$\calc$-valued  Chevalley-Eilenberg cochains of the Lie algebra
$\Der(\calc)$ (with its canonical action  on $\calc$). This
construction admits a generalization to the noncommutative
case;  it is the {\sl derivation-based} differential calculus
(\cite{dv:2},  \cite{dv:3}, \cite{dv:4},\cite{mdv:pm1}
\cite{mdv:pm2}) which will be described  below. As will be
explained (see also \cite{dv:3} and \cite{dv:4}) this is the
right differential calculus for quantum mechanics, in
particular we shall  show that the corresponding noncommutative
symplectic geometry is exactly what is  needed there.\\

For commutative algebras, there is another well-known
generalization  of the calculus of differential forms which is
the K\"ahler differential  calculus \cite{bour}, \cite{hus}, 
\cite{jll}, \cite{ps}. This differential  calculus is
``universal" and consequently functorial for the category of 
(associate unital) commutative algebras. In these lectures we
shall give a  generalization of the K\"ahler differential
calculus for the noncommutative algebras. By  its very
construction, this differential calculus will be functorial for
the algebra-homomorphisms mapping the centers into the centers.
More  precisely this differential calculus will be shown to be
the universal differential  calculus for the category of
algebra $\mathbf{Alg}_Z$ whose objects are the unital 
associative $\mathbb C$-algebras and whose morphisms are the
homomorphisms of  unital algebras mapping the centers into the
centers. This differential calculus  generalizes the K\"ahler
differential calculus in the sense that it reduces to it for  a
commutative (unital associative $\mathbb C$) algebra. This
latter  property is in contrast with what happens for the
so-called {\sl universal  differential calculus}, which is
universal for the category $\mathbf{Alg}$ of  unital
associative $\mathbb C$-algebras and of {\sl all} unital 
algebra-homomorphisms, the construction of which will be
recalled in these lectures.\\

Concerning the generalizations of the notion of module over a 
commutative algebra $\calc$ when one replaces it by a
noncommutative algebra $\cala$,  there are the notion of right
$\cala$-module and the dual notion of left  $\cala$-module, but
since a module over a commutative algebra is also canonically
a  bimodule (of a certain kind) and since a commutative algebra
coincides with its  center, there is a notion of bimodule over
$\cala$ and also the notion of module over  the center
$Z(\cala)$ of $\cala$ which are natural. The ``good choices"
depend  on the kind of problems involved. Again categorial
notions can be of some help.  As will be explained in these
lectures, for each category of algebras there is a  notion of
bimodule over the objects of the category. Furthermore, for
the  category $\mathbf{Algcom}$ of unital commutative
associative $\mathbb  C$-algebras the notion of bimodule just
reduces to the notion of module. Again, like  for the universal
differential calculus, for the notion of bimodule it is 
immaterial for a commutative algebra $\calc$ whether one
considers $\calc$ as an  object of $\mathbf{Algcom}$ or of
$\mathbf{Alg}_Z$ whereas the notion of  bimodule over $\calc$
in $\mathbf{Alg}$ is much wider.\\

This problem of the choice of the generalization of the notion
of  module over a commutative algebra $\calc$ when $\calc$ is
replaced by a  noncommutative algebra $\cala$ is closely
connected with the problem of the noncommutative generalization
of the classical notion of reality. If $\calc$ is the  algebra
of complex continuous functions on a topological space or the
algebra of  complex smooth functions on a smooth manifold, then
it is a $\ast$-algebra  and the (real) algebra of real
functions is the real subspace $\calc^h$ of hermitian  (i.e.
$\ast$-invariant) elements of $\calc$. More generally if
$\calc$ is a  commutative associative complex $\ast$-algebra
the set $\calc^h$ of hermitian  elements of $\calc$ is a
commutative associative real algebra. Conversely if 
$\calc_{\mathbb R}$ is a commutative associative real algebra,
then its  complexification $\calc$ is canonically a commutative
associative complex $\ast$-algebra and  one has
$\calc^h=\calc_{\mathbb R}$. In fact the correspondence
$\calc\mapsto  \calc^h$ defines an equivalence between the
category of commutative  associative complex $\ast$-algebras
and the category of commutative associative real  algebras,
(the morphisms of the first category being the
$\ast$-homomorphisms). This  is in contrast with what happens
for noncommutative algebras. Recall that an associative complex
$\ast$-algebra is an associative complex algebra  $\cala$
equipped with an antilinear involution $x\mapsto x^\ast$ such
that $(xy)^\ast=y^\ast x^\ast$, ($\forall x,y \in \cala$). From
the fact  that the involution reverses the order of the product
it follows that the real  subspace $\cala^h$ of hermitian
elements of a complex associative  $\ast$-algebra is generally
not stable by the product but only by the symmetrized  Jordan
product $x\circ y =\frac{1}{2}(xy+yx)$. Thus $\cala^h$ is not
(generally) an  associative algebra but is a real {\sl Jordan
algebra}. Therefore, one has two  natural choices for the
generalization of an algebra of real functions~:  either the
real Jordan algebra $\cala^h$ of hermitian elements of a
complex  associative $\ast$-algebra $\cala$ which plays the
role of the algebra of complex  functions or a real associative
algebra. In these lectures we take the first  choice which is
dictated by quantum theory (and spectral theory). This choice
has  important consequences on the possible generalizations of
real vector bundles  and, more generally, of modules over
commutative real algebras.\\

Let $\calc$ be a commutative associative $\ast$-algebra and
let  $\calm^h$ be a $\calc^h$-module. The complexified
$\calm=\calm^h\oplus i\calm^h=\calm^h\otimesinf_{\mathbb
R}\mathbb C$ of $\calm^h$ is  canonically a $\calc$-module.
Furthermore there is a canonical antilinear involution
$(\Phi+i\Psi)\mapsto (\Phi+i\Psi)^\ast=\Phi-i\Psi$ 
$(\Phi,\Psi\in  \calm^h)$ for which $\calm^h$ is the set of
$\ast$-invariant elements. This  involution is compatible with
the one of $\calc$ in the sense that one has
$(x\Phi)^\ast=x^\ast\Phi^\ast$ for $x\in \calc$ and $\Phi\in
\calm$;  $\calm$ will be said to be a $\ast$-{\sl module over
the commutative  $\ast$-algebra $\calc$}. In view of the above
discussion what generalizes $\calc$ is a  noncommutative
$\ast$-algebra $\cala$ and we have to generalize the
$\ast$-module  $\calm$ and its ``real part" $\calm^h$. However
it is clear that there is no  noncommutative generalization of
a $\ast$-module over $\cala$ as right or left  module. The
reason is that, since the involution of $\cala$ reverses the
order in  products, it intertwines between actions of $\cala$
and actions of the opposite  algebra $\cala^0$, i.e. between a
structure of right (resp. left) module and  a structure of left
(resp. right) module. Fortunately, as already mentioned, a 
$\calc$-module is canonically a bimodule (of a certain kind)
and the above  compatibility condition can be equivalently
written $(x\Phi)^\ast=\Phi^\ast  x^\ast$. This latter condition
immediately generalizes for $\cala$, namely a {\sl
$\ast$-bimodule over the $\ast$-algebra $\cala$} is a bimodule 
$\calm$ over $\cala$ equipped with an antilinear involution
$\Phi\mapsto\Phi^\ast$  such that $(x\Phi
y)^\ast=y^\ast\Phi^\ast x^\ast$, ($\forall x, y\in \cala$, 
$\forall \Phi\in \calm$). The real subspace $\calm^h=\{\Phi\in
\calm\vert \Phi^\ast=\Phi\}$ of the $\ast$-invariant element of
$\calm$ can play  the role of the sections of a real vector
bundle (for some specific kind of  $\ast$-bimodule $\calm$).
Since a commutative algebra is its center, one can also 
generalize $\ast$-modules over $\calc$ by $\ast$-modules over
the center  $Z(\cala)$ of $\cala$ and modules over $\calc^h$ by
modules over $Z(\cala)^h$. In a  sense these two types of
generalizations of the reality (for modules) are dual
(\cite{mdv:pm1}, \cite{dv:4}) as we shall see later. The main
message  of this little discussion is that notions of reality
force us to consider  bimodules and not only right or left
modules as generalization of vector bundles, \cite{mdv:pm1},
\cite{dv:4}, \cite{connes:05}, \cite{slw}.\\

\noindent\underbar{Remark 2}. One can be more radical. Instead
of  generalizing an associative commutative $\mathbb R$-algebra
$\calc_{\mathbb R}$ by  the Jordan algebra $\cala^h$ of
hermitian elements of an associative complex  $\ast$-algebra
$\cala$, one can more generally choose to generalize
$\calc_{\mathbb  R}$ by a real Jordan algebra $\calj_{\mathbb
R}$ ( not a priori a special  one). The corresponding
generalization of a $\calc_{\mathbb R}$-module could be  then a
{\sl Jordan bimodule over} $\calj_{\mathbb R}$ \cite{nj}
instead of the  real subspace of a $\ast$-bimodule over
$\cala$, (what is a Jordan bimodule will be  explained later).
We however refrain to do that because it is relatively 
complicated technically for a slight generalization
practically.\\

In these lectures we shall be interested in noncommutative
versions of differential geometry where the algebra of smooth
complex functions  on a smooth manifold is replaced by a
noncommutative associative unital complex $\ast$-algebra
$\cala$. Since there are commutative $\ast$-algebras  of this
sort which are not (and cannot be) algebras of smooth functions
on smooth  manifolds, one cannot expect that an arbitrary
$\ast$-algebra as above is a good noncommutative generalization
of an algebra of smooth functions. What  is involved here is
the generalization of the notion of smootheness. It is 
possible to characterize among the unital commutative
associative complex  $\ast$-algebras the ones which are
isomorphic to algebras of smooth functions, however  there are
several inequivalent noncommutative generalizations of this
characterization and no one is universally accepted. Thus
although it is an interesting  subject on which work is
currently in progress \cite {dvkmm}, we shall not  discuss it
here. This means that if the algebra $\cala$ is not ``good
enough", some of our constructions can become a little
trivial.\\

The plan of these notes is the following. After this
introduction, in Section~2 we recall the definition of graded
differential algebras  and of various concepts related to them;
we state in particular the result of D.~Sullivan concerning the
structure of connected finitely generated free graded
comutative differential algebras and we review H.~Cartan's
notion of operation of a Lie algebra in a graded differential
algebra. In Section~3, we explain the equivalence between the
category of finite dimensional Lie algebras and the category of
the free connected graded commutative differential algebras
which are finitely generated in degree 1 (i.e. exterior
algebras of finite dimensional spaces equipped with
differentials); we describe several examples related to Lie
algebras such as the Chevalley-Eilenberg complexes, the Weil
algebra (and we state the result defining the Weil
homomorphism) and we introduce the graded differential algebras
of the derivation-based calculus. In Section~4, we start in an
analogous way as in Section 3, that is we explain the
equivalence between the category of finite dimensional
associative algebras and the category of free connected graded
differential algebras which are generated in degree 1 (i.e.
tensor algebras of finite dimensional spaces equipped with
differentials); we describe examples related to associative
algebras such as Hochschild complexes. In Section~5, we
introduce categories of algebras and we define the associated
notions of bimodules which we follow on several relevant
examples. In Section~6 we recall the notion of first order
differential calculus over an algebra and we introduce our
generalization of the module of K\"ahler differentials and
discuss its functorial properties; we also recall in this
section the definition and properties of the universal first
order calculus. In Section~7 we introduce the higher order
differential calculi and discuss in particular the universal
one as well as our generalization of K\"ahler exterior forms;
we give in particular their universal properties and study
their functorial properties. In Section~8 we introduce another
new differential calculus, the diagonal calculus, which,
although not functorial, is characterized by  a universal
property and we compare it with the other differential calculi
attached to an algebra. In Section~9 we define and study
noncommutative Poisson and symplectic structures and show their
relation with quantum theory. In Section~10 we describe the
theory of connections on modules and on bimodules; in the
latter case we recall in particular the generalization of the
proposal of J. Mourad (concerning linear connections) and
describe its basic properties and its relations with the theory
of first-order operators in bimodules. In Section~11 we discuss
in some examples the relations between connections in the
noncommutative setting and classical Yang-Mills-Higgs models.
Section~12 which serves as conclusion contains some further
remarks concerning in particular the differential calculus on
quantum groups.\\

Apart from in \S 5,  an {\sl algebra} without other
specification shall always mean a unital associative complex
algebra and by a $\ast$-{\sl  algebra} without other
specification we shall mean a unital associative complex
$\ast\mbox{-}$algebra. Given two algebras $\cala$ and $\calb$
in this  sense, a $(\cala,\calb)$-{\sl bimodule} is a vector
space $\calm$ equipped with linear maps $\cala\otimes
\calm\rightarrow \calm$ and $\calm\otimes
\calb\rightarrow~\calm$ denoted by $a\otimes m\mapsto am$ and
$m\otimes b\mapsto mb$ respectively such that $(aa')m=a(a'm)$,
$m(bb')=(mb)b'$,  $(am)b=a(mb)$, $\bbbone m=m$ and
$m\bbbone=m$, $\forall a,a'\in~\cala$, $\forall b,b'\in\calb$,
$\forall m\in \calm$ where $\bbbone$ denotes the unit of
$\cala$ as  well as the one of $\calb$. In Section 5 we shall
define for a more general  algebra $\cala$ a notion of
$\cala$-bimodule which is relative to a category of algebras;
the notion of $(\cala,\cala)$-bimodule as above is the notion
of  $\cala$-bimodule for the category $\alg$ of unital
associative complex algebras. A {\sl  complex} $\fracc$ will be
a $\mathbb Z$-graded vector space (over $\mathbb C$)  equipped
with a homogeneous endomorphism $d$ of degree $\pm 1$ and such
that  $d^2=0$. If $d$ is of degree $-1$, $\fracc$ is said to be
a {\sl chain complex},  its elements are called {\sl chains}
and $d$ is called the {\sl boundary}; if $d$  is of degree
$+1$, $\fracc$ is said to be a {\sl cochain complex}, its
elements  are called {\sl cochains} and $d$ is called the {\sl
coboundary}. The graded  vector space
$H(\fracc)=\ker(d)/\im(d)$ is called the {\sl homology} of
${\fracc}$  if $\fracc$ is a chain complex and the {\sl
cohomology} of $\fracc$ if $\fracc$   is a cochain complex.

\section{Graded differential algebras}

A {\it graded algebra} will be here a unital associative
complex  algebra $\fraca$ which is a $\mathbb Z$-graded vector
space $\fraca=\oplusinf_{n\in  \mathbb Z}\fraca^n$ such that
$\fraca^m. \fraca^n\subset \fraca^{m+n}$. A {\sl homomorphism
of graded algebras} will be a homomorphism of the 
corresponding graded vector spaces (i.e. a homogeneous linear
mapping of degree 0)  which is also a homomorphism of unital
algebras. A graded algebra $\fraca$ is  said to be {\it graded
commutative} if one has $xy=(-1)^{mn}yx$, $\forall x\in 
\fraca^m$ and $\forall y\in \fraca^n$. Most graded algebras
involved in these  lectures will be $\mathbb N$-graded, i.e.
$\fraca^n=0$ for $n\leq -1$. A graded  algebra $\fraca$ is said
to be 0-{\sl connected} or {\sl connected} if it is $\mathbb 
N$-graded with $\fraca^0=\mathbb C\bbbone$, where $\bbbone$
denotes the unit of  $\fraca$. An example of connected graded
algebra is the tensor algebra over  $\mathbb C$ of a complex
vector space $E$ which  will be denoted by $T(E)$. In this 
example, the graduation is the tensorial degree which means
that the degree 1 is  given to the elements of $E$. The
exterior algebra $\bigwedge(E)$ of $E$ is an  example of
connected graded commutative algebra, (the graduation being
again  induced by the tensorial degree).\\

More generally let $C=\oplusinf_n C^n$ be a $\mathbb Z$-graded 
complex vector space and let $T(C)$ be the tensor algebra of
$C$. One has $C\subset  T(C)$ and we equip the algebra $T(C)$
with the unique grading of algebra which  induces on $C$ the
original grading. Since this is not the usual grading of the 
tensor algebra we shall denote the corresponding graded algebra
by $\fracT(C)$. The  graded algebra $\fracT(C)$ is
characterized (uniquely up to an isomorphism)  by the following
universal property: {\sl Any homomorphism of graded vector 
spaces $\alpha:C\rightarrow \fraca$ of the graded vector space
$C$ into a  graded algebra $\fraca$ extends uniquely as a
homomorphism of graded algebras}
$\fracT(\alpha):\fracT(C)\rightarrow \fraca$. Let $\cali$ be
the graded two-sided ideal of $\fracT(C)$ generated by the
graded commutators  $\psi_r\otimes
\varphi_s-(-1)^{rs}\varphi_s\otimes \psi_r$ with $\psi_n$, 
$\varphi_n\in C^n$ and let $\fracF(C)$ denote the quotient
graded algebra  $\fracT(C)/\cali$. Then $\fracF(C)$ is a graded
commutative algebra which contains again $C$  as graded
subspace. The graded commutative algebra $\fracF(C)$ is
characterized  (uniquely up to an isomorphism) by the following
universal property, (which is  the graded commutative
counterpart of the above one): {\sl Any homomorphism of  graded
vector spaces $\alpha:C\rightarrow \fraca$ of the graded vector
space $C$  into a graded commutative algebra $\fraca$ extends
uniquely as a homomorphism of  graded commutative algebras
$\fracF(\alpha):\fracF(C)\rightarrow \fraca$.}  Notice that
$\fracT(C)$ (resp. $\fracF(C)$) is connected if and only if
$C^n=0$  for $n\leq 0$ and that $\fracT(C)=T(C)$ (resp.
$\fracF(C)=\bigwedge(C)$)  as graded  algebras if and only if
$C^n=0$ for $n \not= 1$. Notice also that, as algebra
$\fracF(C)=\bigwedge(\oplusinf_r C^{2r+1})\otimes 
S(\oplusinf_sC^{2s})$ where $S(E)$ denotes the symmetric
algebra of the vector space $E$. The  graded algebra
$\fracT(C)$ will be refered to as {\sl the free graded algebra 
generated by the graded vector space $C$} whereas the graded
algebra $\fracF(C)$ will  be refered to as {\sl the free
graded  commutative algebra generated by the graded vector
space $C$.} Finally, a {\it finitely generated free graded
algebra}  will be a graded algebra of the form $\fracT(C)$ for
some finite dimensional  graded vector space $C$ whereas an
algebra of the form $\fracF(C)$ for some finite  dimensional
graded vector space $C$ will be called a {\sl finitely
generated free  graded commutative algebra}.\\

If $\fraca$ and $\fraca'$ are two graded algebras, their  {\it
tensor  product} $\fraca\otimes \fraca'$ will be here their
{\it skew tensor product}  which means that the product in
$\fraca\otimes \fraca'$ is defined by $(x\otimes  x')(y\otimes
y')=(-1)^{m'n} xy \otimes x'y'$ for $x'\in \fraca^{\prime
m'}$,  $y\in \fraca^n$, $x\in \fraca$ and $y'\in \fraca'$. With
this convention, the tensor  product of two (or more) graded
commutative algebras is again a graded  commutative algebra. If
$C$ and $C'$ are $\mathbb Z$-graded complex vector spaces one
has $\fracF(C\oplus C')=\fracF(C)\otimes \fracF(C')$.\\

By a {\it graded $\ast$-algebra} we here mean a graded algebra
$\fraca=\oplusinf_{n}\fraca^n$ equipped with an involution
$x\mapsto  x^\ast$ satisfying \\ $(i)$ $x\in
\fraca^n\Rightarrow x^\ast \in \fraca^n$  (homogeneity of
degree = 0)\\ $(ii)$ $(\lambda x+y)^\ast=\bar \lambda 
x^\ast+y^\ast,\  \forall x,y\in \fraca $ and $\forall
\lambda\in \mathbb C$ (antilinearity)\\  $(iii)$
$(xy)^\ast=(-1)^{mn}y^\ast x^\ast,\   \forall x\in \fraca^m$
and $\forall y\in\fraca^n$.\\

Notice that Property $(iii)$ implies that if $\fraca$ is
graded  commutative then one has $(xy)^\ast=x^\ast y^\ast$,
$(\forall x,y\in \fraca)$.\\

For a graded algebra $\fraca$, there is, beside the notion of 
derivation, the notion of antiderivation: A linear mapping
$\theta:\fraca\rightarrow  \fraca$ is called an {\it
antiderivation of} $\fraca$ if it satisfies
$\theta(xy)=\theta(x)y+(-1)^m x\theta(y)$ for any $x\in
\fraca^m$ and $y\in\fraca$. However the best generalizations of
the notions of center and of derivations are the following
graded generalizations. The {\it graded  center}
$Z_{\gr}(\fraca)$ of $\fraca$ is the graded subspace of
$\fraca$ generated by the homogeneous elements $x\in\fraca^m$
($m\in \mathbb Z$) satisfying $xy=(-1)^{mn}yx$, $\forall y\in
\fraca^n$ and $\forall n\in \mathbb  Z$, (i.e.
$Z_{\gr}(\fraca)$ is the graded commutant of $\fraca$ in
$\fraca$).  The graded center is a graded subalgebra of
$\fraca$ which is graded  commutative. A {\it graded derivation
of degree $k$ of $\fraca$}, ($k\in\mathbb Z$), is a 
homogeneous linear mapping $X:\fraca\rightarrow \fraca$ which
is of degree $k$  and satisfies $X(xy)=X(x)y+(-1)^{km}xX(y)$
for $x\in\fraca^m$ and $y\in\fraca$.  Thus a homogeneous graded
derivation of even (resp. odd) degree is a  derivation (resp.
antiderivation). The vector space of all these graded
derivations of  degree $k$ will be denoted by
$\Der^k_{\gr}(\fraca)$ and the graded vector space
$\Der_{\gr}(\fraca)=\oplusinf_{k\in\mathbb Z}
\Der^k_{\gr}(\fraca)$  of all graded derivations is a graded
Lie algebra for the {\it graded commutator} $[X,Y]_{\gr}=X
Y-(-1)^{k\ell}Y X$, $X\in\Der^k_{\gr}(\fraca)$,
$Y\in\Der^\ell_{\gr}(\fraca)$. If $x\in\fraca^m$, one defines
a  graded derivation of degree $m$ of $\fraca$, denoted by
$\mbox{ad}_{\gr}(x)$, by setting
$\mbox{ad}_{\gr}(x)y=xy-(-1)^{mn}yx=[x,y]_{\gr}$ for
$y\in\fraca^n$.  The graded subspace of $\Der_{\gr}(\fraca)$
generated by these $\mbox{ad}(x)$,  (when $x$ runs over
$\fraca^m$ and $m$ runs over $\mathbb Z$), is denoted by
$\Int_{\gr}(\fraca)$ and its elements are called {\it inner
graded  derivations of } $\fraca$. It is an ideal of the graded
Lie algebra  $\Der_{\gr}(\fraca)$ and the quotient graded Lie
algebra will be denoted by $\Out_{\gr}(\fraca)$.   Notice that
the graded center $Z_{\gr}(\fraca)$ is stable by the graded 
derivations of $\fraca$ and that this leads to a canonical
homomorphism $\Out_{\gr}(\fraca)\rightarrow
\Der_{\gr}(Z_{\gr}(\fraca))$ since the  inner graded
derivations vanish on $Z_{\gr}(\fraca)$. If $\fraca$ is a 
graded $\ast$-algebra, then $Z_{\gr}(\fraca)$ is stable by the
involution,  (i.e. it is a graded $\ast$-subalgebra of
$\fraca$), one defines in the obvious  manner an involution on
$\Der_{\gr}(\fraca)$ and one has then
$(\mbox{ad}_{\gr}(x))^\ast=-\mbox{ad}_{\gr}(x^\ast)$ for 
$x\in\fraca$. One recovers the usual ungraded notions for an
ordinary (ungraded)  algebra $\cala$ by considering $\cala$ as
a graded algebra which has non zero elements  only in degree
0.\\

Finally a {\it graded differential algebra} is a graded algebra
$\fraca=\oplusinf_n\fraca^n$ equipped with an antiderivation
$d$ of  degree 1 satisfying $d^2=0$, (i.e. $d$ is linear,
$d(xy)=d(x)y+(-1)^mxd(y)$  $\forall x\in \fraca^m$ and $\forall
y\in \fraca$, $d(\fraca^n)\subset  \fraca^{n+1}$ and $d^2=0$);
$d$ is the {\it differential} of the graded differential 
algebra. Notice that then the graded center $Z_\gr(\fraca)$ of
$\fraca$ is stable by the differential $d$ and that it is
therefore a graded differential subalgebra of $\fraca$ which is
graded commutative. A {\it graded differential $\ast$-algebra}
will be a graded  differential algebra $\fraca$ which is also a
graded $\ast$-algebra such that  $d(x^\ast)=(d(x))^\ast,\
\forall x\in \fraca$.\\

Given  a graded differential algebra $\fraca$ its {\it
cohomology}  $H(\fraca)$ is a graded algebra. Indeed the
antiderivation property of $d$ implies  that $\ker(d)$ is a
subalgebra of $\fraca$ and that  $\im (d)$ is a two-sided ideal
of $\ker(d)$ and the homogeneity of $d$ implies  that they are
graded. If $\fraca$ is graded commutative then $H(\fraca)$ is
also  graded commutative and if $\fraca$ is a graded
differential $\ast$-algebra  then $H(\fraca)$ is a graded
$\ast$-algebra.\\

If $\fraca'$ and $\fraca''$ are two graded differential
algebras  their tensor product $\fraca'\otimes \fraca''$ will
be the tensor product of the  graded algebras equipped with the
differential $d$ defined by \[ d(x'\otimes x'')=d(x')\otimes
x''+(-1)^{n'}x'\otimes dx'',\ \forall x'\in  \fraca^{\prime
n'}\ \mbox{and}\ \forall x''\in \fraca''. \] For the
cohomology, one has  the K\"unneth formula \cite{weib} \[
H(\fraca'\otimes\fraca'')=H(\fraca')\otimes  H(\fraca'') \] for
the corresponding graded algebra.\\

\noindent\underbar{Remark 3}. More generally if $\fraca'$ and 
$\fraca''$ are (co)chain complexes of vector spaces with
(co)boundaries denoted by  $d$, then one defines a (co)boundary
$d$ on the graded vector space  $\fraca'\otimes\fraca''$ by the
same formula as above and one has the K\"unneth formula 
$H(\fraca'\otimes \fraca'')=H(\fraca')\otimes H(\fraca'')$ for
the corresponding graded  vector spaces of (co)homologies
\cite{weib}.\\

Let $\fraca$ be a graded differential algebra which is
connected, i.e. such that $\fraca=\mathbb C \bbbone\oplus \fraca^+$ where
$\fraca^+$ is the  ideal of elements of strictly positive
degrees. Then $\fraca$ will be said to  be {\sl minimal} or to
be a {\sl minimal graded differential algebra} if it  satisfies
the {\sl condition of minimality} \cite{sul}: \[ d\fraca
\subset \fraca^+  . \fraca^+\ \ \mbox{(minimal condition).} \]

A {\sl free graded differential algebra} is a graded
differential  algebra which is of the form $\fracT(C)$ for some
graded vector space $C$ as a  graded algebra whereas a {\sl
free graded commutative differential algebra} is a  graded
differential algebra which is of the form $\fracF(C)$ as a
graded  algebra.\\

For instance if $\fracc$ is a cochain complex, its coboundary
extends  uniquely as a differential of $\fracT(\fracc)$ and
also as a differential of $\fracF(\fracc)$. The corresponding
graded differential algebra which  will be again denoted by
$\fracT(\fracc)$ and $\fracF(\fracc)$ when no  confusion arises
will be refered to respectively as {\sl the free graded
differential  algebra generated by the complex $\fracc$} and
{\sl the free graded  commutative differential algebra
generated by the complex $\fracc$}. One can show  (by using the
K\"unneth formula) that one has in cohomology
$H(\fracT(\fracc))=\fracT(H(\fracc))$ and 
$H(\fracF(\fracc))=\fracF(H(\fracc))$. We let the reader guess
the universal properties which characterize $\fracT(\fracc)$
and $\fracF(\fracc)$ and to deduce from these the  functorial
character of the construction. A free graded (resp. graded 
commutative) differential algebra will be said to be {\sl
contractible} if it is  of the form $\fracT(\fracc)$ (resp.
$\fracF(\fracc)$) for a cochain complex (of  vector spaces)
$\fracc$ such that $H(\fracc)=0$ (trivial cohomology). In 
Theorem~1 below we shall be interested in free graded
commutative contractible  differential algebras which are
connected and finitely generated; such a  differential algebra
is a finite tensor product $\otimesinf_\alpha \fracF(\mathbb C 
e_\alpha\oplus \mathbb C de_\alpha)$ with the $e_\alpha$ of
degrees $\geq 1$  (connected property).\\

Concerning the structure of connected finitely generated free
graded  commutative differential algebras, one has the
following result \cite{sul}.

\begin{theo} Every connected finitely generated free graded 
commutative differential algebra is the tensor product of a
unique minimal one  and a unique contractible one. \end{theo}

This result has been for instance an important constructive
ingredient in the computation of the local B.R.S. cohomology of
gauge theory \cite{dvtv}, \cite{dv:1}.\\

There is probably a similar statement for the non graded
commutative  case (i.e. for connected finitely generated free
graded differential algebras)  in which the tensor product is
replaced by the free product of unital algebras.\\

An {\it operation of a Lie algebra $\fracg$ in a graded
differential  algebra} $\fraca$ \cite{cart:01}, \cite{ghv} is a
linear mapping $X\mapsto  i_X$ of $\fracg$ into the space of
antiderivations of degree $-1$ of $\fraca$  such that one has
$(\forall X,Y\in \fracg)$\\ $(i)$ $i_X i_Y +i_Y i_X=0\ 
\mbox{i.e.}\ [i_X,i_Y]_{\gr}=0$\\
$(ii)$$L_Xi_Y-i_YL_X=i_{[X,Y]}\ \mbox{i.e.}\
[L_X,i_Y]_{\gr}=i_{[X,Y]}$ \\ where $L_X$ denotes the
derivation of  degree 0 of $\fraca$ defined by 
\[
L_X=i_X d + d i_X= [d,i_X]_{\gr}
\]
for $X\in  \fraca$. Property $(ii)$ above
implies \\ $(iii)$ $L_XL_Y-L_YL_X=L_{[X,Y]},\   (\forall X,Y\in
\fracg$)\\ which means that $X\mapsto L_X$ is a Lie 
algebra-homomorphism of $\fracg$ into the Lie algebra of
derivations of degree 0 of  $\fraca$. The definition implies
that $L_X$ commutes with the differential $d$ for  any $X\in
\fracg$.\\

Given an operation of $\fracg$ in $\fraca$ as above, an element
$x$  of $\fraca$ is said to be {\it horizontal} if $i_X(x)=0$ 
($\forall X\in  \fracg$), {\it invariant} if $L_X(x)=0$ 
($\forall X\in \fracg$) and {\it basic} if  it is both
horizontal and invariant i.e. if $i_X(x)=0=L_X(x)$  ($\forall
X\in  \fracg$). The set $\fraca_H$ of horizontal elements is a
graded subalgebra of  $\fraca$ stable by the representation
$X\mapsto L_X$ of $\fracg$. The set $\fraca_I$  of invariant
elements is a graded differential subalgebra of $\fraca$ and
the set  $\fraca_B$ of basic elements is a graded differential
subalgebra of $\fraca_I$  (and therefore also of $\fraca$). The
cohomologies of $\fraca_I$ and  $\fraca_B$ are called
respectively {\it invariant cohomology} and {\it basic 
cohomology of} $\fraca$ and are denoted by $H_I(\fraca)$ and
$H_B(\fraca)$.\\

A prototype of graded differential algebra is the graded
differential  algebra $\Omega(M)$ of differential forms on a
smooth manifold $M$. We shall  discuss various generalizations
of it in these lectures. Let $P$ be a smooth  principal bundle
with structure group $G$ and with basis $M$. One defines an 
operation $X\mapsto i_X$ of the Lie algebra $\fracg$ of $G$ in
the graded  differential algebra $\Omega(P)$ of differential
forms on $P$ by letting $i_X$ be  the contraction by the
vertical vector field corresponding to $X\in  \fracg$. Then the
elements of $\Omega(P)_H$ are the horizontal forms in the usual
sense, $\Omega(P)_I$ is the differential algebra of the
differential forms  which are invariant by the action of $G$ on
$P$ whereas the graded differential  algebra $\Omega(P)_B$ is
canonically isomorphic to the graded differential  algebra
$\Omega(M)$ of differential forms on the basis. The
terminology  adopted above for operations comes from this
fundamental example. In \cite{dv:1},  \cite{dv:2} very
different kinds of operations of Lie algebras in graded
differential  algebras have been considered.

\section{Examples related to Lie algebras}

Let $\fracg$ be a finite dimensional complex vector space with
dual  space $\fracg^\ast$. Let $X,Y\mapsto [X,Y]$ be an
antisymmetric bilinear  product on $\fracg$, i.e. a linear
mapping  $[\cdot,\cdot]:\bigwedge^2\fracg\rightarrow \fracg$ of
the second exterior power of $\fracg$ into $\fracg$. The  dual
of the bracket $[\cdot,\cdot]$ is a linear mapping of
$\fracg^\ast$ into
$\bigwedge^2\fracg^\ast(=(\bigwedge^2\fracg)^\ast)$ and such a
linear  mapping of $\fracg^\ast$ into $\bigwedge^2\fracg^\ast$
has a unique extension as  a graded derivation $\delta$ of
degree 1 of the exterior algebra  $\bigwedge\fracg^\ast$.
Conversely, given a graded derivation $\delta$ of degree 1 of
$\bigwedge\fracg^\ast$, the dual of
$\delta:\fracg^\ast\rightarrow \bigwedge^2\fracg^\ast$ is a
bilinear antisymmetric product on $\fracg(=(\fracg^\ast)^\ast)$
and $\delta$ is the unique graded  derivation of degree 1 of
$\bigwedge\fracg^\ast$ which extends the dual of this 
antisymmetric product. Thus to give an antisymmetric product
$[\cdot,\cdot]$ on  $\fracg$ is the same thing as to give a
graded derivation $\delta$ of degree 1 of the  exterior algebra
$\bigwedge\fracg^\ast$. For notational reasons one usually 
introduces the antiderivation $d=-\delta$, i.e. the unique
antiderivation of  $\bigwedge \fracg^\ast$ such that \[
d(\omega)(X,Y)=-\omega([X,Y]) \] for $\omega\in\fracg^\ast$ and
$X,Y\in\fracg$. We shall call $d$ the  antiderivation of
$\bigwedge\fracg^\ast$ corresponding to the bilinear
antisymmetric  product on $\fracg$.

\begin{lemma} The bilinear antisymmetric product
$[\cdot,\cdot]$ on  $\fracg$ satisfies the Jacobi identity if
and only if the corresponding  antiderivation $d$ of
$\bigwedge\fracg^\ast$ satisfies $d^2=0$. \end{lemma} i.e. 
$\fracg$ is a Lie algebra if and only if $\bigwedge\fracg^\ast$
is a graded  differential algebra (for the $d$ corresponding to
the bracket of $\fracg$).\\

\noindent \underbar{Proof}. One has $d^2=\frac{1}{2}[d,d]_{gr}$
so  $d^2$ is a derivation (a graded derivation of degree 2) of
$\bigwedge  \fracg^\ast$. Since, as unital algebra
$\bigwedge\fracg^\ast$ is generated by  $\fracg^\ast$, $d^2=0$
is equivalent to $d^2(\fracg^\ast)=0$. On the other hand by 
definition one has $d(\omega)(X,Y)=-\omega([X,Y])$, for
$\omega\in\fracg^\ast$ and  $X,Y\in \fracg$, and, by the
antiderivation property one has for $X,Y,Z\in \fracg$ \[
3!d^2(\omega)(X,Y,Z)=(d(\omega)(X,[Y,Z])-d(\omega)([X,Y],Z))+\ 
\mbox{cycl}\ (X,Y,Z) \] i.e. $d^2(\omega)(X,Y,Z)=\omega(\left[
[X,Y],Z\right]+\left[[Y,Z],X\right]+\left[[Z,X],Y\right])$.
Therefore $d^2(\omega)=0$ $\forall \omega\in \fracg^\ast$ is
equivalent to the  Jacobi identity for $[\cdot,\cdot]$.
$\square$\\

Thus to give a finite dimensional Lie algebra is the same thing
as to  give the exterior algebra of a finite dimensional vector
space equipped with a differential, that is to give a finitely
generated free graded  commutative differential algebra which
is generated in degree 1. Such a graded  differential algebra
is automatically connected and minimal. This is why, as 
pointed out in \cite{sul}, the connected finitely
generated free graded  commutative differential algebras which
are minimal constitute a natural  categorical closure of finite
dimensional Lie algebras. In fact such generalizations of  Lie
algebras occur in some physical models \cite{sb}.\\

Let $\fracg$ be a finite dimensional Lie algebra, then {\sl
the  cohomology $H(\fracg)$ of $\fracg$} is the cohomology of
$\bigwedge\fracg^\ast$.  More generally, $\bigwedge\fracg^\ast$
is the basic building block to  construct the cochain complexes
for the cohomology of $\fracg$ with values in
representations.\\

Assume that $\fracg$ is the Lie algebra of a Lie group $G$.
Then by  identifying $\fracg$ with the Lie algebra of left
invariant vector fields on $G$  one defines a canonical
homomorphism of $\Lambda\fracg^\ast$ into the graded 
differential algebra $\Omega(G)$ of differential forms on $G$,
(in fact onto the  algebra of left invariant forms). This
induces a homomorphism of $H(\fracg)$ into  the cohomology
$H(G)$ of differential forms on $G$ which is an  isomorphism
when $G$ is compact.\\

In the following, we consider the symmetric algebra
$S\fracg^\ast$, (i.e. the algebra of polynomials on $\fracg$),
to be evenly graded by giving the degree two to its generators,
i.e. by writing $(S\fracg^\ast)^{2n}=S^n\fracg^\ast$ and
$(S\fracg^\ast)^{2n+1}=0$. With this convention $S\fracg^\ast$
is graded commutative and one defines the graded commutative
algebra $W(\fracg)$ by $W(\fracg)=\Lambda\fracg^\ast\otimes
S\fracg^\ast$. Let $(E_\alpha)$ be a basis of $\fracg$ with
dual basis $(E^\alpha)$ and let us define correspondingly
generators $A^\alpha$ and $F^\alpha$ of $W(\fracg)$ by
$A^\alpha=E^\alpha\otimes\bbbone$ and $F^\alpha=\bbbone\otimes
E^\alpha$ so that $W(\fracg)$ is just the free connected graded
commutative algebra (freely) generated by the $A^\alpha$'s in
degree 1 and the $F^\alpha$'s in degree 2. It is convenient to
introduce the elements $A$ and $F$ of $\fracg\otimes W(\fracg)$
defined by $A=E_\alpha\otimes A^\alpha$ and $F=E_\alpha \otimes
F^\alpha$. One then defines the elements $dA^\alpha$ and
$dF^\alpha$ of $W(\fracg)$ by setting \[ \begin{array}{lllll}
dA & = & E_\alpha\otimes dA^\alpha & = & -\frac{1}{2}[A,A]+F\\
dF & = & E_\alpha \otimes dF^\alpha & = & -[A,F] \end{array} \]
where the bracket is the graded Lie bracket obtained by
combining the bracket of $\fracg$ with the graded commutative
product of $W(\fracg)$. One then extends $d$ as an
antiderivation of $W(\fracg)$ of degree 1. One has $d^2=0$, and
since an alternative free system of homogeneous generators of
$W(\fracg)$ is provided by the $A^\alpha$'s and the
$dA^\alpha$'s, $W(\fracg)$ is a connected free graded
commutative differential algebra which is contractible and
which is refered to as the {\sl Weil algebra} of the Lie
algebra $\fracg$ \cite{cart:01}, \cite{ghv}. It is
straightforward to verify that one defines an operation of
$\fracg$ in $W(\fracg)$ by setting $i_X(A^\alpha)=X^\alpha$ and
$i_X(F^\alpha)=0$ for $X=X^\alpha E_\alpha\in \fracg$ and by
extending $i_X$ as an antiderivation of $W(\fracg)$. Since
$W(\fracg)$ is contractible, its cohomology is trivial; the
same is true for the invariant cohomology $H_I(W(\fracg))$ of
$W(\fracg)$, i.e. one has $H^0_I(W(\fracg))=\mathbb C$ and
$H^n_I(W(\fracg))=0$ for $n\geq 1$ \cite{cart:01} (see also in
\cite{dv:1}). The graded subalgebra of horizontal elements of
$W(\fracg)$ is obviously $\bbbone \otimes S\fracg^\ast$ so it
follows that the graded subalgebra of basis elements of
$W(\fracg)$ is just $\bbbone \otimes \cali_S(\fracg)$ where
$\cali_S(\fracg)$ denotes the algebra of invariant polynomials
on $\fracg$ (with the degree 2n given in $W(\fracg)$ to the
homogeneous polynomials of degree n). On the other hand one has
$d(\bbbone\otimes\cali_S(\fracg))=0$ and it is easily seen that
the corresponding homomorphism
$\bbbone\otimes\cali_S(\fracg)\rightarrow H_B(W(\fracg))$ onto 
the basic cohomology of $W(\fracg)$ is an isomorphism.
Therefore, one has $H^{2n}_B(W(\fracg))=\cali^n_S(\fracg)$ and
$H^{2n+1}_B(W(\fracg))=0$, where $\cali^n_S(\fracg)$ denotes
the space of invariant homogeneous polynomials of degree $n$ on
$\fracg$. Let now $P$ be a smooth principal bundle with basis
$M$ and with structure group $G$ such that its Lie algebra is
$\fracg$. One has the canonical operation $X\mapsto i_X$ of
$\fracg$ in $\Omega(P)$ defined at the end of last section.
Given a connection $\omega=E_\alpha\otimes \omega^\alpha\in
\fracg\otimes \Omega^1(P)$ on $P$, there is a unique
homomorphism of graded differential algebras
$\Psi:W(\fracg)\rightarrow \Omega(P)$ such that
$\Psi(A^\alpha)=\omega^\alpha$. This homomorphism satisfies
$\Psi(i_X(w))=i_X(\Psi(w))$ for any $X\in\fracg$ and $w\in
W(\fracg)$. It follows that it induces a homomorphism in basic
cohomomogy $\varphi:H_B(W(\fracg))\rightarrow H_B(P)$, i.e. a
homomorphism of $\cali_S(\fracg)$ into the cohomology $H(M)$ of
the basis $M$ of $P$, such that
$\varphi(\cali^n_S(\fracg))\subset H^{2n}(M)$, (it is an
homormorphism of commutative algebras). One has
$\im(\varphi)\subset H^{ev}(M)=\oplusinf_p H^{2p}(M)$.

\begin{theo} The above homomorphism
$\varphi:\cali_S(\fracg)\rightarrow H^{ev}(M)$ does not depend
on the choice of the connection $\omega$ on $P$. \end{theo}

That is $\varphi$ only depends on $P$; it is called the {\sl
Weil homomorphism} of the principal bundle $P$. Before leaving
this subject, it is worth noticing here that there is a very
interesting noncommutative (or quantized) version of the Weil
algebra of $\fracg$ in the case where $\fracg$ admits a nondegenerate invariant symmetric bilinear form, i.e. for $\fracg$
reductive, where $S\fracg^\ast$ is replaced by the enveloping
algebra $U(\fracg)$ and where $\Lambda\fracg^\ast$ is replaced
by the Clifford algebra $C\ell(\fracg)$ of the bilinear form,
which has been introduced and studied in \cite{am}.\\

In these lectures the Lie algebras involved will be generally
not  finite dimensional and some care must be taken with
respect to duality and  tensor products. For instance, if
$\fracg$ is not finite dimensional then  the dual of the Lie
bracket $[\cdot,\cdot]:\bigwedge^2\fracg\rightarrow \fracg$  is
a linear mapping $\delta:\fracg^\ast\rightarrow
(\bigwedge^2\fracg)^\ast$ and  one only has an inclusion
$\bigwedge^2\fracg^\ast\subset  (\bigwedge^2\fracg)^\ast$. In
the following we give the formulation adapted to this more
general  situation.\\

Let $\fracg$ be a Lie algebra, let $E$ be a representation 
space  of  $\fracg$ (i.e. a $\fracg$-module or, as will be
explained in Section 5, a $\fracg$-bimodule for the category
$\lie$ of Lie algebras) and let  $X\mapsto \pi(X)\in \fin(E)$
denote the action of $\fracg$ on $E$. An $E$-{\sl  valued (Lie
algebra) $n$-cochain of $\fracg$} is a linear mapping 
$X_1\wedge\dots\wedge X_n\mapsto \omega(X_1,\dots,X_n)$ of
$\bigwedge^n\fracg$ into $E$.  The vector space of these
$n$-cochains will be denoted by  $C^n_\wedge(\fracg,E)$. One
defines a homogeneous endomorphism $d$ of degree 1 of the
$\mathbb  N$-graded vector space
$C_\wedge(\fracg,E)=\oplusinf_n C^n_\wedge(\fracg,E)$ of  all
$E$-valued cochains of $\fracg$ by setting \[ \begin{array}{ll}
d(\omega)(X_0,\dots,X_n) & =\sum^n_{k=0}(-1)^
k\pi(X_k)\omega(X_0,\stackrel{k\atop \vee}{\dots},X_n)\\ & 
+\sum_{0\leq r<s\leq
n}(-1)^{r+s}\omega([X_r,X_s],X_0\stackrel{r\atop\vee}{\dots} 
\stackrel{s\atop \vee}{\dots}  X_n) \end{array} \] for
$\omega\in  C^n_\wedge(\fracg,E)$ and $X_i\in\fracg$. It
follows from the Jacobi identity and from
$\pi(X)\pi(Y)-\pi(Y)\pi(X)=\pi([X,Y])$ that $d^2=0$. Thus
equipped  with $d$, $C_\wedge(\fracg,E)$ is a cochain complex
and its cohomology, denoted  by $H(\fracg,E)$, is called the
$E$-{\sl valued cohomology of $\fracg$}.  When $E=\mathbb C$
and $\pi$ is the trivial representation $\pi=0$, it is  the
{\sl cohomology} $H(\fracg)$ of $\fracg$. One verifies that if
$\fracg$ is  finite dimensional, it is the same as the
cohomology of  $\bigwedge\fracg^\ast$; in fact in this case one
has   $C_\wedge(\fracg,E)=E\otimes\bigwedge\fracg^\ast$.\\

Assume now that $E$ is an algebra $\cala$ (unital,
associative,  complex) and that $\fracg$ acts on $\cala$ by
derivations, i.e. that one has
$\pi(X)(xy)=\pi(X)(x)y+x\pi(X)(y)$ for $X\in\fracg$ and 
$x,y\in\cala$. Then $C_\wedge(\fracg,\cala)$ is canonically a
graded differential  algebra. Indeed the product is obtained by
taking the product in $\cala$ after evaluation  and then
antisymmetrizing whereas, the derivation property of the action
of  $\fracg$ implies that $d$ is an antiderivation. The trivial
representation  $\pi=0$ in $\mathbb C$ is of this kind, this is
why $H(\fracg)$ is a graded  algebra.\\

More generally, the vector space $\Der(\cala)$ of all
derivations of $\cala$ into itself is a Lie algebra and
therefore $C_\wedge(\Der(\cala),\cala)$ is a
graded-differential algebra. Furthermore, $\Der(\cala)$ is also
a module over the center $Z(\cala)$ of $\cala$ and one has
$[X,zY]=z[X,Y]+X(z)Y$ from which it follows that the graded
subalgebra $\underline{\Omega}_\der(\cala)$ of
$C_\wedge(\Der(\cala),\cala)$ which consists of
$Z(\cala)$-multilinear cochains is stable by the differential
and is therefore a graded differential subalgebra of
$C_\wedge(\Der(\cala),\cala)$. Since
$\underline{\Omega}^0_\der(\cala)=\cala$, a smaller
differential subalgebra is the smallest differential subalgebra
$\Omega_\der(\cala)$ of $C_\wedge(\Der(\cala),\cala)$
containing $\cala$. When $M$ is a ``good" smooth manifold
(finite dimensional, paracompact, etc.) and $\cala=C^\infty(M)$
then $\underline{\Omega}_\der(\cala)$ and $\Omega_\der(\cala)$
both coincide with the graded differential algebra $\Omega(M)$
of differential forms on $M$. In general, the inclusion
$\Omega_\der(\cala)\subset \underline{\Omega}_\der(\cala)$ is a
strict one even when $\cala$ is commutative (e.g. for the
smooth functions on a $\infty$-dimensional manifold). The
differential calculus over $\cala$ (see in Sections 7, 8) using
$\underline{\Omega}_\der(\cala)$ (or $\Omega_\der(\cala))$ as
generalization of differential forms will be refered to as the
{\sl derivation-based calculus}, \cite{dv:2}, \cite{dv:3},
\cite{dv:4}, \cite{dvkm:1}, \cite{dvkm:2}, \cite{mdv:m2},
\cite{mdv:pm1}, \cite{mdv:pm2}, \cite{mdv:pm3}. If $\cala$ is a
$\ast$-algebra, one defines an involution $X\mapsto X^\ast$ on
$\Der(\cala)$  by setting $X^\ast(a)=(X(a^\ast))^\ast$ and an
involution $\omega\mapsto \omega^\ast$ on
$C_\wedge(\Der(\cala),\cala)$ by setting
$\omega^\ast(X_1,\dots,X_n)=(\omega(X^\ast_1,\dots,X^\ast_n))^\ast$.
So equipped $C_\wedge(\Der(\cala),\cala)$ is a graded
differential $\ast$-algebra and $\Omega_\der(\cala)$ as well as
$\underline{\Omega}_\der(\cala)$ are stable by the involution
and are therefore also graded differential $\ast$-algebras.\\

One defines a linear mapping $X\mapsto i_X$ of $\fracg$ into
the homogeneous endomorphisms of degree $-1$ of
$C_\wedge(\fracg,E)$ by setting
$i_X(\omega)(X_1,\dots,X_{n-1})=\omega(X,X_1,\dots,X_{n-1})$
for  $\omega\in C^n_\wedge(\fracg,E)$ and $X_i\in\fracg$. Then
$X\mapsto  L_X=i_Xd+di_X$ is a representation of $\fracg$ in
$C_\wedge(\fracg,E)$ by homogeneous  endomorphisms of degree 0
which extends the original representation $\pi$ in
$E=C^0_\wedge(\fracg,E)$, i.e. $L_X\restriction E=\pi(X)$ for 
$X\in\fracg$. In the case where $E$ is an algebra $\cala$ and
where $\fracg$ acts by  derivations on $\cala$, we have seen
that $C_\wedge(\fracg,\cala)$ is a graded  differential algebra
and it is easy to show that $X\mapsto i_X$ is an operation of 
the Lie algebra $\fracg$ in the graded differential algebra 
$C_\wedge(\fracg,\cala)$; in fact properties $(i)$ and $(ii)$
of operations (see last section)  hold already in
$C_\wedge(\fracg,E)$ for any $\fracg$-module $E$.\\

In particular one has the operation $X\mapsto i_X$ of the Lie
algebra $\Der(\cala)$ in the graded differential algebra
$C_\wedge(\Der(\cala),\cala)$ defined as above. It is not hard
to verify that the graded differential subalgebras
$\os_\der(\cala)$ and $\Omega_\der(\cala)$ are stable by the
$i_X$ ($X\in\Der(\cala)$). The corresponding operations will be
refered to as {\sl the canonical operations of $\Derth(\cala)$
in $\os_\derth(\cala)$ and in $\Omega_\derth(\cala)$}.

\section{Examples related to associative algebras}

Let $\cala$ be a finite dimensional complex vector space with
dual  space $\cala^\ast$ and let $x,y\mapsto xy$ be an
arbitrary bilinear product  on $\cala$, i.e. a linear mapping
$\otimes^2\cala\rightarrow \cala$ where  $\otimes^2\cala$
denotes the second tensor power of $\cala$. The dual of the
product  is a linear mapping of $\cala^\ast$ into
$\otimes^2\cala^\ast$ and again such a  linear mapping uniquely
extends as a graded derivation $\delta$ of degree~1  of the
tensor algebra $T(\cala^\ast)=\displaystyle{\oplusinf_{n\geq
0}}\otimes^n\cala^\ast$. Conversely, given such a graded
derivation  $\delta$ of degree 1 (i.e. an antiderivation of
degree~1) of $T(\cala^\ast)$, the  dual mapping of the
restriction $\delta:\cala^\ast\rightarrow  \otimes^2\cala^\ast$
of $\delta$ to $\cala^\ast$ is a bilinear product on $\cala$
which is  such that $\delta$ is obtained from it by the above
construction. Thus, to give  a bilinear product on $\cala$ is
the same thing as to give an antiderivation of  degree 1 of
$T(\cala^\ast)$. Again, for notational reasons, it is usual to 
consider the antiderivation $d=-\delta$, i.e. the unique
antiderivation of  $T(\cala^\ast)$ such that \[
d(\omega)(x,y)=-\omega(xy) \] for $\omega\in \cala^\ast$  and
$x,y\in\cala$. We shall call this $d$ the antiderivation of 
$T(\cala^\ast)$ corresponding to the bilinear product of
$\cala$.

\begin{lemma} The bilinear product on $\cala$ is associative if
and  only if the corresponding antiderivation of
$T(\cala^\ast)$ satisfies $d^2=0$.  \end{lemma} i.e. $\cala$ is
an associative algebra if and only if $T(\cala^\ast)$  is a
graded differential algebra (for the $d$ corresponding to the
product of  $\cala$).\\

\noindent\underbar{Proof}. By definition, one has for
$\omega\in  \cala^\ast$ and $x,y,z\in \cala$ \[ 
d(d(\omega))(x,y.z)=d(\omega)(x,yz)-d(\omega)(xy,z)=
\omega((xy)z-x(yz)). \] Therefore the product of $\cala$ is 
associative if and only if $d^2$ vanishes on $\cala^\ast$ but
this is equivalent to  $d^2=0$ since $d^2$ is a derivation and
since the (unital) graded algebra  $T(\cala^\ast)$ is generated
by $\cala^\ast$. $\square$\\

Therefore to give a finite dimensional associative algebra is
the  same thing as to give a finitely generated free graded
differential algebra which  is generated in degree 1. Again
such a graded differential algebra is  automatically connected
and minimal. The situation is very similar to the one of last
section  except that here one has not graded commutativity. So
one can consider in  particular that the connected finitely
generated free graded differential algebras which  are minimal
constitute a natural categorical closure of finite dimensional 
associative algebras, i.e. a natural generalization of the
notion of associative  algebra.\\

Let $\cala$ be a finite dimensional associative algebra; we
shall see  that if $\cala$ has a unit then the cohomology of
the graded differential  algebra $T(\cala^\ast)$ is trivial.
Nevertheless $T(\cala^\ast)$ is the basic  building block to
construct the Hochschild cochain complexes. Namely if  $\calm$
is a $(\cala,\cala)$-bimodule then the graded vector space of 
$\calm$-valued Hochschild cochains of $\cala$ is the graded
space $\calm\otimes  T(\cala^\ast)$ and the Hochschild
coboundary $d_H$ is given by \[ \begin{array}{ll}
d_H(\omega)(x_0,\dots,x_n) & = x_0\omega(x_1,\dots,x_n) + 
(I_\calm\otimes d)(\omega)(x_0,\dots,x_n)\\ & +
(-1)^{n+1}\omega(x_0,\dots,x_{n-1})x_n \end{array} \] for
$\omega\in \calm\otimes(\otimes^n\cala^\ast)$ and
$x_i\in\cala$.\\

In these lectures we shall have to deal with infinite
dimensional  algebras like algebras of smooth functions and
their generalizations so again (as  in last section) one has to
take some care of duality and tensor products.\\

Let $\cala$ be now an arbitrary associative algebra and let
$C(\cala)$  denote the graded vector space of multilinear forms
on $\cala$, i.e.  $C(\cala)=\oplusinf_n C^n(\cala)$ where 
$C^n(\cala)=(\otimes^n\cala)^\ast$ is the dual of the n-th
tensor power of $\cala$. One has  $T(\cala^\ast)\subset
C(\cala)$ and the equality  $T(\cala^\ast)=C(\cala)$ holds if
and only if $\cala$ is finite  dimensional. The product of
$T(\cala^\ast)$ (i.e. the tensor product)  canonically extends
to $C(\cala)$ which so equipped is a graded  algebra.
Furthermore minus the dual of the product of $\cala$ is a 
linear mapping of $C^1(\cala)=\cala^\ast$ into 
$C^2(\cala)=(\cala\otimes \cala)^\ast$ which also canonically
extends  as an antiderivation $d$ of $C(\cala)$ which is a
differential as  consequence of the associativity of the
product of $\cala$. It is  given by: \[
d\omega(x_0,\dots,x_n)=\sum^n_{k=1}(-1)^k\omega(x_0,\dots,x_{i-1}x_i,\dots,x_n)
\] for $\omega \in C^n(\cala)$ and $x_i\in \cala$. The graded
differential algebra $C(\cala)$ is the generalization of the
above $T(\cala^\ast)$ for an infinite dimensional algebra
$\cala$. As announced before the cohomology of $C(\cala)$ is
trivial whenever $\cala$ has a unit.

\begin{lemma} Let $\cala$ be a unital associative algebra (over
$\mathbb C$). Then the cohomology $H(C(\cala))$ of $C(\cala)$
is trivial in the sense that one has: \[ H^0(C(\cala))=\mathbb
C\ \mbox{and}\  H^n(C(\cala))=0\ \mbox{for}\  n\geq 1. \]
\end{lemma}

\noindent\underbar{Proof}. By definition $C(\cala)$ is
connected so $H^0(C(\cala))=\mathbb C$ is obvious. For
$\omega\in C^n(\cala)$ with $n\geq 1 $ let us define
$h(\omega)\in C^{n-1}(\cala)$ by
$h(\omega)(x_1,\dots,x_{n-1})=\omega(\bbbone,
x_1,\dots,x_{n-1})$, $\forall x_i\in\cala$. One has \[
d(h(\omega))+h(d(\omega))=\omega\ \mbox{for any}\  \omega\in
C^n(\cala)\ \mbox{with}\  n\geq 1 \] which implies
$H^n(C(\cala))=0$ for $n\geq 1$. $\square$\\

If $\calm$ is a $(\cala,\cala)$-bimodule, then the graded
vector space of $\calm$-{\sl valued Hochschild  cochains of}
$\cala$ is the graded vector space $C(\cala,\calm)$ of
multilinear mappings of $\cala$ into $\calm$, i.e.
$C^n(\cala,\calm)$ is the space of linear mappings of
$\otimes^n\cala$ into $\calm$, equipped with the Hochschild
coboundary $d_H$ defined by \[ \begin{array}{ll}
d_H(\omega)(x_0,\dots,x_n)=x_0\omega(x_1,\dots,x_n)&+d(\omega)(x_0,\dots,x_n)\\
&+(-1)^{n+1}\omega(x_0,\dots,x_{n-1})x_n \end{array} \] for
$\omega\in C^n(\cala,\calm)$, $x_i\in\cala$ and where $d$ is
``the obvious extension" to $C(\cala,\calm)$ of the 
differential $d$ of $C(\cala)$. When $\cala$ is finite
dimensional all this reduces to the previous definitions, in
particular in this case one has $C(\cala,\calm)=\calm\otimes
T(\cala^\ast)$. The cohomology $H(\cala,\calm)$ of
$C(\cala,\calm)$ is the $\calm$-{\sl valued Hochschild
cohomology of} $\cala$ or the {\sl Hochschild cohomology of
$\cala$ with coefficients in $\calm$.} The $\calm$-valued
Hochschild cochains of $\cala$ which vanishes whenever one of
their arguments is the unit $\bbbone$ of $\cala$ are said to be
{\sl normalized Hochschild cochains}. The graded vector space
$C_0(\cala,\calm)$ of $\calm$-valued normalized Hochschild
cochains is stable by the Hochschild coboundary $d_H$ and it is
well known and easy to show that the injection of
$C_0(\cala,\calm)$ into $C(\cala,\calm)$ induces an isomorphism
in cohomology, i.e. the cohomology of $C_0(\cala,\calm)$ is
again $H(\cala,\calm)$. Notice that a $\calm$-valued Hochschild
1-cocycle (i.e. an element of $C^1(\cala,\calm)$ in
$\ker(d_H)$) is a derivation $\delta$ of $\cala$ in $\calm$,
and that it is automatically normalized. If $\caln$ is another
$(\cala,\cala)$-bimodule then the tensor product over $\cala$
of $\calm$ and $\caln$,
$(\calm,\caln)\mapsto\calm\otimesinf_\cala\caln$, induces a
product $(\alpha,\beta)\mapsto \alpha\cup \beta$, {\sl the cup
product} $\cup:C(\cala,\calm)\otimes C(\cala,\caln)\rightarrow
C(\cala,\calm\otimesinf_\cala\caln)$ such that
$C^m(\cala,\calm)\cup C^n(\cala,\caln)\subset
C^{m+n}(\cala,\calm\otimesinf_\cala \caln)$ defined by \[
(\alpha\cup\beta)(x_1,\dots,x_{m+n})=\alpha(x_1,\dots,x_m)\otimesinf_\cala
\beta(x_{m+1},\dots,x_{m+n}) \] for $\alpha\in
C^m(\cala,\calm),\beta\in C^n(\cala,\caln)$ and $x_i\in\cala$.
If $\calp$ is another $(\cala,\cala)$-bimodule and if
$\gamma\in C^p(\cala,\calp)$, one
has:$(\alpha\cup\beta)\cup\gamma=\alpha\cup(\beta\cup\gamma)$.
Furthermore one has
$d_H(\alpha\cup\beta)=d_H(\alpha)\cup\beta+(-1)^m \alpha\cup
d(\beta)$ for $\alpha\in C^m(\cala,\calm)$, $\beta\in
C(\cala,\caln)$. This implies in particular that
$C(\cala,\cala)$ is a graded differential algebra (when
equipped with the cup product and with $d_H$). In fact,
$C(\cala,\cala)$ has a very rich structure which was first described in \cite{gerst}.  As pointed out in
\cite{gerst}, its cohomology $H(\cala,\cala)$ which inherits from this structure is graded commutative (as graded algebra for the
cup product). The cohomology $H(\cala,\cala)$ is a sort of
graded commutative Poisson algebra.\\

A unital associative algebra $\cala$ is said to be of {\sl
Hochschild dimension} $n$ if $n$ is the smaller integer such
that $H^k(\cala,\calm)=0$ for any $k\geq n+1$ and any
$(\cala,\cala)$-bimodule $\calm$. The Hochschild dimension of
the algebra $\mathbb C[X_1,\dots,X_n]$ of complex polynomials
with $n$ indeterminates is $n$. If one considers $\cala$ as the
generalization of the algebra of smooth functions on a
noncommutative space then its Hochschild dimension $n$ is the
analog of the dimension of the noncommutative space.\\

In spite of the triviality of the cohomology of $C(\cala)$,
several complexes with nontrivial cohomologies can be extracted
from it. Let $\cals:C(\cala)\rightarrow C(\cala)$ and
$\calc:C(\cala)\rightarrow C(\cala)$ be linear mappings defined
by \[
\cals(\omega)(x_1,\dots,x_n)=\sum_{\pi\in\cals_n}\varepsilon(\pi)\omega(x_{\pi(1)},\dots,x_{\pi(n)})
\] and \[
\calc(\omega)(x_1,\dots,x_n)=\sum_{\gamma\in\calc_n}\varepsilon(\gamma)\omega(x_{\gamma(1)},\dots,x_{\gamma(n)})
\] for $\omega\in C^n(\cala)$, $x_i\in\cala$ and where
$\cals_n$ is the group of permutations of $\{1,\dots,n\}$ and
$\calc_n$ is the subgroup of cyclic permutations,
$(\varepsilon(\pi)$ denoting the signature of the permutation
$\pi$). The mapping
$C(\cala)\stackrel{\cals}{\rightarrow}\cals(C(\cala))$ is a
homomorphism of graded differential algebras of $C(\cala)$ onto
the graded differential algebra $C_\wedge(\cala_\lied)$ of Lie
algebra cochains of the underlying Lie algebra $\cala_\lied$
with values in the trivial representation of $\cala_\lied$ in
$\mathbb C$; (Notice that the product of
$C_\wedge(\cala_\lied)$ is not induced by the inclusion
$C_\wedge(\cala_\lied)\subset C(\cala))$. The cohomology of
$\im(\cals)=C_\wedge(\cala_\lied)$ is therefore the Lie algebra
cohomology of $\cala_\lied$. On the other hand, (see Lemma 3 in
\cite{connes:01} part II), one has $\calc\circ d=d_H\circ
\calc$ where $d_H$ is the Hochschild coboundary of
$C(\cala,\cala^\ast)$ and therefore $(\im(\calc),d_H)$ is a
complex the cohomology of which coincides with {\sl the cyclic
cohomology $H_\lambda(\cala)$ of} $\cala$ up to a shift $-1$ in
degree \cite{connes:01}.\\

Let us define for $a\in\cala$ the homogeneous linear mapping
$i_a$ of degree $-1$ of $C(\cala)$ into itself by setting \[
i_a(\omega)(x_1,\dots,x_{n-1})=\sum^{n-1}_{k=0}(-1)^k\omega(x_1,\dots,x_k,a,x_{k+1},\dots,x_{n-1})
\] for $\omega\in C^n(\cala)$ with $n\geq 1$ and $x_i\in\cala$,
and by setting $i_a(C^0(\cala))=0$. For each $a\in\cala$, $i_a$
is an antiderivation of $C(\cala)$ and it is easy to verify
that $a\mapsto i_a$ is an operation of the Lie algebra
$\cala_\lied$ in the graded differential algebra $C(\cala)$.
The homotopy $h$ used in the proof of Lemma 3 commutes with the
$L_a$'s which implies that the invariant cohomology
$H_I(C(\cala))$ of $C(\cala)$ is also trivial. The basic cohomology
of $C(\cala)$ for this operation has been called {\sl basic
cohomology of} $\cala$ and denoted by $H_B(\cala)$ in
\cite{dv:tm}. It is given by the following theorem \cite{dv:tm}

\begin{theo} The basic cohomology $H_B(\cala)$ of $\cala$
identifies with the algebra $\cali_S(\cala_\lied)$ of invariant
polynomials on the Lie algebra $\cala_\lied$ where the degree
$2n$ is given to the homogeneous polynomials of degree $n$,
that is $H^{2n}_B(\cala)=\cali^n_S(\cala_\lied)$ and
$H^{2n+1}_B(\cala)=0.$ \end{theo}

The proof of this theorem which is not straightforward uses a
familiar trick in equivariant cohomology to convert the
operation $i$ of $\cala_\lied$ into a differential.\\

Two algebras $\cala$ and $\calb$ (associative unital, etc.) are
said to be {\sl Morita equivalent} if there is a
$(\cala,\calb)$-bimodule $\calu$ and a $(\calb,\cala)$-bimodule
$\calv$ such that one has an isomorphism of
$(\cala,\cala)$-bimodules $\calu\otimesinf_\calb \calv \simeq
\cala$ and an isomorphism of $(\calb,\calb)$-bimodules
$\calv\otimesinf_\cala \calu\simeq \calb$. This is an
equivalence relation and this induces an equivalence between
the category of right $\cala$-modules (resp. left
$\cala$-modules, $(\cala,\cala)$-bimodules) and the category of
right $\calb$-modules (resp. left $\calb$-modules,
$(\calb,\calb)$-bimodules). The algebras $M_m(\cala)$ and
$M_n(\cala)$ of $m\times m$ matrices and of $n\times n$
matrices with entries in $\cala$ are Morita equivalent for any
$m, n\in\mathbb N$; in fact the
$(M_m(\cala),M_n(\cala))$-bimodule  $M_{mn}(\cala)$ of
rectangular $m\times n$ matrices and the
$(M_n(\cala),M_m(\cala))$-bimodule $M_{nm}(\cala)$ of
rectangular $n\times m$ matrices with entries in $\cala$ are
such that
$M_n(\cala)=M_{nm}(\cala)\otimesinf_{M_m(\cala)}M_{mn}(\cala)$
and
$M_m(\cala)=M_{mn}(\cala)\otimesinf_{M_n(\cala)}M_{nm}(\cala)$,
(the tensor products over $M_m(\cala)$ and $M_n(\cala)$ being
canonically the usual matricial products).\\

An important property of Hochschild cohomology and cyclic
cohomology (and of the corresponding homologies) is their
Morita invariance \cite{jac}, \cite{jll}, \cite{weib}. More
precisely if $\cala$ and $\calb$ are Morita equivalent with
$\calu$ and $\calv$ as above and if $\calm$ is a
$(\cala,\cala)$-bimodule (resp. $\caln$ is a
$(\calb,\calb)$-bimodule) one has a canonical isomorphism
$H(\cala,\calm)\simeq H(\calb,\calv
\otimesinf_\cala\calm\otimesinf_\cala\calu)$, (resp.
$H(\calb,\caln)\simeq
H(\cala,\calu\otimesinf_\calb\caln\otimesinf_\calb \calv)$) in
 Hochschild cohomology and also
$H_\lambda(\cala)\simeq H_\lambda(\calb)$ in cyclic cohomology.
In contrast, the Lie algebra cohomology $H(\cala_\lied)$ and
the basic cohomology $H_B(\cala)$ are not Morita invariant
since for instance for $\cala=M_n(\mathbb C)$ they depend on
the number $n\in \mathbb N$ whereas $M_n(\mathbb C)$ is Morita
equivalent to $\mathbb C$.

\section{Categories of algebras}

In this section we consider general algebras over $\mathbb C$.
That  is by an {\sl algebra} we here mean a complex vector
space $\cala$ equipped with a  bilinear product $m:\cala\otimes
\cala\rightarrow\cala$. Given two such  algebras $\cala$ and
$\calb$, an {\sl algebra homomorphism} of $\cala$ into $\calb$
is  a linear mapping $\varphi:\cala\rightarrow \calb$ such that
$\varphi(m(x\otimes y))=m(\varphi(x)\otimes\varphi(y))$,
$(\forall x, y\in \cala)$, i.e. $\varphi\circ
m=m\circ(\varphi\otimes\varphi)$. \\

Let us define the category $\mathbf A$ to be the category such
that  the class $\mbox{Ob}(\mathbf A)$ of its objects is the
class of all algebras  (in the above sense) and such that for
any $\cala, \calb \in \mbox{Ob}(\mathbf A)$  the set
$\hom_{\mathbf A}(\cala,\calb)$ of morphisms from $\cala$ to
$\calb$  is the set of all algebra homomorphisms of $\cala$
into $\calb$.\\

A subcategory of $\mathbf A$ will be called a {\sl category of 
algebras}. Thus a category $\mathbf C$ is a category of
algebras if $\mbox{Ob}(\mathbf  C)$ is a subclass of
$\mbox{Ob}(\mathbf A)$ and if, for any $\cala, \calb\in
\mbox{Ob}(\mathbf C)$, one has $\hom_{\mathbf
C}(\cala,\calb)\subset \hom_{\mathbf A}(\cala,\calb)$. We now
list some categories of  algebras which will be used later.\\

\noindent 1. The category $\mathbf{Alg}$ of unital associative 
algebras: $\mbox{Ob}(\mathbf{Alg})$ is the class of all complex
unital  associative algebras and for any $\cala,\calb\in
\mbox{Ob}(\mathbf{Alg})$, $\mbox{Hom}_{\alg}(\cala,\calb)$ is
the set of all algebra  homomorphisms mapping the unit of
$\cala$ onto the unit of $\calb$.\\

\noindent 2. The category $\mathbf{Alg}_Z$ is the subcategory
of  $\mathbf{Alg}$ defined by
$\mbox{Ob}(\mathbf{Alg}_Z)=\mbox{Ob}(\mathbf{Alg})$ and  for
any $\cala,\calb\in \mbox{Ob}(\mathbf{Alg}_Z)$, 
$\hom_{\mathbf{Alg}_Z}(\cala,\calb)$ is the set of all
$\varphi\in \hom_{\mathbf{Alg}}(\cala,\calb)$  mapping the
center $Z(\cala)$ of $\cala$ into the center $Z(\calb)$ of
$\calb$,  i.e. such that $\varphi(Z(\cala))\subset Z(\calb)$.\\

\noindent 3. The category $\mathbf{Jord}$ of complex unital
Jordan  algebras: $\mbox{Ob}(\mathbf{Jord})$ is the class of
all complex unital Jordan  algebras and for any\linebreak[4]
$\cala,\calb\in\mbox{Ob}(\mathbf{Jord})$,
$\hom_{\mathbf{Jord}}(\cala,\calb)$ is the set of all algebra 
homomorphisms mapping the unit of $\cala$ onto the unit of
$\calb$.\\

\noindent 4. The category $\mathbf{Algcom}$ of unital
associative and commutative algebras: $\mbox{Ob}(\mathbf{Algcom})$
is the class of all complex  unital associative commutative
algebras and for any $\cala,\calb\in
\mbox{Ob}(\mathbf{Algcom})$,
$\hom_{\mathbf{Algcom}}(\cala,\calb)=\hom_{\mathbf{Alg}}(\cala,\calb)$.\\

\noindent 5. The category $\mathbf{Lie}$ of Lie algebras: 
$\mbox{Ob}(\mathbf {Lie})$ is the class of all complex Lie
algebras and for any $\cala,  \calb\in
\mbox{Ob}(\mathbf{Lie})$,
$\hom_{\mathbf{Lie}}(\cala,\calb)=\hom_{\mathbf{A}}(\cala,\calb)$.\\

\noindent \underbar{Remark 4}. If $\cala\in\ob(\mathbf{Alg})$
and  $\calb\in \ob(\mathbf{Algcom})$, one has\linebreak[4]
$\hom_{\mathbf{Alg}}(\cala,\calb)=\hom_{\mathbf{Alg}_Z}(\cala,\calb)$.

On the other hand if $\cala$ and $\calb$ are objects of 
$\mathbf{Algcom}$ then \[
\hom_{\mathbf{Algcom}}(\cala,\calb)=\hom_{\mathbf{Jord}}(\cala,\calb).
\]

Thus $\algcom$ is {\sl a full subcategory} of $\alg$, of
$\alg_Z$ and  of $\jord$, i.e. for any $\cala,\calb\in
\ob(\algcom)$ one has :\\ $
\hom_{\algcom}(\cala,\calb)=\hom_{\alg}(\cala,\calb)=\hom_{{\alg}_Z}(\cala,\calb)=\hom_{\jord}(\cala,\calb)$\\

In order to discuss reality conditions we shall also need
categories  of $\ast$-algebras. By a $\ast$-algebra we here
mean a general complex  algebra $\cala$ as above equipped with
an antilinear involution $x\mapsto  x^\ast$ such that
$m(x\otimes y)^\ast=m(y^\ast\otimes x^\ast)$, (i.e. such that
it  reverses the order in the product). If $\cala$ and $\calb$
are  $\ast$-algebras,  a $\ast$-algebra homomorphism of $\cala$
into $\calb$ is an algebra  homomorphism $\varphi$ of $\cala$
into $\calb$ which preserves the involutions,  i.e.
$\varphi(x^\ast)=\varphi(x)^\ast$ for $x\in \cala$. One defines
the  category of algebras $\ast$-${\mathbf A}$ to be the
category where  $\ob(\ast$-${\mathbf A})$ is the class of
 $\ast$-algebras and such that for any  $\cala,\calb\in
\ob(\ast\mbox{-}{\mathbf  A})$, 
$\hom_{\ast\mbox{-}{\mathbf{A}}}(\cala,\calb)$ is the set of
$\ast$-algebra homomorphisms of $\cala$ into $\calb$. A 
subcategory of $\ast$-${\mathbf A}$ will be called a {\sl
category of  $\ast$-algebras} and one defines in the obvious
manner the categories of $\ast$-algebras  $\ast$-$\alg$,
$\ast$-$\alg_Z$, $\ast$-$\jord$,\linebreak[4] $\ast$-$\algcom$,
$\ast$-$\lie$  corresponding to the above examples 1, 2, 3, 4,
5.\\

Let $\mathbf C$ be a category of algebras and let $\cala$ be
an  object of $\mathbf C$ with product denoted by $a\otimes
a'\mapsto aa'$  ($a,a'\in a$). A complex vector space $\cale$
will be said to be a $\cala$-{\sl  bimodule for $\mathbf C$} if
there are linear mappings  $\cala\otimes\cale\rightarrow \cale$
and $\cale\otimes \cala\rightarrow \cale$, denoted by
$a\otimes  e\mapsto ae$ and $e\otimes a\mapsto ea$ ($a\in
\cala,e\in \cale$) respectively, such  that the direct sum
$\cala\oplus \cale$ equipped with the product \[ (a\oplus 
e)\otimes (a'\oplus e')\mapsto aa'\oplus (ae'+ea') \] is an
object of $\mathbf  C$ and such that the canonical linear
mappings \[ i:\cala\rightarrow \cala\oplus  \cale\ \mbox{and}\ 
p:\cala\oplus \cale\rightarrow \cala \] defined by 
$i(a)=a\oplus 0$ and $p(a\oplus e)=a$  $(\forall a\in \cala$
and $\forall e\in \cale)$  are morphisms of $\mathbf C$. In
other words $\cale$ is a  $\cala$-bimodule for $\mathbf C$ if
$\cala\oplus \cale$ is equipped with a bilinear  product
vanishing on $\cale\otimes \cale$ and such that
$\cala\oplus\cale \in \ob  (\mathbf C)$~, $i\in \hom_{\mathbf
C}(\cala,\cala\oplus\cale)$ and $p\in  \hom_{\mathbf
C}(\cala\oplus\cale,\cala)$.\\

For the category $\ag$ this notion of bimodule is not very 
restrictive. In fact, if $\cala$ is an algebra (i.e.
$\cala\in\ob(\ag))$ then a  $\cala$-bimodule for $\ag$ is
simply a complex vector space $\cale$ with two bilinear 
mappings corresponding to linear mappings
$\cala\otimes\cale\rightarrow \cale$  and $\cale\otimes
\cala\rightarrow \cale$ as above. These two linear  mappings
will be always denoted by $a\otimes e\mapsto ae$ and $e\otimes
a\mapsto ea$  and called {\sl left} and {\sl right action of
$\cala$ on $\cale$}. Let us  describe what restrictions occur
for the categories of algebras of examples 1, 2,  3, 4, 5.\\

\noindent 1. Let $\cala$ be a unital associative complex
algebra with  product denoted by $a\otimes a'\mapsto aa'$ and
unit denoted by $\bbbone$.  Then, $\cale$ is a $\cala$-bimodule
for $\alg$ if and only if one has

\[ \begin{array}{lllll} (i)& & (aa')e=a(a'e) & \mbox{and} &
\bbbone  e=e\\ (ii) & & e(aa')=(ea)a' &  \mbox{and} &
e\bbbone=e\\ (iii) & & (ae)a'=a(ea')  & & \end{array} \] for
any $a,a'\in\cala$ and $e\in\cale$. Conditions  $(i)$ express
the fact that $\cale$ is a left $\cala$-module in the usual
sense,  conditions $(ii)$ express the fact that $\cale$ is a
right $\cala$-module in the  usual sense whereas, completed
with the compatibility condition $(iii)$, all  these conditions
express the fact that $\cale$ is a $(\cala,\cala)$-bimodule in
the  usual sense for unital associative algebras.\\

\noindent 2. Let $\cala$ be as in 1 above. Then $\cale$ is a 
$\cala$-bimodule for $\alg_Z$ if and only if it is a
$\cala$-bimodule for $\alg$ such that  one has \[ ze = ez \]
for any element $z$ of the center $Z(\cala)$ of $\cala$  and
$e\in \cale$. This condition expresses that as 
$(Z(\cala),Z(\cala))$-bimodule, $\cale$ is the underlying
bimodule of a $Z(\cala)$-module. Such 
$(\cala,\cala)$-bimodules were called {\sl central bimodules}
over $\cala$ in \cite{mdv:pm1}, \cite{mdv:pm2} (see also in
\cite{dv:4}). We shall keep this  terminology here and call
central bimodule a bimodule for $\alg_Z$.\\

Let $\cale$ be a $\cala$-bimodule for $\alg$ (i.e. a
$(\cala,\cala)$-bimodule). One can associate to $\cale$ two
$\cala$-bimodules for $\alg_Z$ (i.e. two central bimodules)
$\cale^Z$ and $\cale_Z$. The bimodule $\cale^Z$ is the biggest
$(\cala,\cala)$-subbimodule of $\cale$ which is central and we
denote by $i^Z$ the canonical inclusion of $\cale^Z$ into
$\cale$ whereas $\cale_Z$ is the quotient of $\cale$ by the
$(\cala,\cala)$-subbimodule $[Z(\cala),\cale]$ generated by the
$ze-ez$ where $z$ is in the center $Z(\cala)$ of $\cala$,
$e\in\cale$ and we denote by $p_Z$ the canonical projection of
$\cale$ onto $\cale_Z$. The pair $(\cale^Z,i^Z)$ is
characterized by the following universal property: {\sl For any
$(\cala,\cala)$-bimodule homomorphism $\Phi:\caln\rightarrow
\cale$ of a central bimodule $\caln$ into $\cale$, there is a
unique $(\cala,\cala)$-bimodule homomorphism
$\Phi^Z:\caln\rightarrow \cale^Z$ such that $\Phi=i^Z\circ
\Phi^Z$.} The pair $(\cale_Z,p_Z)$ is characterized by the
following universal property: {\sl For any
$(\cala,\cala)$-bimodule homomorphism $\varphi:\cale\rightarrow
\calm$ of $\cale$ into a central bimodule $\calm$ there is a
unique $(\cala,\cala)$-bimodule homomorphism
$\varphi_Z:\cale_Z\rightarrow \calm$ such that
$\varphi=\varphi_Z\circ p_Z$.} In functorial language, this
means that $\cale\mapsto\cale^Z$ is a right adjoint and that
$\cale\mapsto \cale_Z$ is a left adjoint of the canonical
functor $I_Z$ from the category of $\cala$-bimodules for
$\alg_Z$ in the category of $\cala$-bimodules for $\alg$.
Notice also that $\cale$ is central if and only if
$\cale=\cale_Z$ which is equivalent to $\cale=\cale^Z$ and that
if $\calm$ and $\caln$ are two $\cala$-bimodules for $\alg_Z$
(i.e. two central bimodules) then one has
$(\calm\otimes\caln)_Z=\calm\otimesinf_{Z(\cala)}\caln$. One
has the further following stability properties for the
$\cala$-bimodules for $\alg_Z$: Every subbimodule of a central
bimodule is central, every quotient of a central bimodule is
central and any product of central bimodules is central. For
all this, we refer to \cite{mdv:pm2}.\\

\noindent 3. Let $\calj$ be a complex unital Jordan algebra
with  product denoted by $x\otimes y\mapsto x\bullet y$
($x,y\in \calj$) and unit  $\bbbone$. Then $\cale$ is a
$\calj$-bimodule for $\jord$ if and only if one has \[
\begin{array}{lll} (i)& & xe = ex\hspace{0,3cm}  
\mbox{and}\hspace{0,3cm} \bbbone e=e\\ (ii)& & x((x\bullet
x)e)=(x\bullet x)(xe)\\ (iii)& & ((x\bullet  x)\bullet
y)e-(x\bullet x)(ye)=2((x\bullet y)(xe)-x(y(xe))) \end{array}
\] for  any $x,y\in\calj$ and $e\in\cale$. Such a bimodule for
$\jord$ is called  a {\sl Jordan module} over $\calj$ \cite{nj}
which is natural since, in view  of $(i)$, there is only one
bilinear mapping of $\calj\times \cale$ into  $\cale$.\\

\noindent 4. Let $\calc$ be a unital associative commutative
complex  algebra. Then $\cale$ is a $\calc$-bimodule for
$\algcom$ if and only if it is  a $\calc$-bimodule for $\alg$
such that one has \[ ce = ec \] for any  $c\in \calc$ and $e\in
\cale$. This means that a $\calc$-bimodule for $\algcom$ is 
the same thing as (the underlying bimodule of) a $\calc$-module
in the usual  sense. Since the center of $\calc$ coincides with
$\calc$, $Z(\calc)=\calc$, this  implies that it is also the
same thing as a $\calc$-bimodule for $\alg_Z$, as  announced in
the introduction. Notice that in the case of a $\calc$-bimodule
for  $\alg$ one generally has $ce\not=ec$.\\

\noindent 5. Let $\fracg$ be a complex Lie algebra with product
(Lie  bracket) denoted by $X\otimes Y\mapsto [X,Y]$ for $X,Y\in
\fracg$. Then,  $\cale$ is a $\fracg$-bimodule for $\lie$ if
and only if one has \[  \begin{array}{lll} (i) & & Xe=-eX\\
(ii) && [X,Y]e=X(Ye)-Y(Xe) \end{array} \] for any $X,Y\in 
\fracg$ and $e\in \cale$. Condition $(i)$ shows that again
there is only one  bilinear mapping of $\fracg\times\cale$ into
$\cale$ and $(ii)$ means that $\cale$ is  the space of a linear
representation of $\fracg$; Thus a $\fracg$-bimodule for 
$\lie$ is what is usually called a $\fracg$-module (or a linear
representation of  $\fracg$).\\

One defines in a similar way the notion of $\ast$-bimodule for a
category  $\ast$-$\cg$ of $\ast\mbox{-}$algebras. Namely, if
$\cala\in \ob(\ast$-$\cg)$, a  complex vector space $\cale$
will be said to be a $\cala$-$\ast$-{\sl bimodule for} 
$\ast$-$\cg$ if $\cala\oplus\cale$ is equipped with a structure
of $\ast$-algebra  with product vanishing on
$\cale\otimes\cale$ such that $\cala\oplus\cale\in 
\ob(\ast$-$\cg)$,
$i\in\hom_{\ast\mbox{-}\cg}(\cala,\cala\oplus\cale)$ and
$p\in\hom_{\ast\mbox{-}\cg}(\cala\oplus\cale,\cala)$.\\

One can easily describe what is a $\ast$-bimodule for the
various  categories of $\ast$-algebras. If $\cala$ is a
$\ast$-algebra, we also denote by  $\cala$ the algebra obtained
by ``forgetting the involution". If $\cala$ is an  object of
$\ast$-$\alg$ then a $\cala$-$\ast$-bimodule for $\ast$-$\alg$
is a $\cala$-bimodule $\cale$ for $\alg$ which is equipped with
an  antilinear involution $e\mapsto e^\ast$ such that
$(xey)^\ast=y^\ast e^\ast  x^\ast$ for $x,y\in \cala$ and $e\in
\cale$, i.e. it is what has been called in  the introduction a
$\ast$-bimodule over the (unital associative complex)
$\ast$-algebra $\cala$. A $\cala$-$\ast$-bimodule
for $\ast$-$\alg_Z$ is then just such a $\ast$-bimodule over
$\cala$  which is central. If $\calc$ is a unital associative
complex commutative  $\ast$-algebra, then a
$\calc$-$\ast$-bimodule for $\ast$-$\algcom$ is just what has 
been called a $\ast$-module over the (unital associative
complex) commutative  $\ast$-algebra $\calc$.\\

One can proceed similarily with real algebras. However to be
in  conformity with the point of view of the introduction
concerning reality, we shall  work with $\ast$-algebras and,
eventually, extract their hermitian parts as  well as the
hermitian parts of the $\ast$-bimodules over them.

\section{First order differential calculi}

Throughout the following $\cala$ denotes a unital associative
complex  algebra. A pair $(\Omega^1,d)$ where $\Omega^1$ is a
$(\cala,\cala)$-bimodule  (i.e. a $\cala$-bimodule for $\alg$)
and where $d:\cala\rightarrow \Omega^1$  is a derivation of
$\cala$ into $\Omega^1$, that is a linear mapping which 
satisfies (the Leibniz rule) \[ d(xy)=d(x)y+xd(y) \] for any
$x,y\in \cala$,  will be called a {\sl first order differential
calculus over $\cala$ for $\alg$} or  simply a {\sl first order
differential calculus over $\cala$} \cite{slw}. If  furthermore
$\Omega^1$ is a central bimodule (i.e. a $\cala$-bimodule for 
$\alg_Z$), we shall say that $(\Omega^1,d)$ is a {\sl first
order differential calculus  over $\cala$ for $\alg_Z$}. One
can more generally define the notion of first order
differential calculus over $\cala$ for any category $\cg$ of
algebras  such that $\cala\in\ob(\cg)$.\\

\noindent \underbar{Remark 5}. If $\Omega^1$ is a
$\cala$-bimodule  for $\cg$ a derivation $d:\cala\rightarrow
\Omega^1$ can be defined to be a  linear mapping such that
$a\mapsto a\oplus d(a)$ is in 
$\hom_{\cg}(\cala,\cala\oplus\Omega^1)$. However, for the
category $\alg_Z$ this does not impose restrictions  on first
order differential calculus. Indeed if $\Omega^1$ is a central 
bimodule and if $d:\cala\rightarrow \Omega^1$ is a derivation
one has $d(z)a+zd(a)=d(za)=d(az)=ad(z)+d(a)z$ for any
$a\in\cala$ and $z$ in  the center $Z(\cala)$ of $\cala$, i.e.
$d(z)a=ad(z)$ since, by ``centrality",  $zd(a)=d(a)z$; again,
by centrality $z\omega=\omega z$, $\forall z\in Z (\cala)$ and 
$\forall \omega \in \Omega^1$, which finally implies $(z\oplus
d(z))(a\oplus\omega)=(a\oplus\omega)(z\oplus d(z))$ and
therefore  $z\oplus d(z)\in Z (\cala\oplus \Omega^1)$ for any
$z\in Z(\cala)$ which means  that the linear mapping $a\mapsto
a\oplus d(a)$ is in
$\hom_{\alg_Z}(\cala,\cala\oplus\Omega^1)$.\\

We shall refer to $d$ as the {\sl first order differential};
by  definition it is a $\Omega^1$-valued Hochschild cocycle of
degree 1 of $\cala$, i.e.  $d\in Z^1_H(\cala,\Omega^1)$.
Examples of first order differentials are  thus provided by
Hochschild coboundaries i.e. given by $d(x)=\tau x-x\tau$
($\forall  x\in \cala$) for some $\tau\in\Omega^1$. We shall
now explain that there are  ``universal first order
differential calculi" for $\alg$ and for $\alg_Z$ which define 
respectively functors from $\alg$ and from $\alg_Z$ in the
corresponding  categories of first order differential calculi.
For the case of a commutative algebra, there is also a
well-known universal first order differential calculus for
$\algcom$ which is the universal derivation into the module of
K\"ahler differentials (\cite{bour}, \cite{jll}, \cite{ps}). We
shall see however that it reduces to the universal calculus for
$\alg_Z$ (Corollary 1).\\

Let $m$ be the product of $\cala,(x,y)\mapsto m(x\otimes y)=xy$
and  let $\Omega^1_u(\cala)$ be the kernel of $m$, i.e. one has
the short  exact sequence \[ 0\rightarrow
\Omega^1_u(\cala)\stackrel{\subset}{\rightarrow}
\cala\otimes\cala\stackrel{m}{\rightarrow}\cala\rightarrow 0 \]
of $(\cala,\cala)$-bimodules ($\cala$-bimodules for $\alg$).
Define $d_u:\cala\rightarrow \Omega^1_u(\cala)$ by
$d_u(x)=\bbbone \otimes  x-x\otimes \bbbone$, $\forall x\in
\cala$. One verifies easily that $d_u$ is a  derivation. The
first order differential calculus $(\Omega^1_u(\cala),d_u)$
over  $\cala$ is characterized uniquely (up to an isomorphism)
by the following  universal property \cite{cart:02},
\cite{bour}.\\

\begin{proposition} For any first order differential calculus 
$(\Omega^1,d)$ over $\cala$, there is a unique bimodule
homomorphism $i_d$ of  $\Omega^1_u(\cala)$ into $\Omega^1$ such
that $d=i_d\circ d_u$. \end{proposition}

\noindent \underbar{Proof}. $\Omega^1_u(\cala)$ is generated
by  $d_u(\cala)$ as left module since $x^\alpha\otimes
y_\alpha$ with $x^\alpha  y_\alpha=0$ is the same thing as
$x^\alpha d(y_\alpha)$. On the other hand
$d_u(\bbbone)=0(=d_u(\bbbone^2)=2d_u(\bbbone))$. Therefore one
has a  surjective left $\cala$-module homomorphism of
$\cala\otimes (\cala/\mathbb  C\bbbone)$ onto
$\Omega^1_u(\cala)$, $x\otimes \dot y\mapsto xd_u(y)$, which
is  easily shown to be an isomorphism. Then $xd_u(y)\mapsto
xd(y)$ defines a left  $\cala$-module homomorphism $i_d$ of
$\Omega^1_u(\cala)$ into $\Omega^1$ which is  easily shown to
be a bimodule homomorphism by using the Leibniz rule for $d_u$
and  for $d$. One clearly has $d=i_d\circ d_u$. Uniqueness is
straightforward.  $\square$\\

Concerning the image of $i_d$, let us notice the following easy
lemma.

\begin{lemma} Let $(\Omega^1,d)$ be a first order differential
calculus over $\cala$. The following conditions are
equivalent.\\ \begin{tabular}{ll} $(i)$ & $\Omega^1$ is
generated by $d\cala$ as left $\cala$-module.\\ $(ii)$ &
$\Omega^1$ is generated by $d\cala$ as right $\cala$-module.\\
$(iii)$ & $\Omega^1$ is generated by $d\cala$ as 
$(\cala,\cala)$-bimodule.\\ $(iiii)$ & The homomorphism $i_d$
is surjective, i.e. $\Omega^1=i_d(\Omega^1_u(\cala))$.
\end{tabular} \end{lemma}

\noindent \underbar{Proof}. The equivalences
$(i)\Leftrightarrow(ii)\Leftrightarrow(iii)$ follows from
(Leibniz rule) \[ ud(v)w=ud(vw)-uvd(w)=d(uv)w-d(u)vw \] for
$u,v,w\in \cala$ whereas the equivalence $(iii)\Leftrightarrow
(iiii)$ is straightforward from the definitions. $\square$\\

\noindent \underbar{Remark 6}. Proposition 1 claims that there
is a  unique bimodule homomorphism $i_d$ of $\Omega^1_u(\cala)$
into $\Omega^1$  mapping the $\Omega^1_u(\cala)$-valued
Hochschild 1-cocycle $d_u$ on the  $\Omega^1$-valued Hochschild
1-cocycle $d$. One can complete the statement by the 
following: {\sl The $\Omega^1$-valued Hochschild 1-cocycle $d$
is a Hochschild  coboundary, (i.e. there is a $\tau\in\Omega^1$
such that $d(a)=\tau a-a\tau$ for any  $a\in\cala$), if and
only if $i_d$ has an extension $\tilde \imath_d$ as a bimodule 
homomorphism of $\cala\otimes\cala$ into $\Omega^1$},
\cite{ckmv}. In fact $\tau$ is  then $\tilde
\imath_d(\bbbone\otimes\bbbone)$, which gives essentially the
proof.\\

The first order differential calculus $(\Omega^1_u(\cala),d_u)$
{\sl  is universal for $\alg$}, it is usually simply called
{\sl the universal first  order differential calculus over
$\cala$}. From Proposition 1 follows the  functorial property.

\begin{proposition} Let $\cala$ and $\calb$ be algebras and let
$\varphi:\cala\rightarrow \calb$ be a homomorphism, (i.e. let 
$\cala,\calb\in \obth(\alg)$ and let
$\varphi\in\homth_{\alg}(\cala,\calb))$, then there  is a
unique linear mapping $\Omega^1_u(\varphi)$ of
$\Omega^1_u(\cala)$ into $\Omega^1_u(\calb)$ satisfying
$\Omega^1_u(\varphi)(x\omega
y)=\varphi(x)\Omega^1_u(\varphi)(\omega) \varphi(y)$ for any
$x,y\in  \cala$ and $\omega\in\Omega^1_u(\cala)$ and such that
$d_u\circ\varphi=\Omega^1_u(\varphi)\circ d_u$.
\end{proposition}

\noindent \underbar{Proof}. One equips $\Omega^1_u(\calb)$ of
a  structure of $(\cala,\cala)$-bimodule by setting $x\lambda 
y=\varphi(x)\lambda\varphi(y)$ for $x,y\in\cala$ and $\lambda
\in\Omega^1_u(\calb)$. Then  $d=d_u\circ\varphi$ is a
derivation of $\cala$ into the $(\cala,\cala)$-bimodule 
$\Omega^1_u(\calb)$, i.e. $(\Omega^1_u(\calb),d)$ is a first
order differential calculus over  $\cala$, and the result
follows from Proposition 1 with $\Omega^1_u(\varphi)=i_d$. 
$\square$\\

One can summarize the content of Proposition 2 by the
following: {\sl  For any homomorphism
$\varphi:\cala\rightarrow\calb$ (of unital associative 
$\mathbb C$-algebras) there is a unique
$(\cala,\cala)$-bimodule homomorphism
$\Omega^1_u(\varphi):\Omega^1_u(\cala)\rightarrow
\Omega^1_u(\calb)$  for which the diagram

\[ \begin{diagram} \node{\cala} \arrow{e,t}{\varphi}
\arrow{s,r}{d_u} \node{\calb}\arrow{s,r}{d_u}\\
\node{\Omega^1_u(\cala)} \arrow{e,t}{\Omega^1_u(\varphi)}
\node{\Omega^1_u(\calb)}  \end{diagram} \] is commutative}. All
this was for the category $\alg$, we now pass to  $\alg_Z$.\\

Let $[Z(\cala),\Omega^1_u(\cala)]$ be the subbimodule of 
$\Omega^1_u(\cala)$ defined by \[
[Z(\cala),\Omega^1_u(\cala)]=\{z\omega-\omega z\vert  z\in
Z(\cala),\omega\in\Omega^1_u(\cala)\}. \] By definition the
quotient
$\Omega^1_Z(\cala)=\Omega^1_u(\cala)/[Z(\cala),\Omega^1_u(\cala)]$
is  a central bimodule i.e. a $\cala$-bimodule for $\alg_Z$.
Let $p_Z:\Omega^1_u(\cala)\rightarrow \Omega^1_Z(\cala)$ be the
canonical  projection and let $d_Z:\cala\rightarrow
\Omega^1_Z(\cala)$ be defined by  $d_Z=p_Z\circ d_u$. Then
$d_Z$ is again a derivation so $(\Omega^1_Z(\cala),d_Z)$  is a
first order differential calculus over $\cala$ for $\alg_Z$. It
is  characterized uniquely (up to an isomorphism) among the
first order differential  calculi over $\cala$ for $\alg_Z$ by
the following universal property  \cite{mdv:pm2}.

\begin{proposition} For any first order differential calculus 
$(\Omega^1,d)$ over $\cala$ for $\alg_Z$, there is a unique
bimodule homomorphism $i_d$ of $\Omega^1_Z(\cala)$ into
$\Omega^1$ such that $d=i_d\circ d_Z$; i.e.  there is a unique
morphism of first order differential calculi over $\cala$ for 
$\alg_Z$ from $(\Omega^1_Z(\cala),d_Z)$ to $(\Omega^1,d)$.
\end{proposition}

\noindent\underbar{Proof}. The unique bimodule homomorphism
$i_d:\Omega^1_u(\cala)\rightarrow \Omega^1$ of Proposition 1
vanishes  on $[Z(\cala),\Omega^1_u(\cala)]$ since $\Omega^1$ is
central. Therefore  it factorizes as
$\Omega^1_u(\cala)\stackrel{p_Z}{\rightarrow}\Omega^1_Z(\cala)\rightarrow\Omega^1$
through a unique bimodule homomorphism, again denoted  $i_d$,
of $\Omega^1_Z(\cala)$ into $\Omega^1$ for which one has
$d=i_d\circ d_Z$. Again, uniqueness is obvious. $\square$\\

\noindent\underbar{Remark 7}. Proposition 3 can be slightly 
improved. One can replace the assumption ``$(\Omega^1,d)$ over
$\cala$ for $\alg_Z$" by ``$(\Omega^1,d)$ over $\cala$ such
that $zd(a)=d(a)z$ for any  $a\in\cala$ and $z\in Z(\cala)$" in
the statement. That is, what is important is that  the
subbimodule of $\Omega^1$ generated by $d\cala$ is central.\\

The first order differential calculus $(\Omega^1_Z(\cala),d_Z)$
will  be called {\sl the universal first order differential
calculus over $\cala$ for  $\alg_Z$}. Concerning the functorial
property of this first order differential  calculus,
Proposition 2 has the following counterpart for $\alg_Z$.

\begin{proposition} Let $\cala$ and $\calb$ be algebras as
above and let\linebreak[4] $\varphi:\cala\rightarrow\calb$ be a
homomorphism  such that $\varphi(Z(\cala))\subset Z(\calb)$,
(i.e. let  $\cala,\calb\in\obth(\alg_Z)$ and let
$\varphi\in\homth_{\alg_Z}(\cala,\calb))$, then there is a
unique  linear mapping $\Omega^1_Z(\varphi)$ of
$\Omega^1_Z(\cala)$ into  $\Omega^1_Z(\calb)$ satisfying
$\Omega^1_Z(\varphi)(x\omega
y)=\varphi(x)\Omega^1_Z(\varphi)(\omega)\varphi(y)$ for any 
$x,y\in\cala$ and $\omega\in\Omega^1_Z(\cala)$ and such that
$d_Z\circ \varphi=\Omega^1_Z(\varphi)\circ d_Z$.
\end{proposition}

\noindent \underbar{Proof}. Again, as in the proof of
Proposition 2, $\Omega^1_Z(\calb)$ is a
$(\cala,\cala)$-bimodule by setting $x\lambda
y=\varphi(x)\lambda\varphi(y)$ for $x,y\in\cala$ and
$\lambda\in\Omega^1_Z(\calb)$. Thus Proposition 4 follows from 
Proposition 3 if one can show that this bimodule is central
i.e. if $\varphi(z)\lambda=\lambda\varphi(z)$ for any $z\in
Z(\cala)$ and $\lambda\in\Omega^1_Z(\calb)$. This however
follows from the fact that $\Omega^1_Z(\calb)$ is central over
$\calb$ and that $\varphi$ maps  the center $Z(\cala)$ of
$\cala$ into the center $Z(\calb)$ of $\calb$.  $\square$\\

Again this can be summarized (by identifying
$\Omega^1_Z(\calb)$ with  a central bimodule over $\cala$ via
$\varphi$) as : {\sl For any $\varphi\in
\hom_{\alg_Z}(\cala,\calb)$, there is a unique homomorphism of 
$\cala$-bimodules for $\alg_Z$,
$\Omega^1_Z(\varphi)~:~\Omega^1_Z(\cala)~\rightarrow~
\Omega^1_Z(\calb)$,  for which the diagram

\[ \begin{diagram} \node{\cala} \arrow{e,t}{\varphi}
\arrow{s,r}{d_Z} \node{\calb}\arrow{s,r}{d_Z}\\
\node{\Omega^1_Z(\cala)} \arrow{e,t}{\Omega^1_Z(\varphi)}
\node{\Omega^1_Z(\calb)}  \end{diagram} \] is commutative}.\\

Proposition 3 has the following corollary

\begin{corol} If $\cala$ is commutative, $\Omega^1_Z(\cala)$ 
identifies canonically with the module of K\"ahler
differentials  $\Omega^1_{\cala\vert \mathbb C}$ and $d_Z$
identifies with the corresponding universal  derivation.
\end{corol}

\noindent\underbar{Proof}. The proof is straightforward since,
for a  commutative algebra $\cala$, a central bimodule is just
(the underlying bimodule  of) a $\cala$-module and then,
Proposition 3 just reduces to the universal  property which
characterizes  the first order K\"ahler differential calculus 
(see e.g. in \cite{bour}, \cite{jll}, \cite{ps}). $\square$\\

\noindent \underbar{Remark 8}. If $\cala$ is commutative, the
module of K\"ahler differentials $\Omega^1_{\cala\vert\mathbb
C}$ is known to be a version of differential 1-forms. There is
however a little subtility. In fact
$\Omega^1_{\cala\vert\mathbb C}$ is the quotient of
$\Omega^1_u(\cala)$ which is a commutative algebra (a
subalgebra of $\cala\otimes\cala$) by the ideal
$(\Omega^1_u(\cala))^2$. If $\cala$ is the algebra of smooth
functions $C^\infty(M)$ on a manifold $M$, this means that
$\Omega^1_{\cala\vert\mathbb C}$ is the algebra of functions in
$\cala\otimes\cala=C^\infty(M)\otimes C^\infty (M)$ vanishing
on the diagonal of $M\times M$ modulo functions vanishing to
order one on the diagonal of\linebreak[4] $M\times M$. On the other hand it
is clear (by using the Taylor expansion around the diagonal)
that the ordinary differential 1-forms are smooth
functions of $C^\infty(M\times M)$ vanishing on the diagonal of
$M\times M$ modulo the functions vanishing to order one on the
diagonal of $M\times M$. The subtility here lies in the fact
that without completion of the tensor product, the inclusion
$C^\infty(M)\otimes C^\infty(M)\subset C^\infty(M\times M)$ is
a strict one so there is generally a slight difference between
$\Omega^1_{C^\infty(M)\vert\mathbb C}$ and the module
$\Omega^1(M)$ of smooth 1-forms on $M$. Apart from this, one
can consider that $(\Omega^1_Z(\cala),d_Z)$ generalizes the
ordinary first order differential calculus. This is in contrast
to what happens for $(\Omega^1_u(\cala),d_u)$. Indeed if 
$\cala$ is an algebra of functions on a set $S$ containing more
than one element, $(\card(S)>1)$, then $\Omega^1_u(\cala)$
consists of functions on  $S\times S$ which vanish on the
diagonal and is therefore {\sl not} the  underlying bimodule of
a module (nonlocality) whereas $(d_uf)(x,y)=f(y)-f(x)$
$(x,y\in  S)$ is the finite difference.

\section{Higher order differential calculi}

Let $\cala$ be as before a unital associative complex algebra.
A  $\mathbb N$-graded differential algebra
$\Omega=\oplusinf_{n\geq 0} \Omega^n$  such that the subalgebra
$\Omega^0$ of its elements of degree 0 coincides with  $\cala$,
$\Omega^0=\cala$, will be called a {\sl differential calculus
over  $\cala$ for $\alg$} or simply a {\sl differential
calculus over $\cala$}. If  furthermore the $\Omega^n$
($n\in\mathbb N$) are central bimodules over $\cala$, (i.e.
$\cala$-bimodules for $\alg_Z$), $\Omega$ will be said to be a
{\sl  differential calculus over $\cala$ for $\alg_Z$}.\\

Let us define the $(\cala,\cala)$-bimodules $\Omega^n_u(\cala)$
for  $n\geq 0$ by $\Omega^0_u(\cala)=\cala$ and by
$\Omega^n_u(\cala)=\Omega^1_u(\cala)\otimesinf_{\cala}\dots
\otimesinf_{\cala}\Omega^1_u(\cala)$ ($n$ factors) for $n\geq
1$. The  direct sum $\Omega_u(\cala)=\oplusinf_{n\geq
0}\Omega^n_u(\cala)$ is canonically  a graded algebra, it is
{\sl the tensor algebra over $\cala$,
$T_\cala(\Omega^1_u(\cala))$,  of the $(\cala,\cala)$-bimodule
$\Omega^1_u(\cala)$}. The derivation $d_u:\cala\rightarrow 
\Omega^1_u(\cala)$ has a unique extension, again denoted by
$d_u$, as a differential of $\Omega_u(\cala)$: In fact, it is
known on $\cala=\Omega^0_u(\cala)$  and $d^2_u=0$ fixes it on
$d_u(\cala)$ to be 0 so it is known on the  generators of
$\Omega_u(\cala)$ and the extension by the antiderivation
property to  the whole $\Omega_u(\cala)$ is well defined and
unique; moreover, $d^2_u$ is a  derivation vanishing on the
generators and therefore $d^2_u=0$. So equipped,
$\Omega_u(\cala)$ is a graded differential algebra \cite{kar}
which is characterized uniquely (up to an isomorphism) by the
following  universal property.

\begin{proposition} Any homomorphism $\varphi$ of unital
algebras of  $\cala$ into the subalgebra $\Omega^0$  of
elements of degree 0 of a graded  differential algebra
$\Omega$  has a unique extension
$\tilde\varphi:\Omega_u(\cala)\rightarrow \Omega$ as a
homomorphism  of graded differential algebras.
\end{proposition}

\noindent\underbar{Proof}. The $(\Omega^0,\Omega^0)$-bimodule 
$\Omega^1$ can be considered as a $(\cala,\cala)$-bimodule by
setting $x\lambda  y=\varphi(x)\lambda \varphi(y)$ for
$x,y\in\cala$ and $\lambda\in\Omega^1$ and then  $d\circ
\varphi$ defines a derivation of $\cala$ into $\Omega^1$.
Therefore, by  Proposition 1, there is a unique bimodule
homomorphism  $\varphi^1:\Omega^1_u(\cala)\rightarrow \Omega^1$
such that $d\circ \varphi=\varphi^1\circ  d_u:\cala\rightarrow
\Omega^1$ (namely $\varphi^1=i_{d\circ\varphi}$). The  property
of  $\Omega_u(\cala)$ to be the tensor algebra
$T_\cala(\Omega^1_u(\cala))$ implies that  $\varphi$ and
$\varphi^1$ uniquely extend as a homomorphism
$\tilde\varphi:\Omega_u(\cala)\rightarrow \Omega$ of graded
algebras.  By construction one has $\tilde\varphi\circ
d_u=d\circ\tilde\varphi$ on  $\cala$ and on $d_u\cala$ where it
vanishes which implies $\tilde\varphi\circ
d_u=d\circ\tilde\varphi$ everywhere by the antiderivation
property of  $d_u$ and $d$. $\square$\\

The graded differential algebra $\Omega_u(\cala)$ is called (in
view  of the above universal property) {\sl the universal
differential calculus over  $\cala$} (it is universal for
$\alg$). The functorial property follows immediately:  {\sl For
any homomorphism $\varphi:\cala\rightarrow \calb$ (of unital
$\mathbb  C$-algebras), there is a unique homomorphism 
$\Omega_u(\varphi):\Omega_u(\cala)\rightarrow \Omega_u(\calb)$
of graded differential algebra which extends  $\varphi$ (i.e.
$\varphi=\Omega_u(\varphi)\restriction \cala$)}. This defines
the  covariant functor $\Omega_u$ from the category $\alg$ to
the category $\dif$ of  graded differential algebras (the
morphisms being the homomorphisms of graded differential
algebras preserving the units).\\

Proposition 5 is clearly a generalization of Proposition 1.
There is  another useful generalization of the universality of
the Hochschild  $1\mbox{-}$cocycle $a\mapsto d_u(a)$ (which is
the content of Proposition 1) and of  Remark 6 which is
described in \cite{ckmv} (see also in \cite{cq}) and which we
now  review (Proposition 6 below). First, notice that
$(a_1,\dots,a_n)\mapsto  d_u(a_1)\dots d_u(a_n)$ is a
$\Omega^n_u(\cala)$-valued Hochschild $n$-cocycle  which is
normalized (i.e. which vanishes whenever one of the $a_i$ is
the unit  $\bbbone$ of $\cala$). Second, notice that the short
exact sequence of Section  6 (before Proposition 1) has the
following generalization for $n\geq 1$  \[  0\rightarrow
\Omega^n_u(\cala)\stackrel{\subset}{\rightarrow}\cala\otimes\Omega^{n-1}_u(\cala)\stackrel{m}{\rightarrow}\Omega^{n-1}_u(\cala)\rightarrow0 
\]  as short exact sequence of $(\cala,\cala)$-bimodules, where
$m$  is the left multiplication by elements of $\cala$ of
elements of  $\Omega^{n-1}_u(\cala)$, (the inclusion is
canonical). One has the following \cite{ckmv}.

\begin{proposition} Let $\calm$ be a $(\cala,\cala)$-bimodule
and let $(a_1,\dots,a_n)\mapsto c(a_1,\dots,a_n)$ be a
normalized  $\calm$-valued Hochschild $n$-cocyle. Then, there
is a unique bimodule homomorphism
$i_c:\Omega^n_u(\cala)\rightarrow \calm$ such that \[
c(a_1,\dots,a_n)=i_c(d_u(a_1)\dots d_u(a_n)),\  \forall a_i\in
\cala.  \] Furthermore, $c$ is a Hochschild coboundary if and
only if $i_c$ has  an extension $\tilde \imath_c$ as a bimodule
homomorphism of  $\cala\otimes\Omega^{n-1}_u(\cala)$ into
$\calm$. \end{proposition}

\noindent\underbar{Proof}. We only give here some indications
and  refer to \cite{ckmv} for the detailed proof. The proof of
the first part  proceeds exactly as the proof of Proposition 1:
One first shows that the mapping  $a_0\otimes
a_1\otimes\dots\otimes a_n\mapsto a_0d_u(a_1)\dots d_u(a_n)$
induces  a left module isomorphism of $\cala\otimes
(\otimes^n(\cala/\mathbb  C\bbbone))$ onto $\Omega^n_u(\cala)$
which implies that $a_0 d_u(a_1)\dots  d_u(a_n)\mapsto a_0
c(a_1,\dots,a_n)$ defines a left module homomorphism $i_c$ of 
$\Omega^n_u(\cala)$ into $\calm$; the cocycle property of $c$
then implies that $i_c$ is  a bimodule homomorphism. Again
uniqueness is straightforward. Concerning the  last part, if
there is an extension $\tilde \imath_c$ to 
$\cala\otimes\Omega^{n-1}_u(\cala)$, then $c$ is the Hochschild
coboundary of $(a_1,\dots,a_{n-1})\mapsto \tilde
\imath_c(\bbbone\otimes d_u(a_1)\dots d_u(a_{n-1}))$ and conversely,
if  $c$ is the coboundary of a normalized $(n-1)$-cochain $c'$
then one defines an  extension $\tilde \imath_c$ by setting $\tilde
\imath_c(\bbbone\otimes d_u(a_1)\dots
d_u(a_{n-1}))=c'(a_1,\dots,a_{n-1})$.~$\square$\\

Thus, for each integer $n\geq 1$, the normalized
$n\mbox{-}$cocycle  $d^{\cup n}_u$, defined by  $d^{\cup
n}_u(a_1,\dots,a_n)=d_u(a_1)\dots  d_u(a_n)$, is universal
among the normalized Hochschild $n$-cocyles.\\

By its very construction, $\Omega_u(\cala)$ is a graded
subalgebra of the tensor algebra over $\cala$,
$T_\cala(\cala\otimes\cala)$, of  the $(\cala,\cala)$-bimodule
$\cala\otimes \cala$. Indeed $T^n_\cala(\cala\otimes \cala)$ is
the $(n+1)$-th tensor power (over $\mathbb C$)
$\otimes^{n+1}\cala$ of $\cala$ whereas
$\Omega^n_u(\cala)=T^n_\cala(\Omega^1_u(\cala))$ is the
intersection of the kernels of the $(\cala,\cala)$-bimodule
homomorphisms $m_k:\otimes^{n+1}\cala\rightarrow
\otimes^n\cala$ defined by \[ \begin{array}{ll}
m_1(x_0\otimes\dots\otimes x_n) & = x_0x_1\otimes
x_2\otimes\dots \otimes x_n\\ \dotfill& \dotfill \\
m_n(x_0\otimes\dots\otimes x_n) &=x_1\otimes\dots\otimes
x_{n-2}\otimes x_{n-1}x_n \end{array} \] (i.e. $m_k$ is the
product in $\cala$ of the consecutive factors $x_{k-1}$ and
$x_k$). It turns out that the differential of $\Omega_u(\cala)$
has an extension, again denoted by $d_u$, as a differential of
$T_\cala(\cala\otimes\cala)$ which is defined by \[
d_u(x_0\otimes\dots \otimes x_n)=\sum^{n+1}_{k=0}(-1)^k\ \ \ 
x_0\otimes \dots \otimes x_{k-1}\otimes \bbbone \otimes
x_k\otimes \dots \otimes x_n \] for $x_i\in \cala$ where the first term of the summation is $\bbbone\otimes
x_0\otimes\dots \otimes x_n$ and the last term is
$(-1)^{n+1}x_0\otimes\dots\otimes x_n\otimes\bbbone$ (by convention). So
equipped $T_\cala(\cala\otimes\cala)$ is a graded differential
algebra, in fact a differential calculus over $\cala$, and
$\Omega_u(\cala)$ is a graded-differential subalgebra.

\begin{lemma} The cohomologies of $T_\cala(\cala\otimes\cala)$
and of $\Omega_u(\cala)$ are trivial in the sense that one has
: $H^0(T_\cala(\cala\otimes\cala))=H^0(\Omega_u(\cala))=\mathbb
C$ and \linebreak[4] $H^n(T_\cala(\cala\otimes \cala)) =
H^n(\Omega_u(\cala))=0$ for $n\geq 1$. \end{lemma}

\noindent\underbar{Proof}. Define $\delta:\mathbb C\rightarrow
\cala$ by $\delta(\lambda)=\lambda\bbbone$, one has $d_u\circ
\delta=0$ so \[ 0\rightarrow \mathbb
C\stackrel{\delta}{\rightarrow}\cala\stackrel{d_u}{\rightarrow}\cala\otimes
\cala \stackrel{d_u}{\rightarrow}\dots
\stackrel{d_u}{\rightarrow}\otimes^{n+1}\cala
\stackrel{d_u}{\rightarrow} \otimes^{n+2}\cala
\stackrel{d_u}{\rightarrow}\dots \] is a cochain complex with
coboundary $d$ being $d_u$ or $\delta$. Let $\omega$ be a
linear form on $\cala$ such that $\omega(\bbbone)=1$ and define
$k$ by $k(\mathbb C)=0$ and by $k(x_0\otimes\dots\otimes
x_n)=\omega(x_0)x_1\otimes\dots \otimes x_n$ for $n\geq 0$. One
has $kd+dk=I$ which implies $H^n(T_\cala(\cala\otimes\cala))=0$
for $n\geq 1$ and
$H^0(T_\cala(\cala\otimes\cala))=H^0(\Omega_u(\cala))=\mathbb
C$. Then $H^n(\Omega_u(\cala))=0$ for $n\geq 1$ follows from
the fact that one has $k(\Omega^n_u(\cala))\subset
\Omega^{n-1}_u(\cala)$ for $n\geq 1$. $\square$\\

\noindent\underbar{Remark 9}. The graded differential algebra
$C(\cala,\cala)$ of $\cala$-valued\linebreak[4] Hochschild
cochains of $\cala$ (see in Section 4) is a differential
calculus over $\cala$. Therefore, by Proposition 5 there is a
unique homomorphism $\Phi$ of $\Omega_u(\cala)$ into
$C(\cala,\cala)$ of graded differential algebras which induces
the identity mapping of $\cala$ onto itself. This homomorphism
extends to $T_\cala(\cala\otimes\cala)$ i.e. as an homomorphism
$\Phi:T_\cala(\cala\otimes\cala)\rightarrow C(\cala,\cala)$ of
graded differential algebras which is given by
$\Phi(x_0\otimes\dots\otimes x_n)(y_1,\dots,y_n)=x_0 y_1
x_1\dots y_n x_n$, \cite{tm}. Notice that
$\Phi(\Omega_u(\cala))$ is contained in the graded differential
subalgebra $C_0(\cala,\cala)$ of the normalized cochains of
$C(\cala,\cala)$.\\

In Section 6 we have defined the central bimodule
$\Omega^1_Z(\cala)$  to be the quotient of $\Omega^1_u(\cala)$
by the bimodule  $[Z(\cala),\Omega^1_u(\cala)]$ and the
derivation $d_Z$ of $\cala$ into $\Omega^1_Z(\cala)$ to be the 
image of $d_u:\cala\rightarrow \Omega^1_u(\cala)$. Let $I_Z$ be
the closed  two-sided ideal of $\Omega_u(\cala)$ generated by
$[Z(\cala),\Omega^1_u(\cala)]$ i.e.  the two-sided ideal
generated by $[Z(\cala),\Omega^1_u(\cala)]$ and
$d_u([Z(\cala),\Omega^1_u(\cala)])$. The space $I_Z$ is a
graded  ideal which is closed and such that $I_Z\cap 
\Omega^1_u(\cala)=[Z(\cala),\Omega^1_u(\cala)]$ which implies
that the quotient $\Omega_Z(\cala)$ is a graded  differential
algebra which coincides in degree 1 with the above 
$\Omega^1_Z(\cala)$ and that its differential (the image of
$d_u$) extends $d_Z:\cala\rightarrow \Omega^1_Z(\cala)$; this
differential will be also denoted by $d_Z$.  By construction,
$\Omega_Z(\cala)$ is, as graded algebra, a quotient of  the
tensor algebra over $\cala$ of the central bimodule
$\Omega^1_Z(\cala)$; on  the other hand it is easily seen that
tensor products over $\cala$ of central  bimodules and
quotients of central bimodules are again central bimodules 
\cite{mdv:pm2} so the $(\cala,\cala)$-bimodules
$\Omega^n_Z(\cala)$ are central bimodules
$(\Omega_Z(\cala)=\oplusinf_n\Omega^n_Z(\cala))$ and therefore
the  graded differential algebra $\Omega_Z(\cala)$ is a 
differential calculus  over $\cala$ for $\alg_Z$. Proposition 5
has the following counterpart for  $\Omega_Z(\cala)$.

\begin{proposition} Any homomorphism $\varphi$ of unital
algebras of  $\cala$ into the subalgebra $\Omega^0$ of elements
of degree 0 of a graded  differential algebra $\Omega$ which is
such that $\varphi(z)d(\varphi(x))=d(\varphi(x))\varphi(z)$ 
for any $z\in Z(\cala)$ and $x\in \cala$, ($d$ being the
differential of  $\Omega$), has a unique extension $\tilde
\varphi_Z:\Omega_Z(\cala)\rightarrow \Omega$  as a homomorphism
of graded differential algebras. \end{proposition}

\noindent \underbar{Proof}. By Proposition 5,  $\varphi$ has a
unique  extension $\tilde \varphi:\Omega_u(\cala)\rightarrow
\Omega$ as homomorphism of  graded differential algebras. On
the other hand $\varphi(z) 
d(\varphi(x))=d(\varphi(x))\varphi(z)$ for $z\in Z(\cala)$ and
$x\in\cala$ implies that $\tilde\varphi$  vanishes on
$[Z(\cala),\Omega^1_u(\cala)]$ and therefore also on $I_Z$
since it  is a homomorphism of graded differential algebras.
Thus $\tilde \varphi$  factorizes through a homomorphism
$\tilde\varphi_Z:\Omega_Z(\cala)\rightarrow  \Omega$ of graded
differential algebras which extends $\varphi$. Uniqueness is 
also straightforward here. $\square$\\

Proposition 7 has the following corollaries.

\begin{corol} For any differential calculus $\Omega$ over
$\cala$ for  $\alg_Z$, there is a unique homomorphism
$j_\Omega:\Omega_Z(\cala)\rightarrow  \Omega$ of differential
algebras which induces the identity mapping of $\cala$  onto
itself. \end{corol}

In other words $\Omega_Z(\cala)$ is universal among the
differential  calculi over $\cala$ for $\alg_Z$ and this
universal property characterizes it (up  to an isomorphism).
This is why we shall refer to $\Omega_Z(\cala)$ as {\sl  the
universal differential calculus over $\cala$ for $\alg_Z$.}

\begin{corol} Any homomorphism $\varphi:\cala\rightarrow \calb$
of  unital algebras mapping the center $Z(\cala)$ of $\cala$
into the center  $Z(\calb)$ of $\calb$ has a unique extension 
$\Omega_Z(\varphi):\Omega_Z(\cala)\rightarrow \Omega_Z(\calb)$
as a homomorphism of graded differential algebras.  \end{corol}

In fact $\Omega_Z$ is a covariant functor from the category
$\alg_Z$  to the category $\dif$ of graded differential
algebras. \\

In Section 2 it was pointed out that the graded  center of a
graded  algebra is stable by the graded derivations. This
implies in particular that the  graded center $Z_{\gr}(\Omega)$
of a graded differential algebra $\Omega$ is  a graded
differential subalgebra of $\Omega$ which is graded
commutative. We  have defined a differential calculus over
$\cala$ for $\alg_Z$ to be a graded differential algebra
$\Omega$ such that $\Omega^0=\cala$ and such that the center 
$Z(\cala)$ of $\cala(=\Omega^0)$ is contained in the center of
$\Omega$ i.e. in  its graded center $Z_\gr(\Omega)$ since its
elements are of degree zero in  $\Omega$. It follows that if
$\Omega$ is a differential calculus over $\cala$ for 
$\alg_{Z}$ then the center $Z(\cala)$ of $\cala$ generates a
graded differential  subalgebra of $\Omega$ which is graded
commutative and is in fact a graded  differential subalgebra of
the graded center $Z_\gr(\Omega)$ of $\Omega$. This  applies in
particular to $\Omega_{Z}$. If $\cala$ is commutative then 
$\Omega_Z(\cala)$ is graded commutative since it is generated
by $\cala$ which coincides  then with its center. In this case
Proposition 7  has the following corollary.

\begin{corol} If $\cala$ is commutative $\Omega_{Z}(\cala)$
identifies canonically with the graded differential algebra
$\Omega_{\cala\vert  \mathbb C}$ of Cartan-de Rham-K\"ahler
exterior differential forms. \end{corol}

\noindent \underbar{Proof}. Let us recall that
$\Omega_{\cala\vert  \mathbb C}$ is the exterior algebra over
$\cala$ of the module  $\Omega^1_{\cala\vert \mathbb C}$ of
K\"ahler differential,  $\Lambda_{\cala} \Omega^1_{\cala\vert 
\mathbb C}$, equipped with the unique differential extending
the universal  derivation of $\cala$ into $\Omega^1_{\cala\vert
\mathbb C}$. From this definition  and the universality of the
derivation of $\cala$ into $\Omega^1_{\cala\vert  \mathbb C}$
(which identifies, in view of Corollary 1, with
$d_Z:\cala\rightarrow \Omega^1_{Z}(\cala)$) it follows that
$\Omega_{\cala\vert \mathbb C}$  is characterized by the
following universal property: {\sl Any  homomorphism $\psi$ of
$\cala$ into the subalgebra $\Omega^0$ of the elements of
degree 0  of a graded commutative differential algebra $\Omega$
has a unique extension $\tilde\psi:\Omega_{\cala\vert\mathbb
C}\rightarrow \Omega$ as a  homomorphism of graded
(commutative) differential algebras}.\\ Let us come back to 
the proof of Corollary 4. Since $\Omega_Z(\cala)$ is graded
commutative with $\Omega^0_Z(\cala)=\cala$, the above universal
property implies that  there is a unique homomorphism of graded
differential algebras of  $\Omega_{\cala\vert\mathbb C}$ into
$\Omega_Z(\cala)$ which induces the identity mapping of 
$\cala$ onto itself. On the other hand Proposition 7 (or
Corollary 2) implies that  there is a unique homomorphism of
graded differential algebras of  $\Omega_Z(\cala)$ into
$\Omega_{\cala\vert\mathbb C}$ which induces the identity of
$\cala$  onto itself. Using again these two universal
properties, it follows that the above homomorphisms are inverse
isomorphisms. $\square$\\

If $\cala$ is commutative the cohomology of
$\Omega_Z(\cala)=\Omega_{\cala\vert\mathbb C}$ if often called
the de Rham cohomology (\cite{jll}, \cite{hus}) in spite of the
fact that as explained in Remark 8, for $\cala=C^\infty(M)$,
$\Omega_{\cala\vert\mathbb C}$ can be slightly different from
the algebra of smooth differential forms and that therefore
there is an ambiguity. Nevertheless $\Omega_Z(\cala)$ can be
considered as  a generalization of the graded differential 
algebra of differential forms which has the great advantage
that the  correspondence $\cala\mapsto \Omega_Z(\cala)$ is
functorial (Corollary~3). In contrast to the cohomology of
$\Omega_u(\cala)$, (see Lemma 5), the cohomology $H_Z(\cala)$
of $\Omega_Z(\cala)$ is generally non trivial. Since
$H_Z(\cala)$ is a noncommutative generalization of the de Rham
cohomology and since, by construction, $\cala\mapsto
H_Z(\cala)$ is a covariant functor from the category $\alg_Z$
to the category of graded algebras, it is natural to study the
properties of this cohomology.\\

Let $\Der(\cala)$ denote the vector space of all derivations
of  $\cala$ into itself. This vector space is a Lie algebra for
the bracket  $[\cdot,\cdot]$ defined by
$[X,Y](a)=X(Y(a))-Y(X(a))$ for $X,Y\in \Der(\cala)$ and  $a\in
\cala$. In view of Proposition~1, (universal property of 
$(\Omega^1_u(\cala),d_u)$), for each $X\in \Der(\cala)$ there
is a unique bimodule homomorphism
$i_X:\Omega^1_u(\cala)\rightarrow \cala$ for which $X=i_X\circ
d_u$.  This homomorphism of $\Omega^1_u(\cala)$ into
$\cala=\Omega^0_u(\cala)$  has a unique extension as an
antiderivation of 
$\Omega_u(\cala)=T_\cala(\Omega^1_u(\cala))$. This
antiderivation which will be again denoted by $i_X$ is of
degree  $-1$, (i.e. it is a graded derivation of degree $-1$).
It is not hard to verify  that $X\mapsto i_X$ is an operation
of the Lie algebra $\Der(\cala)$ in  the graded differential
algebra $\Omega_u(\cala)$, (see Section 2 for the  definition).
The corresponding Lie derivative $L_X=i_Xd_u+d_ui_X$ is for 
$X\in\Der(\cala)$ a derivation of degree 0 of $\Omega_u(\cala)$
which extends $X$. This  operation will be refered to as {\sl
the canonical operation of $\Derth(\cala)$ in
$\Omega_u(\cala)$}.\\

Let $X\in\Der(\cala)$ be a derivation of $\cala$ and let $z\in 
Z(\cala)$ and $\omega\in\Omega^1_u(\cala)$ one has \[ 
i_X([z,\omega])=[z,i_X(\omega)]=0 \] and \[
i_X(d([z,\omega]))=L_X([z,\omega])=[X(z),\omega]+[z,L_X(\omega)]=[z,L_X(\omega)]
\]  since $Z(\cala)$ is stable by the derivations of $\cala$.
This  implies that $i_X(I_Z)\subset I_Z$ and therefore that the
antiderivation $i_X$  passes to the quotient and defines an
antiderivation of degree $-1$ of  $\Omega_Z(\cala)$ which will
be again denoted by $i_X$. Notice that this (abuse of)
notation  is coherent with the one used in Proposition 3,
($\cala$ is obviously a central  bimodule). The corresponding
mapping $X\mapsto i_X$ of $\Der(\cala)$ into the
antiderivations of degree $-1$ of $\Omega_Z(\cala)$ is again
an  operation (the quotient of the one in $\Omega_u(\cala)$)
which will be refered to as  {\sl the canonical operation of
$\Derth(\cala)$ in $\Omega_Z(\cala)$}.\\

Finally if $\cala$ is a $\ast$-algebra,
$T_\cala(\cala\otimes\cala)$ is a graded differential
$\ast$-algebra if one equips it with the involution defined by
$(x_0\otimes\dots\otimes x_n)^\ast =$\linebreak[4] $(-1)^{\frac{n(n+1)}{2}}x^\ast_n\otimes\dots \otimes
x_0^\ast$. Since $\Omega_u(\cala)$ is stable by this
involution, it is also a graded differential $\ast$-algebra,
\cite{slw}. Furthermore $[Z(\cala),\Omega^1_u(\cala)]$ is
$\ast$-invariant which implies that the involution of
$\Omega_u(\cala)$ passes to the quotient and induces an
involution on $\Omega_Z(\cala)$ for which $\Omega_Z(\cala)$
also becomes a graded differential $\ast$-algebra. More generally in this case, a differential calculus $\Omega$ over $\cala$ will always be assumed to be equipped with an involution extending the involution of $\cala$ and such that it is a graded differential $\ast$-algebra, (notice that if $\Omega$ is generated by $\cala$ such an involution is unique).

\section{Diagonal and derivation-based calculi}

Let $\cala$ be a unital associative complex algebra and let
$\calm$ be an arbitrary $(\cala,\cala)$-bimodule. Then the set
$\hom^\cala_\cala(\calm,\cala)$ of all bimodule homomorphisms
of $\calm$ into $\cala$ is a module over the center $Z(\cala)$
of $\cala$ which will be refered to as the $\cala$-{\sl dual
of} $\calm$ and denoted by $\calm^{\ast_\cala}$,
\cite{mdv:pm1}, \cite{dv:4}. Conversely, if $\caln$ is a
$Z(\cala)$-module the set $\hom_{Z(\cala)}(\caln,\cala)$ of all
$Z(\cala)$-module homorphisms of $\caln$ into $\cala$ is
canonically a $(\cala,\cala)$-bimodule which will be also
refered to as the $\cala$-{\sl dual of} $\caln$ and denoted by
$\caln^{\ast_\cala}$. The $\cala$-dual of a $Z(\cala)$-module
is clearly a central bimodule over $\cala$ so the above duality
between $(\cala,\cala)$-bimodules and $Z(\cala)$-modules can be
restricted to a duality between the central bimodules over
$\cala$ and the $Z(\cala)$-modules. This latter duality
generalizes the duality between modules over a commutative
algebra, \cite{mdv:pm1}, \cite{dv:4}. Indeed, if $\cala$ is
commutative both central bimodules over $\cala$ and
$Z(\cala)$-modules coincide with $\cala$-modules and the above
duality is then the usual duality between $\cala$-modules. Let
us come back to the general situation and let $\calm$ be a
$(\cala,\cala)$-bimodule; then one obtains by evaluation a {\sl
canonical homomorphism} of $(\cala,\cala)$-bimodule
$c:\calm\rightarrow \calm^{\ast_\cala\ast_\cala}$ of $\calm$
into {\sl its $\cala$-bidual
$\calm^{\ast_\cala\ast_\cala}=(\calm^{\ast_\cala})^{\ast_\cala}$}.

\begin{lemma} The following properties (a) and (b) are
equivalent for a\linebreak[4] $(\cala,\cala)$-bimodule $\calm$.\\ (a) The
canonical homomorphism $c:\calm\rightarrow
\calm^{\ast_\cala\ast_\cala}$ is injective.\\ (b) $\calm$ is
isomorphic to a subbimodule of $\cala^I$ for some set $I$.
\end{lemma}

\noindent\underbar{Proof}. $(a)\Rightarrow (b)$. By definition
$\calm^{\ast_\cala\ast_\cala}$ is a subbimodule of $\cala^I$
with $I=\hom_{Z(\cala)}(\calm^{\ast_\cala},\cala)$ so if $c$ is
injective $\calm$ is isomorphic  to a subbimodule of
$\calm^{\ast_\cala\ast_\cala}$ and therefore also to a
subbimodule of $\cala^I$.

\noindent\phantom{\underbar{Proof}.}               
$(b)\Rightarrow (a)$. Let $\varphi$ be a bimodule homomorphism
of $\cala$ into itself. One has
$\varphi(a)=a\varphi(\bbbone)=\varphi(\bbbone)a$ which implies
$\varphi(\bbbone)\in Z(\cala)$. Conversely any $z\in Z(\cala)$
defines a bimodule homomorphism $\varphi$ of $\cala$ into
itself by setting $\varphi(a)=az$ (i.e. $\varphi(\bbbone)=z)$.
It follows that $\cala^{\ast_\cala}=Z(\cala)$. Let $\Phi$ be a
$Z(\cala)$-module homomorphism of $Z(\cala)$ into $\cala$. Then
$\Phi(z)=z\Phi(\bbbone)$ with $\Phi(\bbbone)\in\cala$.
Conversely any $a\in\cala$ defines such a $Z(\cala)$-module
homomorphism $\Phi$ by setting $\Phi(z)=za$ (i.e.
$\Phi(\bbbone)=a)$. It follows that
$Z(\cala)^{\ast_\cala}=\cala$ and therefore
$\cala^{\ast_\cala\ast_\cala}=\cala$. This immediately implies
that if $\calm\subset\cala^I$ as subbimodule then
$c:\calm\rightarrow\calm^{\ast_\cala\ast_\cala}$ is injective.
$\square$\\

An $(\cala,\cala)$-bimodule $\calm$ satisfying the equivalent
conditions of Lemma 6 will be said to be a {\sl diagonal
bimodule} over $\cala$, \cite{mdv:pm1}, \cite{mdv:pm2} (see
also in \cite{dv:4}). A diagonal bimodule is central but the
converse is not generally true. The $\cala$-dual of an
arbitrary  $Z(\cala)$-module is a diagonal bimodule. Every
subbimodule of a diagonal bimodule is diagonal, every product
of diagonal bimodules is diagonal and the tensor product over
$\cala$ of two diagonal bimodules is diagonal.\\

If $\cala$ is commutative, a diagonal bimodule over $\cala$ is
simply a $\cala$-module such that the canonical homomorphism in
its bidual is injective. In particular in this case a
projective module is diagonal (as a bimodule for the underlying
structure).\\

\noindent\underbar{Remark 10}. It is a fortunate circumstance
which is easy to verify that, for a $Z(\cala)$-module $\caln$,
the biduality does not depend on $\cala$ but only on
$Z(\cala)$. That is one has
$\caln^{\ast_\cala\ast_\cala}=\caln^{\ast\ast}$ and the
canonical homomorphism $c:\caln\rightarrow\caln^{\ast\ast}$
obtained by evaluation for the $\cala$-duality reduces to the
usual one for a module over the commutative algebra
$Z(\cala)$.\\

Let $\calm$ be a $(\cala,\cala)$-bimodule then the canonical
image $c(\calm)$ of $\calm$ in its $\cala$-bidual is a diagonal
bimodule. The diagonal bimodule $c(\calm)$ is the universal
``diagonalization" of $\calm$ in the sense that it is
characterized (among the diagonal bimodules over $\cala$) by
the following universal property, \cite{mdv:pm1},
\cite{mdv:pm2}. 

\begin{proposition} For any homomorphism of
$(\cala,\cala)$-bimodules\linebreak[4] 
$\varphi:\calm\rightarrow\caln$ of $\calm$ into a diagonal
bimodule $\caln$ over $\cala$, there is a unique homomorphism
of $(\cala,\cala)$-bimodules $\varphi_c:c(\calm)\rightarrow
\caln$ such that $\varphi=\varphi_c\circ c$. \end{proposition}

\noindent\underbar{Proof}. In view of the definition and Lemma
6 (b), it is sufficient to prove the statement for
$\caln=\cala^I$ for some set $I$, which is then equivalent to
the statement for $\caln=\cala$. On the other hand, for
$\caln=\cala$,
$\varphi\in\hom^\cala_\cala(\calm,\cala)=\calm^{\ast_\cala}$
and one has $\varphi(m)=<c(m),\varphi>=\varphi_c(c(m))$ for
$m\in\calm$ (by the definitions of
$\calm^{\ast_\cala\ast_\cala}$ and of the evaluation $c$) which
defines $\varphi_c$ uniquely. $\square$\\

One has $c(\Omega^1_u(\cala))=c(\Omega^1_Z(\cala))$ and we
shall denote by $\Omega^1_\diag(\cala)$ this diagonal bimodule
and by $d_\diag$ the derivation $c\circ d_u$ (or equivalently
$c\circ d_Z$) of $\cala$ into $\Omega^1_\diag(\cala)$.

\begin{proposition} For any first order differential calculus
$(\Omega^1,d)$ over $\cala$ such that $\Omega^1$ is diagonal,
there is a unique bimodule homomorphism $i_d$ of
$\Omega^1_{\diagth}(\cala)$ into $\Omega^1$ such that
$d=i_d\circ d_\diagth$. \end{proposition}

\noindent\underbar{Proof}. In view of the above universal
property of $c(\Omega^1_u(\cala))$, the corresponding canonical
homomorphism of $\Omega^1_u(\cala)$ into $\Omega^1$ (as in
Proposition 1) factorizes through a unique homomorphism
$i_d:\Omega^1_\diag(\cala)\rightarrow \Omega^1$. $\square$\\

In other words, the derivation
$d_\diag:\cala\rightarrow\Omega^1_\diag(\cala)$ of $\cala$ into
the diagonal bimodule $\Omega^1_\diag(\cala)$ is universal for
the derivations of $\cala$ into diagonal bimodules over
$\cala$.\\

Let us recall (see Section 3) that the vector space
$\Der(\cala)$ of all derivations of $\cala$ into itself is a
Lie algebra and also a $Z(\cala)$-module and that
$\Omega_\der(\cala)$ was defined to be the graded differential
subalgebra of $C_\wedge(\Der(\cala),\cala)$ generated by
$\cala$ whereas $\os_\der(\cala)$ was defined to be the graded
differential subalgebra of $C_\wedge(\Der(\cala),\cala)$ which
consists of cochains of $\Der(\cala)$ which are
$Z(\cala)$-multilinear. Clearly $C_\wedge^n(\Der(\cala),\cala)$
is  diagonal so the first order differential calculus
$(C^1_\wedge(\Der(\cala),\cala),d)$ satisfies the conditions of
Proposition 9 which implies that there is a unique bimodule
homomorphism $i_d$ of $\Omega^1_\diag(\cala)$ into
$C^1_\wedge(\Der(\cala),\cala)$ for which $d=i_d\circ d_\diag$.

\begin{proposition} The homomorphism $i_d:\Omega^1_\diagth(\cala)\rightarrow
C^1_\wedge(\Der(\cala),\cala)$ is injective, so by
identifying $\Omega^1_\diagth(\cala)$ with its image (by
$i_d$), one has canonically:
\[
\Omega^1_\diagth(\cala)=\Omega^1_\derth(\cala),(\Omega^1_\diagth(\cala))^{\ast_\cala}=\Derth(\cala)\ 
\mbox{and}\ 
(\Omega^1_\diagth(\cala))^{\ast_\cala\ast_\cala}=\os^1_\derth(\cala).
\]
\end{proposition}

\noindent\underbar{Proof}. Applying Proposition 9  for
$\Omega^1=\cala$ leads to the identification\linebreak[4]
$\hom^\cala_\cala(\Omega^1_\diag(\cala),\cala)=\Der(\cala)$
that is $(\Omega^1_\diag(\cala))^{\ast_\cala}=\Der(\cala)$,
(notice that one has also
$(\Omega^1_u(\cala))^{\ast_\cala}=\Der(\cala))$. By definition
one has
$\os^1_\derth(\cala)=$ \linebreak[4] $\hom_{Z(\cala)}(\Der(\cala),\cala)$ that
is $\os^1_\der(\cala)=(\Der(\cala))^{\ast_\cala}$ and therefore
\linebreak[4]
$(\Omega^1_\diag(\cala))^{\ast_\cala\ast_\cala}=\os^1_\der(\cala)$.
On the other hand one has
$i_d(\Omega^1_\diag(\cala))=\Omega^1_\der(\cala)$ since
$\Omega^1_\der(\cala)$ is generated by $\cala$ (as bimodule).
The injectivity of $i_d$ follows from the fact that
$\Omega^1_\diag(\cala)$ is diagonal i.e. that the canonical
homomorphism in its $\cala$-bidual is injective. $\square$\\

Notice that by definition one also has
$(\bigwedge^n_{Z(\cala)}\Der(\cala))^{\ast_\cala}=\os^n_\der(\cala)$.\\

Let $I_\diag$ be the closed two-sided ideal of
$\Omega_u(\cala)$ generated by the kernel of the canonical
homomorphism $c$ of $\Omega^1_u(\cala)$ into its
$\cala$-bidual. The ideal $I_\diag$ is graded such that
$I_\diag\cap \Omega^0_u(\cala)=0$ and $I_\diag\cap
\Omega^1_u(\cala)=\ker(c)$ which implies that the quotient
$\Omega^1_u(\cala)/I_\diag$ is a graded differential algebra
which is a differential calculus over $\cala$ and coincides in
degree 1 with $c(\Omega^1_u(\cala))=\Omega^1_\diag(\cala)$. This
differential calculus will be refered to as the {\sl diagonal
calculus} and denoted by $\Omega_\diag(\cala)$. The
differential of $\Omega_\diag(\cala)$ is the image of $d_u$ and
extends the derivation $d_\diag:\cala\rightarrow
\Omega^1_\diag(\cala)$; this differential will be also denoted
by $d_\diag$. Proposition 5 and Proposition 7 have the
following counterpart for $\Omega_\diag(\cala)$.

\begin{proposition}  Any homomorphism $\varphi$ of unital
algebras of  $\cala$ into the subalgebra $\Omega^0$ of elements
of degree 0 of a graded  differential algebra $\Omega$ which is
such that $d(\cala)$ spans a diagonal bimodule over $\cala$
(for the $(\cala,\cala)$-bimodule structure on $\Omega^1$
induced by $\varphi$) has a unique extension $\tilde
\varphi_\diagth:\Omega_\diagth(\cala)\rightarrow \Omega$  as a
homomorphism of graded differential algebras. \end{proposition}

\noindent \underbar{Proof}. By Proposition 5,  $\varphi$ has a
unique  extension $\tilde \varphi:\Omega_u(\cala)\rightarrow
\Omega$ as homomorphism of  graded differential algebras. On
the other hand the assumption means that $d:\cala\rightarrow
\tilde\varphi(\Omega^1_u(\cala))$ is a derivation and that
$\tilde\varphi(\Omega^1_u(\cala))$ is a diagonal bimodule over
$\cala$ so, in view of Proposition 9, the homomorphism
$\tilde\varphi:\Omega^1_u(\cala)\rightarrow\Omega^1$ factorizes
through a homomorphism
$\tilde\varphi^1_\diag:\Omega^1_\diag(\cala)\rightarrow
\Omega^1$. Thus $\tilde\varphi$ vanishes on $\ker(c)$ and
therefore on $I_Z$ since it is a homomorphism of graded
differential algebras so it factorizes through a homomorphism
$\tilde\varphi_\diag:\Omega_\diag(\cala)\rightarrow\Omega$ of
graded differential algebras. Uniqueness is again
straightforward. $\square$\\

Thus $\Omega_\diag(\cala)$ is also characterized by a universal
property like $\Omega_u(\cala)$ and $\Omega_Z(\cala)$ but in
contrast to the cases of $\Omega_u(\cala)$ and
$\Omega_Z(\cala)$, the correspondence
$\cala\mapsto\Omega_\diag(\cala)$ has no obvious functorial
property. The reason for this is the fact that the diagonal
bimodules are not the bimodules for a category of algebras in
the sense explained in Section 5.\\

Proposition 11 implies in particular that one has a unique
homomorphism of graded differential algebra of
$\Omega_\diag(\cala)$ into $\Omega_\der(\cala)$ which extends
the identity mapping of $\cala$ onto itself. This homomorphism
$\Omega_\diag(\cala)\rightarrow \Omega_\der(\cala)$ is
surjective since $\Omega_\der(\cala)$ is generated by $\cala$
as differential algebra. Furthermore in degree 1 it is, in view
of Proposition 10, a bimodule isomorphism of
$\Omega^1_\diag(\cala)$ onto $\Omega^1_\der(\cala)$. However,
for $m\geq 2$, the corresponding bimodule homomorphism of
$\Omega^m_\diag(\cala)$ onto $\Omega^m_\der(\cala)$ is not
generally injective (i.e. it has a non trivial kernel).\\

For instance when $\cala$ coincides with the algebra
$M_n(\mathbb C)$ of complex $n\times n$ matrices one has
\[
\begin{array}{ll}
&\Omega_u(M_n(\mathbb C))=\Omega_Z(M_n(\mathbb
C))=\Omega_\diag(M_n(\mathbb C))\simeq\\
\\
& C_0(M_n(\mathbb C),M_n(\mathbb C))=M_n(\mathbb C)\otimes T
\slg(n,\mathbb C)^\ast
\end{array}
\]
whereas 
\[
\Omega_\der(M_n(\mathbb
C))=C_\wedge(\slg(n,\mathbb C),M_n(\mathbb C))=M_n(\mathbb
C)\otimes \bigwedge \slg(n,\mathbb C)^\ast.
\]
In fact, in this
case, the homomorphism $\Phi$ of Remark 9 is an isomorphism
which induces the isomorphism of $\Omega_u(M_n(\mathbb C))$
onto the differential algebra $C_0(M_n(\mathbb C),M_n(\mathbb
C))$ of normalized Hochschild cochains; the latter being
identical as graded algebra to the tensor product $M_n(\mathbb
C)\otimes T \slg(n,\mathbb C)^\ast$ of $M_n(\mathbb C)$ with
the tensor algebra over $\mathbb C$ of the dual of
$\slg(n,\mathbb C)$, (concerning $\Omega^1_\der(M_n(\mathbb
C))=\Omega^1_u(M_n(\mathbb C))$, and $\Omega_\der(M_n(\mathbb
C))=C_\wedge(\slg(n,\mathbb C),M_n(\mathbb C))$, see in
\cite{dv:2}).\\

In the case where $\cala$ is the algebra $C^\infty(M)$ of
smooth functions on a good smooth manifold (finite dimensional
paracompact, etc.) then one has
$\Omega_\diag(C^\infty(M))=\Omega_\der(C^\infty(M))$
$(=\os_\der(C^\infty(M)))$.\\

It is not hard to show that the operations of the Lie algebra
$\Der(\cala)$ in $\Omega_u(\cala)$ and in $\Omega_Z(\cala)$
pass to the quotient to define an operation of $\Der(\cala)$ in
the graded differential algebra $\Omega_\diag(\cala)$ which will be again refered to as {\sl the canonical operation of $\Derth(\cala)$ in $\Omega_\diagth(\cala)$}.
Furthermore, all these operations of $\Der(\cala)$ pass to the
quotient to define an operation of $\Der(\cala)$ in
$\Omega_\der(\cala)$ which coincides with the canonical
operation of $\Der(\cala)$ in $\Omega_\der(\cala)$ defined in
Section 3.\\

One has the following commutative diagram of surjective
homomorphisms of graded differential algebras which is also a
diagram of homomorphisms of the operations of $\Der(\cala)$.

\[ \begin{diagram} \node{\Omega_u(\cala)}\arrow{e} \arrow{s}
\arrow{se}  \node{\Omega_Z(\cala)} \arrow{s} \arrow{sw} \\   
\node{\Omega_\diag(\cala)} \arrow{e} \node{\Omega_\der(\cala)}
\end{diagram} \]

Furthermore, if $\cala$ is a $\ast$-algebra there is a
canonical involution on $\Omega_\diag(\cala)$ such that this
diagram is also a diagram of graded differential
$\ast$-algebras, (the involutions of
$\Omega_u(\cala),\Omega_Z(\cala)$ and $\Omega_\der(\cala)$ have
been defined previously in Section 7 and Section 3).

\section{Noncommutative symplectic geometry and quantum
mechanics}

Let $\cala$ be as before a unital associative complex algebra.
A {\sl Poisson bracket} on $\cala$ is a Lie bracket which is a
biderivation on $\cala$ (for its associative product). That is
$(a,b)\mapsto\{a,b\}$ is a Poisson bracket if it is a bilinear
antisymmetric mapping of $\cala\times\cala$ into $\cala$ (i.e.
a linear mapping of $\bigwedge^2\cala$ into $\cala$) which
satisfies \[ \begin{array}{ll}
\{\{a,b\},c\}+\{\{b,c\},a\}+\{\{c,a\},b\}=0 & \mbox{(Jacobi
identity)}\\ \{a,bc\}=\{a,b\}c+b\{a,c\} & \mbox{(derivation
property)} \end{array} \] for any elements $a,b,c$ of $\cala$.
Equipped with such a Poisson bracket, $\cala$ is refered to as
a {\sl Poisson algebra}, \cite{fl}.\\

There is  a lot of classical commutative Poisson algebras, for
instance the symmetric algebra $S(\fracg)$ (over $\mathbb C$)
of a (complex) Lie algebra $\fracg$, the algebra $C^\infty(M)$
of smooth functions on a symplectic manifold, etc.. For a
noncommutative algebra $\cala$, a generic type of Poisson
bracket $\{\cdot,\cdot\}$ is obtained by setting for $a,b\in
\cala$  \[ \{a,b\}=\frac{i}{\hbar}[a,b] \] where $[a,b]$ denotes
the commutator in $\cala$, i.e. $[a,b]=ab-ba$, and where
$\hbar\in \mathbb C$ is any non zero complex number. We have
put a $i\in \mathbb C$ in front of the right-hand side of the
above formula in order that in the case where $\cala$ is a
$\ast$-algebra {\sl the Poisson bracket is real}, i.e.
satisfies $\{a,b\}^\ast=\{a^\ast,b^\ast\}$, whenever $\hbar$ is
real. The reason why the Poisson brackets proportional to the
commutator are quite generic (in the noncommutative case) is
connected to the following lemma \cite{fl}.

\begin{lemma} Let $\cala$ be a Poisson algebra, then one has
$[a,b]\{c,d\}=\{a,b\}[c,d]$ and more generally
$[a,b]x\{c,d\}=\{a,b\}x[c,d]$ for any elements $a,b,c,d$ and
$x$ of $\cala$. \end{lemma}

\noindent\underbar{Proof}. The first identity is obtained by
developing $\{ac,bd\}$ in two different orders by using the
biderivation property. The second (more general since
$\bbbone\in \cala$) identity is obtained by replacing $c$ by
$xc$ in the first identity, by developing and by using again
the first identity. $\square$\\

For more details concerning the ``generic side" of Poisson
brackets proportional to the commutator we refer to \cite{fl}.
We simply observe here that this is the type of Poisson
brackets which occurs in quantum mechanics.\\

Our aim is now to develop a (noncommutative) generalization of
symplectic structures related to the above Poisson brackets.
One should start from a notion of differential form i.e. from a
differential calculus $\Omega$ over $\cala$. Since for a
Poisson bracket $x\mapsto\{a,x\}$ is an element of
$\Der(\cala)$ for any $a\in \cala$, it is natural to assume
that one has an operation $X\mapsto i_X$ of the Lie algebra
$\Der(\cala)$ in the graded differential algebra $\Omega$.
Furthermore we wish to take into account the structure of
$Z(\cala)$-module of $\Der(\cala)$ so we require that $\Omega$
is a central bimodule over $\cala$ and that $X\mapsto i_X$ is a
$Z(\cala)$-linear mapping of $\Der(\cala)$ into
$\Der^{-1}_\gr(\Omega)$. Notice that this $Z(\cala)$-linearity
is well defined since $\Omega$ central is equivalent to
$Z(\cala)\subset Z^0_\gr(\Omega)$, (see in Section 2 for the
notations). Having such a differential calculus, one defines a
homomorphism $\lambda$ of $\Omega$ into $\os_\der(\cala)$ by
setting $\lambda(\omega)(X_1,\dots,X_n)=i_{X_n}\dots
i_{X_1}\omega$ for $\omega\in\Omega^n$. The fact that this
defines a homomorphism of graded differential algebra of
$\Omega$ into $C_\wedge(\Der(\cala),\cala)$ follows from the
general properties of operations whereas the fact that the
image of $\lambda$ is contained in $\os_\der(\cala)$ follows
from the $Z(\cala)$-linearity. It turns out that even if one
uses a general differential calculus $\Omega$ for the
symplectic structures, the only relevant parts for the
corresponding Poisson structures are the images by $\lambda$ in
$\os_\der(\cala)$, (see e.g. in \cite{fl}). One is then led to
the definitions of \cite{dv:3}, or more precisely to the
following slight generalizations \cite{dv:4}.\\

An element $\omega$ of $\os^2_{\der}(\cala)$ will be said to be
{\sl nondegenerate} if, for any $x\in\cala$, there is a
derivation $\ham(x)\in \Der(\cala)$ such that one has
$\omega(X,\ham(x))=X(x)$ for any $X\in \Der(\cala)$.  Notice
that if $\omega$ is nondegenerate then $X\mapsto i_X\omega$ is
an injective linear mapping of $\Der(\cala)$ into
$\os^1_{\der}(\cala)$ but that the converse is not true; the
condition for $\omega$ to be nondegenerate is stronger than the
injectivity of $X\mapsto i_X\omega$. If $M$ is a manifold, an
element $\omega\in \os^2_{\der}(C^\infty(M))$ is an ordinary
2-form on $M$ and it is nondegenerate in the above sense if and
only if the 2-form $\omega$ is nondegenerate in the classical
sense (i.e. everywhere nondegenerate).\\

Let $\omega\in\os^2_{\der}(\cala)$ be nondegenerate, then for a
given $x\in\cala$ the derivation $\ham(x)$ is unique and
$x\mapsto\ham(x)$ is a linear mapping of $\cala$ into
$\Der(\cala)$.\\

A closed nondegenerate element $\omega$ of
$\os^2_{\der}(\cala)$ will be called {\sl a symplectic
structure for} $\cala$.

\begin{lemma} Let $\omega$ be a symplectic structure for
$\cala$ and let us define an antisymmetric bilinear bracket on
$\cala$ by $\{x,y\}=\omega(\hamth(x),\hamth(y))$. Then
$(x,y)\mapsto \{x,y\}$ is a Poisson bracket on $\cala$.
\end{lemma}

\noindent\underbar{Proof}.  One has
$\{x,yz\}=\{x,y\}z+y\{x,y\}$ for $x,y,z\in \cala$. Furthermore
one has the identity  \[
d\omega(\ham(x),\ham(y),\ham(z))=-\{x,\{y,z\}\}-\{y,\{z,x\}\}
-\{z,\{x,y\}\} \] which implies the Jacobi identity since
$d\omega=0$. $\square$\\

Let $\omega$ be a symplectic structure for $\cala$, then one
has \[ [\ham(x),\ham(y)]=\ham(\{x,y\}), \] i.e. $\ham$ is a
Lie-algebra homomorphism of $(\cala,\{,\})$ into
$\Der(\cala)$.  We shall refer to the above bracket as {\sl the
Poisson bracket associated to the symplectic structure
$\omega$}. If $\cala$ is a $\ast$-algebra and if furthermore
$\omega$ is real, i.e. $\omega=\omega^\ast$, then this Poisson
bracket is real and $\ham(x^\ast)=(\ham(x))^\ast$ for any
$x\in\cala$.\\

An algebra $\cala$ equipped with a symplectic structure will be
refered to as a {\sl symplectic algebra}. Thus, symplectic
algebras are particular Poisson algebras.\\

\noindent\underbar{Remark 11}. Let $\cala$ be an arbitrary
Poisson algebra with Poisson bracket $(x,y)\mapsto\{x,y\}$; one
defines a linear mapping $\ham:\cala\rightarrow \Der(\cala)$ by
$\ham(x)(y)=\{x,y\}$, (i.e. $\ham(x)=\{x,\cdot\}$), for $x,y\in
\cala$. In this general setting one also has the identity
$[\ham(x),\ham(y)]=\ham(\{x,y\})$ since it is equivalent to the
Jacobi identity for the Poisson bracket.\\

If $M$ is a manifold, a symplectic structure for $C^\infty(M)$
is just a symplectic form on $M$. Since there are manifolds
which do not admit symplectic form, one cannot expect that an
arbitrary $\cala$ admits a symplectic structure.\\

Assume that $\cala$ has a trivial center $Z(\cala)=\mathbb 
C\bbbone$ and that all its derivations are inner (i.e. of the
form $ad(x),x\in\cala$). Then one defines an element $\omega$
of $\os^2_{\der}(\cala)$ by setting
$\omega(ad(ix),ad(iy))=i[x,y]$. It is easily seen that $\omega$
is a symplectic structure for which one has $\ham(x)=ad(ix)$
and $\{x,y\}=i[x,y]$. If furthermore $\cala$ is a
$\ast$-algebra, then this symplectic structure is real
($\omega=\omega^\ast$). Although a little tautological, this
construction is relevant for quantum mechanics.\\

Let $\cala$ be, as above, a complex unital $\ast$-algebra with
a trivial center and only inner derivations and assume that
there exists a linear form $\tau$ on $\cala$ which is central,
i.e. $\tau(xy)=\tau(yx)$, and normalized by $\tau(\bbbone)=1$.
Then one defines an element $\theta\in\os^1_{\der}(\cala)$ by
$\theta(ad(ix))=x-\tau(x)\bbbone$. One has 
$(d\theta)(ad(ix),ad(iy))=i[x,y]$, i.e. $\omega=d\theta$, so in
this case the symplectic form $\omega$ is exact. As examples of
such algebras one can take $\cala=M_n(\mathbb C)$, (a factor of
type I$_n$), with $\tau=\frac{1}{n}$ trace, or $\cala=\calr$, a
von Neumann algebra which is a factor of type II$_1$ with
$\tau$  equal to the  normalized trace. The algebra
$M_n(\mathbb C)$ is the algebra of observables of a quantum
spin $s=\frac{n-1}{2}$ while $\calr$ is the  algebra used to
describe the observables of an infinite assembly of quantum
spins; two typical types of quantum systems with no classical
counterpart.\\

Let us now consider the C.C.R. algebra (canonical commutation
relations) $\cala_{CCR}$ \cite{dv:3}. This is the complex
unital $\ast$-algebra generated by two hermitian elements $q$
and $p$ satisfying the relation $[q,p]=i\hbar\bbbone$. This
algebra is the algebra of observables of the quantum
counterpart of a classical system with one degree of freedom.
We keep here the positive constant $\hbar$ (the Planck
constant) in the formula for comparison with classical
mechanics, although the algebra for $\hbar\not=0$ is isomorphic
to the one with $\hbar=1$. We restrict attention to one degree
of freedom to simplify the notations but the discussion extends
easily to a finite number of degrees of freedom. This algebra
has again only inner derivations and a trivial center so
$\omega(ad(\frac{i}{\hbar}x),ad(\frac{i}{\hbar}y))=\frac{i}{\hbar}[x,y]$
defines a symplectic structure for which
$\ham(x)=ad(\frac{i}{\hbar}x)$ and
$\{x,y\}=\frac{i}{\hbar}[x,y]$ which is the standard quantum
Poisson bracket. In this case one can express $\omega$ in terms
of the generators $q$ and $p$ and their differentials
\cite{dv:3}, \cite{dv:4}: $$\omega=\sum_{n\geq
0}\left(\frac{1}{i\hbar}\right)^n \frac{1}{(n+1)!}[\dots
[dp,\underbrace{p],\dots,p]}_n[\dots[dq,\underbrace{q],\dots,q}_n]$$
Notice that this formula is meaningful because if one inserts 
two derivations $ad(ix),ad(iy)$ in it, only a finite number of
terms contribute to the sum. In contrast to the preceding case,
here the symplectic form is not exact, i.e. it corresponds to a
non vanishing element of $H^2(\os_\der(\cala_{CCR}))$ which is
therefore non trivial. This was guessed in \cite{dv:3} on the
basis of the nonexistence of a finite trace (i.e. central
linear form) on $\cala_{CCR}$ and finally proved in \cite{fl}.
For $\hbar=0$, $q$ and $p$ commute and the algebra reduces to
the algebra of complex polynomial functions on the phase space
$\mathbb R^2$. Furthermore the limit of
$\{x,y\}=\frac{i}{\hbar}[x,y]$ at $\hbar=0$ reduces to the
usual classical Poisson bracket as well known and, by using the
above formula, one sees that the formal limit of $\omega$ at
$\hbar=0$ is $dpdq$. This limit is however very singular since
the limit algebra is the algebra of complex polynomials in two
indeterminates, the limit symplectic form is exact and not
every derivation is hamiltonian in contrast to what happens for
$\cala_{CCR}$ (i.e. for $\hbar\not= 0$).

\section{Theory of connections}

Throughout this section, $\cala$ is a unital associative
complex algebra and $\Omega$ is a differential calculus over
$\cala$, that is a graded differential algebra such that
$\Omega^0=\cala$ with differential denoted by $d$.\\

Let $\calm$ be a left $\cala$-module; a $\Omega$-{\sl
connection} on $\calm$ (or simply a {\sl connection} on $\calm$
if no confusion arises) is a linear mapping
$\nabla:\calm\rightarrow\Omega^1\otimes_\cala\calm$ such that
one has \[ \nabla(am)=a\nabla(m)+d(a)\otimes_\cala m \] for any
$a\in\cala$ and $m\in\calm$, ($\Omega^1\otimes_\cala\calm$
being equipped with its canonical structure of left
$\cala$-module). One extends $\nabla$ to
$\Omega\otimes_\cala\calm$ by setting
$\nabla(\omega\otimes_\cala
m)=(-1)^n\omega\nabla(m)+d(\omega)\otimes_\cala m$ for
$\omega\in\Omega^n$ and $m\in\calm$ ($\Omega\otimes_\cala\calm$
is canonically a left $\Omega$-module). It then follows from
the definitions that $\nabla^2$ is a left $\Omega$-module
endomorphism of $\Omega\otimes_\cala\calm$ which implies that
its restriction $\nabla^2:\calm\rightarrow
\Omega^2\otimes_\cala\calm$ to $\calm$ is a homomorphism of
left $\cala$-modules; this homomorphism is called the {\sl
curvature} of the connection $\nabla$.\\

Not every left $\cala$-module admits a connection. If $\calm$
is the free $\cala$-module $\cala\otimes E$, where $E$ is some
complex vector space, then $\nabla=d\otimes I_E$ is a
connection on $\cala\otimes E$ which has a vanishing curvature
(such a connection with zero curvature is said to be flat). If
$\calm\subset \cala\otimes E$ is a direct summand of a free
$\cala$-module $\cala\otimes E$ and if $P:\cala\otimes
E\rightarrow \calm$ is the corresponding projection, then
$\nabla=P\circ(d\otimes I_E)$ is a connection on $\calm$. Thus
a projective module admits (at least one) a connection. In the
case where $\Omega$ is the universal differential calculus
$\Omega_u(\cala)$ the converse is also true: It was shown in
\cite{cq} that a (left) $\cala$-module admits a
$\Omega_u(\cala)$-connection if and only if it is projective.\\

One defines in a similar manner $\Omega$-connections on right
modules. Namely if $\caln$ is a right $\cala$-module, a
$\Omega$-connection on $\caln$ is a linear mapping
$\nabla$ of $\caln$ into $\caln\otimes_\cala\Omega^1$ such that
$\nabla(na)=\nabla(n)a+n\otimes_\cala d(a)$ for any $n\in
\caln$ and $a\in \cala$.\\

Let $\calm$ be a left $\cala$-module, then its dual
$\calm^\ast$ (i.e. the set of left $\cala$-module homomorphisms
of $\calm$ into $\cala$) is a right $\cala$-module. We denote
by $<m,n>\in\cala$ the evaluation of $n\in\calm^\ast$ on
$m\in\calm$. Let $\nabla$ be a $\Omega$-connection on $\calm$,
then one defines a unique linear mapping $\nabla^\ast$ of
$\calm^\ast$ into $\calm^\ast\otimes_\cala\Omega^1$ by setting
(with obvious notations)  \[
<m,\nabla^\ast(n)>=d(<m,n>)-<\nabla(m),n> \] for any
$m\in\calm$ and $n\in\calm^\ast$. It is easy to verify that
$\nabla^\ast$ is a $\Omega$-connection on the right module
$\calm^\ast$ which will be refered to as {\sl the dual
connection} of $\nabla$. One defines in the same way the dual
connection of a connection on a right module.\\

Our aim is now to recall the definitions of hermitian modules
over a\linebreak[4] $\ast$-algebra $\cala$ and of hermitian connections. We
assume that $\cala$ is a $\ast$-algebra such that the convex
cone $\cala^+$ generated by the $a^\ast a$ ($a\in\cala$) is a
strict cone i.e. such that $\cala^+\cap(-\cala^+)=0$. This
property is satisfied for instance by $\ast$-algebras of
operators in Hilbert spaces. A {\sl hermitian structure} on a
right $\cala$-module $\calm$ \cite{connes:02} is a sesquilinear
mapping $h:\calm\times \calm\rightarrow \cala$ such that one
has:\\ \[ \begin{array}{ll} (i) & h(ma,nb)=a^\ast h(m,n)b,\ \
\forall m,n\in \calm \ \mbox{and}\ \forall a,b\in \cala\\ (ii)
& h(m,m)\in \cala^+, \ \ \forall m\in \calm \ \mbox{and}\ 
h(m,m)=0\Rightarrow m=0. \end{array} \] A right $\cala$-module
$\calm$ equipped with a hermitian structure will be refered to
as a {\sl hermitian module} over $\cala$. If $\calm$ is a
hermitian module over $\cala$, a {\sl hermitian connection} on
$\calm$ is a connection $\nabla$ on the right $\cala$-module
$\calm$ such that one has \[ d(h(m,n))=h(\nabla m,n)+h(m,\nabla
n) \] for any $m,n\in \calm$ with obvious notations. We have
chosen to define hermitian structures on right modules for
notational reasons, (we prefer the convention of physicists for
sesquilinearity, i.e. linearity in the second argument); one
can define similarily hermitian structures and connections for
left modules.\\

Let $\calm$ be a right $\cala$-module. The group $\aut(\calm)$
of all module automorphisms of $\calm$ acts on the affine space
of all connections on $\calm$ via $\nabla\mapsto
\nabla^U=U\circ \nabla \circ U^{-1}$, $U\in\aut(\calm)$, (one
canonically has $\aut(\calm)\subset
\aut(\calm\otimes_\cala\Omega^1)$). If furthermore $\cala$ is a
$\ast$-algebra as above and if $h$ is a hermitian structure on
$\calm$, then the subgroup of $\aut(\calm)$ of all
automorphisms $U$ which preserve $h$, i.e. such that
$h(Um,Un)=h(m,n)$ for $m,n\in \calm$, will be denoted by
$\aut(\calm,h)$ and called the {\sl gauge group} whereas its
elements will be called {\sl gauge transformations}; it acts on
the real affine space of hermitian connections on $\calm$.\\

As pointed out before, one-sided modules are not sufficient and
one needs bimodules for a lot of reasons. Firstly, in the case
where $\cala$ is a $\ast$-algebra, one needs $\ast$-bimodules
to formulate and discuss reality conditions \cite{mdv:pm1},
\cite{connes:05}, \cite{dv:4} (see also in the introduction).
Secondly, a natural noncommutative generalization of linear
connections should be connections on $\Omega^1$, since $\Omega$
is taken as an analog of differential forms, but this is a
$(\cala,\cala)$-bimodule in an essential way. Thirdly, in order
to have an analog of local couplings, one needs to have a
tensor product over $\cala$ since the latter is the
noncommutative version of the local tensor product of tensor
fields. In short one needs a theory of connections for
bimodules and any of the above quoted problems shows that
one-sided connections on bimodules (i.e. on bimodules
considered as left or right modules) are of no help. The
difficulty to define a $\Omega$-connection on a
$(\cala,\cala)$-bimodule $\calm$ lies in the fact that a left
$\cala$-module connection on $\calm$ sends $\calm$ into
$\Omega^1\otimes_\cala\calm$ whereas a right $\cala$-module
connection on $\calm$ sends $\calm$ into
$\calm\otimes_\cala\Omega^1$. A solution of this problem
adapted to the case where $\calm=\Omega^1$ has been given in
\cite{jm} and generalized in \cite{mdv:m} for arbitrary
$(\cala,\cala)$-bimodules on the basis of an analysis of first
order differential operators in bimodules. This led to the
following definition \cite{mdv:m}.\\

Let $\calm$ be a $(\cala,\cala)$-bimodule; a {\sl left bimodule
$\Omega$-connection} on $\calm$ is a left $\cala$-module
$\Omega$-connection $\nabla$ on $\calm$ for which there is a
bimodule homomorphism
$\sigma:\calm\otimes_\cala\Omega^1\rightarrow\Omega^1\otimes_\cala\calm$
such that \[ \nabla(ma)=\nabla(m)a+\sigma(m\otimes_\cala d(a))
\] for any $a\in \cala$ and $m\in \calm$. Clearly $\sigma$ is
then unique under these conditions. One defines similarily a
{\sl right bimodule $\Omega$-connection} on $\calm$ to be a
right $\cala$-module $\Omega$-connection $\nabla$ on $\calm$
for which there is a bimodule homomorphism
$\sigma:\Omega^1\otimes_\cala\calm\rightarrow
\calm\otimes_\cala \Omega^1$ such that \[
\nabla(am)=a\nabla(m)+\sigma(d(a)\otimes_\cala m) \] for any
$a\in\cala$ and $m\in\calm$. When no confusion arises on
$\Omega$ and on ``left-right" we simply refer to this notion as
{\sl bimodule connection}.\\

In the case where $\calm$ is the bimodule $\Omega^1$ itself, a
left bimodule $\Omega$-connection is just the first part of the
proposal of \cite{jm} for the definition of linear connections
in noncommutative geometry; the second part of the proposal of
\cite{jm} relates $\sigma$ and the product
$\Omega^1\otimes_\cala\Omega^1\rightarrow \Omega^2$ so it makes
sense only for $\calm=\Omega^1$ and is there necessary to
define the generalization of torsion.\\

It has been shown in \cite{bmds} (Appendix A of \cite{bmds})
that on general grounds, the above definition is just what is
needed to define tensor products over $\cala$ of bimodule
connections and of left (right) bimodule connections with left
(right) module connections. In fact, let $\nabla'$ be a left
bimodule connection on the bimodule $\calm'$ and let $\nabla''$
be a connection on a left module $\calm''$. Then one defines a
connection $\nabla$ on the left module 
$\calm'\otimes_\cala\calm''$ by setting \[
\nabla=\nabla'\otimes_\cala I_{\calm''}+(\sigma' \otimes_\cala
I_{\calm''})\circ (I_{\calm'}\otimes_\cala \nabla'') \] where
$\sigma':\calm'\otimes_\cala\Omega^1\rightarrow
\Omega^1\otimes_\cala\calm'$ is the bimodule homomorphism
corresponding to $\nabla'$. If furthermore $\calm''$ is a
$(\cala,\cala)$-bimodule and if $\nabla''$ is a left bimodule
connection with corresponding bimodule homomorphism
$\sigma'':\calm''\otimes_\cala\Omega^1\rightarrow
\Omega^1\otimes_\cala\calm''$, then $\nabla$ is also a left
bimodule connection with corresponding bimodule homomorphism
$\sigma$ given by \[ \sigma = (\sigma'\otimes_{\cala}
I_{\calm''})\circ (I_{\calm'}\otimes_\cala\sigma'') \] of
$\calm'\otimes_\cala\calm'' \otimes_\cala \Omega^1$ into
$\Omega^1\otimes_\cala \calm'\otimes_\cala \calm''.$\\

Let $\calm$ be  a $(\cala,\cala)$-bimodule and let $\calm^\ast$
denote the dual of $\calm$ considered as a left $\cala$-module.
Then $\calm^\ast$ is a right $\cala$-module as dual of a left
$\cala$-module, but it is in fact a bimodule if one defines the
left action $m'\mapsto am'$ of $\cala$ on $\calm^\ast$ by
$<m,am'>=<ma,m'>$ for any $m\in\calm, a\in\cala, m'\in
\calm^\ast$. If $\nabla$ is a left bimodule $\Omega$-connection
on $\calm$  then one verifies that $\nabla^\ast$ is a right
bimodule $\Omega$-connection on $\calm^\ast$ \cite{bmds}
(Appendix B of \cite{bmds}). Notice that this kind of duality
between bimodules is different of the $\cala$-duality between
bimodules over $\cala$ and modules over $Z(\cala)$ discussed in
Section 8.\\

When $\cala$ is a $\ast$-algebra, there is also a generalization of hermitian forms on $(\cala,\cala)$-bimodules which has been introduced on \cite{rief} and called {\sl right hermitian forms} in \cite{mdv:pm1} which is adapted for tensor products over $\cala$. If $\calm$ is a $(\cala,\cala)$-bimodule, then a right hermitian form on $\calm$ is a hermitian form $h$ on $\calm$ considered as a right $\cala$-module which is such that for the left multiplication by $a\in\cala$ one has $h(m,an)=h(a^\ast m,n)$. One can then define the notion of right hermitian bimodule connection, (which is in particular a right bimodule connection).\\

We now explain the relation between the above notion of
bimodule connection and the theory of first order operators in
bimodules. Let $\cala$ and $\calb$ be unital associative
complex algebras and let $\calm$ and $\caln$ be two
$(\cala,\calb)$-bimodules. We denote by $l_a$ the left
multiplication by $a\in\cala$ in $\calm$ and in $\caln$ and we
denote by $r_b$ the right multiplication by $b\in \calb$ in
$\calm$ and in $\caln$. A linear mapping $D$ of $\calm$ into
$\caln$ which is such that one has $[[D,l_a],r_b]=0$ for any
$a\in \cala$ and $b\in\calb$ is called {\sl a first-order
operator} or {\sl an operator of order 1} of $\calm$ into
$\caln$ \cite{connes:03}. Notice that homomorphisms of left
$\cala$-modules of $\calm$ into $\caln$ as well as
homomorphisms of right $\calb$-modules of $\calm$ into $\caln$
are first-order operators of $\calm$ into $\caln$. The
structure of first-order operators is given by the following
theorem \cite{mdv:m}.

\begin{theo} Let $\calm$ and $\caln$ be two
$(\cala,\calb)$-bimodules and let $D$ be a first order operator
of $\calm$ into $\caln$. Then, there is a unique
$(\cala,\calb)$-bimodule homomorphism $\sigma_L(D)$ of
$\Omega^1_u(\cala) \otimesinf_{\cala}\calm$ into $\caln$ and
there is a unique $(\cala,\calb)$-bimodule homomorphism
$\sigma_R(D)$ of $\calm\otimesinf_{\calb}\Omega^1_u(\calb)$
into $\caln$ such that one has: \[ D(amb)=aD(m)b +
\sigma_L(D)(d_ua\otimes m)b + a\sigma_R(D)(m\otimes d_ub) \]
for any $m\in \calm$, $a\in\cala$ and $b\in \calb$. \end{theo}
For the proof and further informations, see in \cite{mdv:m}. It
is clear that $\sigma_L(D)$ and $\sigma_R(D)$ are the
appropriate generalization of the notion of symbol in this
setting. We shall refer to them as the {\sl left} and the {\sl
right universal symbol of} $D$ respectively.\\

\noindent\underbar{Remark 12}. The converse of Theorem 4 is
also true. More precisely, let $(\Omega^1_L,d_L)$ be a first
order differential calculus over $\cala$, let
$(\Omega^1_R,d_R)$ be a first order differential calculus over
$\calb$ and let $D:\calm\rightarrow \caln$ be a linear mapping
then any of the following condition $(1)$ or $(2)$ implies that
$D$ is a first-order operator.\\

$(1)$ There is a $(\cala,\calb)$-bimodule homomorphism
$\sigma_L:\Omega^1_L\otimes_\cala\calm\rightarrow \caln$ such
that  \[ D(am)=aD(m)+\sigma_L(d_L(a)\otimes m),\ \ \forall
m\in\calm\ \mbox{and}\ \forall a\in\cala \]

 $(2)$ There is a $(\cala,\calb)$-bimodule homomorphism
$\sigma_R:\calm\otimes_\calb\Omega^1_R\rightarrow \caln$ such
that \[ D(mb)=D(m)b +\sigma_R(m\otimes d_R(b)),\ \ \forall m\in
\calm\ \mbox{and}\  \forall b \in \calb. \]

Let $\calm$ be a $(\cala,\cala)$-bimodule and let $\nabla$ be a
left $\cala$-module $\Omega$-connection on $\calm$. It is
obvious that $\nabla$ is a first-order operator of the
$(\cala,\cala)$-bimodule $\calm$ into the
$(\cala,\cala)$-bimodule $\Omega^1\otimes_\cala \calm$. It
follows therefore from the above theorem that there is a unique
$(\cala,\cala)$-bimodule homomorphism
$\sigma_R(\nabla)$ of $\calm\otimes_\cala\Omega^1_u(\cala)$ into
 $\Omega^1\otimes_\cala\calm$ such that one has \[
\nabla(ma)=\nabla(m)a+\sigma_R(\nabla)(m\otimes_\cala d_u(a))
\] for any $m\in\calm$ and $a\in\cala$. Therefore, $\nabla$ is
a left bimodule $\Omega$-connection on $\calm$ if and only if
$\sigma_R(\nabla)$ factorizes through a
$(\cala,\cala)$-bimodule homomorphism
$\sigma:\calm\otimes_\cala\Omega^1\rightarrow
\Omega^1\otimes_\cala\calm$ as $\sigma_R(\nabla)=\sigma\circ
(I_\calm\otimes i_d)$ where $I_\calm$ is the identity mapping
of $\calm$ onto itself and $i_d$ is the unique
$(\cala,\cala)$-bimodule homomorphism of $\Omega^1_u(\cala)$
into $\Omega^1$ such that $d=i_d\circ d_u$ (see Proposition 1).
This implies in particular that any left $\cala$-module
$\Omega_u(\cala)$-connection is a left bimodule
$\Omega_u(\cala)$-connection.\\

In the case of the derivation-based differential calculus,
there is an easy natural way to define connections on left and
right modules and on central bimodules over $\cala$,
\cite{mdv:pm1}. We describe it in the case of central bimodules
(for left and for right modules, just forget multiplications on
the right and on the left respectively). Let $\calm$ be a
central bimodule over $\cala$, i.e. a $\cala$-bimodule for
$\alg_Z$, {\sl a (derivation-based) connection on} $\calm$ is a
linear maping $\nabla,X\mapsto\nabla_X$, of $\Der(\cala)$ into
the linear endomorphisms of $\calm$ such that  \[
\nabla_{zX}(m)=z\nabla_X(m),
\nabla_X(amb)=a\nabla_X(m)b+X(a)mb+amX(b) \] for any
$m\in\calm$, $X\in\Der(\cala)$, $z\in Z(\cala)$ and
$a,b\in\cala$. One verifies that such a connection on the
central bimodule $\calm$ is a bimodule
$\os_\der(\cala)$-connection on $\calm$ in the previous sense
with a well defined $\sigma$, (modulo some technical problems
of completion of the tensor products
$\os^1_\der(\cala)\otimes_\cala\calm$ and
$\calm\otimes_\cala\os^1_\der(\cala)$). The interest of this
formulation is that curvature is straightforwardly defined and
is a bimodule homomorphism \cite{mdv:pm1}. We refer to
\cite{mdv:pm1} (and also to \cite{dv:4}) for more details and
in particular for the relation with $\cala$-duality.
Furthermore, in this frame the notion of reality on connections
is obvious. Assume that $\cala$ is a $\ast$-algebra and that
$\calm$ is a central bimodule which is a $\ast$-bimodule over
$\cala$ then a (derivation-based) connection $\nabla$ on
$\calm$ will be said to be {\sl real} if one has
$\nabla_X(m^\ast)=(\nabla_X(m))^\ast$ for any $m\in \calm$ and
any $X\in \Der_{\mathbb R}(\cala)$, i.e. $X\in\Der(\cala)$ with
$X=X^\ast$.\\

The notion of reality in the general frame of bimodule
$\Omega$-connections is slightly more involved and will not be
discussed here.

\section{Classical Yang-Mills-Higgs models}

An aspect with no counterpart in ordinary differential geometry
of the theory of $\Omega$-connections on $\cala$-modules for a
differential calculus $\Omega$ which is not graded commutative
is the generic occurrence of inequivalent $\Omega$-connections
with vanishing curvature (on a fixed $\cala$-module). By taking
as algebra $\cala$ the algebra of functions on space-time with
values in some algebra $\cala_0$, i.e. $\cala=C^\infty(\mathbb
R^{s+1})\otimes\cala_0$, this led to classical Yang-Mills-Higgs
models based on noncommutative geometry in which the Higgs
field is the part of the connection which is in the
``noncommutative directions".\\

In the following, we display the case of $\Omega$-connections
on right modules over the algebra $\cala=C^\infty(\mathbb
R^{s+1})\otimes M_n(\mathbb C)$ of smooth $M_n(\mathbb
C)$-valued functions on $\mathbb R^{s+1}$ for
$\Omega=\Omega_\der(\cala)$.\\

Let us first describe the situation for ${\cal A}=M_n(\mathbb
C)$. The derivations of $M_n(\mathbb C)$ are all inner so the
complex Lie algebra ${\rm Der}(M_n(\mathbb C))$ reduces to
$\slg (n)$ and the real Lie algebra ${\rm Der}_\mathbb
R(M_n(\mathbb C))$ reduces to $\sug(n)$. As already mentioned
in Section 8, one has  \[ \Omega_\der(M_n(\mathbb
C))=C_\wedge({\rm Der} M_n(\mathbb C), M_n(\mathbb
C))=C_\wedge(\slg(n), M_n(\mathbb C)) \] as can be shown
directly \cite{dv:2} and as also follows from the formulas below.
Let $E_k, k\in \{1,2,\dots,n^2-1\}$ be a base of self--adjoint
traceless $n\times n$--matrices. The $\partial_k={\rm
ad}(iE_k)$ form a basis of real derivations i.e. of
${\rm Der}_\mathbb R(M_n(\mathbb C))=\sug(n)$. One has
$[\partial_k,\partial_\ell]=C^m_{k\ell}\partial_m$, the
$C^m_{k\ell}$ are the corresponding structure constants of
$\sug(n)$, (or $\slg(n)$). Define $\theta^k \in
\Omega^1_\der(M_n(\mathbb C))$ by
$\theta^k(\partial_\ell)=\delta^k_\ell \bbbone$.The following
formulas give a presentation of the graded differential algebra
$\Omega_\der(M_n(\mathbb C))$ \cite{dvkm:1}, \cite{dv:3}: \[
\begin{array}{lll} E_kE_\ell & = &
g_{k\ell}\bbbone+(S^m_{k\ell}-{i\over2} C^m_{k\ell})E_m \\ \\
E_k\theta^\ell & = & \theta^\ell E_k \\ \\ \theta^k \theta^\ell
& = & -\theta^\ell \theta^k\\ \\  dE_k & = &
-C^m_{k\ell}E_m\theta^\ell\\ \\ d\theta^k & = & -{1\over 2}
C^k_{\ell m}\theta^\ell \theta^m  \end{array} \] where
$g_{k\ell}=g_{\ell k}$,$\;$ $S^m_{k\ell}=S^m_{\ell k}$ are
real, $g_{k\ell}$ are the components of the Killing form of
$\sug(n)$ and $C^m_{k\ell}=- C^m_{\ell k}$ are as above the
(real) structure constants of $\sug(n)$.Formula giving the 
$dE_k$ can be inverted and one has  
\[
\theta^k=-{i\over n^2}
g^{\ell m}g^{kr} E_\ell E_r dE_m
\]
where $g^{k\ell}$ are the
components of the inverse matrix of $(g_{k\ell})$. The element
$\theta=E_k\theta^k$ of $\Omega^1_\der(M_n(\mathbb C))$ is
real, $\theta=\theta^\ast$, and independent of the choice of
the $E_k$, in fact we already met $\theta$ in Section~9: 
$\theta({\rm ad}(iA))=A-{1\over n} {\rm tr} (A)\bbbone$ and
$\omega=d\theta$ is the natural symplectic structure for
$M_n(\mathbb C)$. Furthermore $\theta$ is invariant, $L_X
\theta=0$, and any invariant element of
$\Omega^1_\der(M_n(\mathbb C))$ is a scalar multiple of
$\theta$. We call $\theta$ the {\it canonical invariant
element} of $\Omega^1_\der(M_n(\mathbb C))$. One has
$$dM=i[\theta, M], \; \forall M\in M_n(\mathbb C) $$ $$d(-
i\theta)+(- i\theta)^2=0.$$

The $\ast$--algebra $M_n(\mathbb C)$ is simple with only one
irreductible representation in $\mathbb C^n$. A general finite
right--module (which is projective) is the space
$M_{Kn}(\mathbb C)$ of $K \times n$--matrices with right action
of $M_n(\mathbb C)$. Then ${\rm Aut}(M_{Kn}(\mathbb
C))$ is the group $GL(K)$ acting by left matrix multiplication. The module
$M_{Kn}(\mathbb C)$ is naturally hermitian with $h(\Phi,
\Psi)=\Phi^\ast \Psi$ where $\Phi^\ast$ is the $n \times K$
matrix hermitian conjugate to $\Phi$. The gauge group is then
the unitary group $U(K) (\subset GL(K))$. Here, there is a
natural origin $\Onabla$ in the space of connections given by
$\Onabla\Phi=- i\Phi \theta$ where $\Phi \in M_{Kn}(\mathbb C)$
and where $\theta$ is the canonical invariant element of
$\Omega^1_\der(M_n(\mathbb C))$. The fact that this defines a
connection follows from  \[ \Onabla(\Phi M) = (\Onabla \Phi)M
+\Phi i[\theta,M] \] and from the above expression of $dM$ for
$M\in M_n(\mathbb C)$. This connection is hermitian and its
follows from the above expression for $d\theta$ that its
curvature vanishes, i.e. $(\Onabla)^2=0$. Any connection
$\nabla$ is of the form $\nabla \Phi= \Onabla \Phi + A \Phi$
where $A=A_k \theta^k$ with $A_k \in M_K(\mathbb C)$ and
$A\Phi$ means $A_k \Phi \otimes \theta^k$. The connection
$\nabla$ is hermitian if and only if the $A_k$ are
antihermitian i.e. $A_k^\ast=-A_k$. The curvature of $\nabla$
is given by $\nabla^2 \Phi=F\Phi\; (= F_{k\ell} \Phi \otimes
\theta^k \theta^\ell)$ with  $$F ={1\over 2}([A_k, A_\ell] -
C^m_{k\ell} A_m) \theta^k \theta^\ell.$$ Thus $\nabla^2=0$ if
and only if the $A_k$ form a representation of the Lie algebra
$\slg(n)$ in $\mathbb C^K$ and two such connections are in the
same ${\rm Aut}(M_{Kn}(\mathbb C))$--orbit if and only if the
corresponding representations of $\slg(n)$ are equivalent. This
implies that the {\it gauge orbits of flat} $(\nabla ^2=0)$
{\it hermitian connections are in one--to--one correspondence
with unitary classes of representations of $\sug(n)$ in
$\mathbb C^K$}, \cite{dvkm:1}. For instance if $n=2$, these
orbits are labelled by the number of partitions of the integer
$K$.\\

We now come to the case $\cala=C^\infty(\mathbb R^{s+1})\otimes
M_n(\mathbb C)$. Let $x^\mu, \; \mu \in \{0,1,\dots,s\}$, be
the canonical coordinates of $\mathbb R^{s+1}$. One has\break 
$\Omega_\der(C^\infty(\mathbb R^{s+1}) \otimes M_n(\mathbb
C))=\Omega_\der(C^\infty(\mathbb R^{s+1}))\otimes
\Omega_\der(M_n(\mathbb C))$ so one can split the differential
as $d=d'+d''$ where $d'$ is the differential along $\mathbb
R^{s+1}$ and $d''$ is the differential of
$\Omega_\der(M_n(\mathbb C))$. A typical finite projective
right module is $C^\infty(\mathbb R^{s+1}) \otimes
M_{Kn}(\mathbb C)$. This is an hermitian module with hermitian
structure given by $h(\Phi,\Psi)(x)=\Phi(x)^\ast \Psi(x),\;
(x\in \mathbb R^{s+1})$. As a $C^\infty(\mathbb
R^{s+1})$--module, this module is free (of rank $K.n$), so
$d'\Phi$ is well defined for $\Phi \in C^\infty(\mathbb
R^{s+1}) \otimes M_{Kn}(\mathbb C)$. In fact,
$d'\Phi(x)={\partial \Phi \over \partial x^\mu}(x) dx^\mu$. A
connection on the $C^\infty(\mathbb R^{s+1}) \otimes
M_n(\mathbb C)$--module $C^\infty(\mathbb R^{s+1}) \otimes
M_{Kn}(\mathbb C)$ is of the form $\nabla \Phi =d'\Phi -i\Phi
\theta + A\Phi$ with $A=A_\mu dx^\mu + A_k\theta^k$, where the
$A_\mu$ and the $A_k$ are $K \times K$ matrix valued functions
on $\mathbb R^{s+1}$ (i.e. elements of $C^\infty (\mathbb
R^{s+1}) \otimes M_K(\mathbb C))$ and where
$A\Phi(x)=A_\mu(x)\Phi(x) dx^\mu +A_k(x)\Phi(x)\theta^k$. Such
a connection is hermitian if and only if the $A_\mu(x)$ and the
$A_k(x)$ are antihermitian, $\forall x \in \mathbb R^{s+1}$.
The curvature of $\nabla$ is given by $\nabla^2\Phi=F\Phi$
where  \[ \begin{array}{ll} F={1\over 2}(\partial_\mu
A_\nu-\partial_\nu A_\mu  &+[A_\mu,A_\nu])dx^\mu dx^\nu\\ \\
&+(\partial_\mu A_k + [A_\mu ,A_k])dx^\mu \theta^k\\ \\
&+{1\over 2}([A_k,A_\ell]-C^m_{k\ell}A_m)\theta^k \theta^\ell 
\end{array} \]

The connection $\nabla$ is flat (i.e. $\nabla^2=0$) if and only
if each term of the above formula vanishes which implies that
$\nabla$ is gauge equivalent to a connection for which one has
$A_\mu=0, \ \  \partial_\mu A_k=0$ and
$[A_k,A_\ell]=C^m_{k\ell}A$. Furthermore two such connections
are equivalent if and only if the corresponding representations
of $\sug(n)$ in $\mathbb C^K$ (given by the constant $K \times
K$--matrices $A_\ell$) are equivalent. So again, {\it the gauge
orbits of flat hermitian connections are in one-to-one
correspondence with the unitary classes of (antihermitian)
representations of $\sug(n)$ in $\mathbb C^K$}. Again, in the
case n=2, the number of such orbits is the number of partitions
of the integer$K$ i.e. $${\rm card}\{(n_r) \vert \sum_r \ \ 
n_r.r=K\}.$$ \\

If we consider $\mathbb R^{s+1}$ as the ($s+1$)--dimensional
space--time and if  we replace the algebra of smooth functions
on $\mathbb R^{s+1}$ by $C^\infty(\mathbb R^{s+1}) \otimes
M_n(\mathbb C)$ which we interpret as the algebra of ``smooth
functions on a noncommutative generalized space-time". It is
clear, from the above expression for the curvature that the
generalization of the (euclidean) Yang--Mills action for a
hermitian connection $\nabla$ on $C^\infty (\mathbb R^{s+1})
\otimes M_{Kn}(\mathbb C)$ is \[  \begin{array}{ll} \Vert
F\Vert^2 &= \int d^{s+1}x \, {\rm tr}\Bigl\lbrace {1\over 4}
\sum (\partial_\mu A_\nu-\partial_\nu A_\mu +[A_\mu,A_\nu])^2\\
&+{1\over 2} \sum(\partial_\mu A_k + [A_\mu,A_k])^2 + {1\over
4} \sum ([A_k,A_\ell]-C^m_{k\ell} A_m)^2\Bigl\rbrace 
\end{array} \] where the metrics of space-time is $g_{\mu
\nu}=\delta_{\mu \nu}$ and where the basis $E_k$ of hermitian
traceless $n \times n$--matrices  is chosen in such a way that
$g_{k\ell}=\delta_{k\ell}$, i.e. ${\rm
tr}(E_kE_\ell)=n\delta_{k\ell}$. This can be more deeply
justified by introducing the analog of the Hodge involution on
$\Omega_\der(M_n(\mathbb C))$, the analog of the integration of
elements of $\Omega^{n^2-1}_\der(M_n(\mathbb C))$ (essentially
the trace) and by combining these operations with the
corresponding one on $\mathbb R^{s+1}$ to obtain a scalar
product on $\Omega_\der(C^\infty(\mathbb R^{s+1}) \otimes
M_n(\mathbb C))$ etc. See in \cite{dvkm:1}, \cite{dvkm:2}  for
more details.\\

The above action is the Yang--Mills action on the
noncommutative space corresponding to $C^\infty(\mathbb
R^{s+1}) \otimes M_n(\mathbb C)$. However it can be interpreted
as the action of a field theory on the $(s+1)$--dimensional
space--time $\mathbb R^{s+1}$. At first sight, this field
theory consists of a $U(K)$-Yang-Mills potential $A_\mu(x)$
minimally coupled with scalar fields $A_k(x)$ with values in
the adjoint representation which interact among themselves
through a quartic potential. The action  is positive and
vanishes for $A_\mu=0$ and $A_k=0$, but is also vanishes on
other gauge orbits. Indeed $\Vert F \Vert^2=0$ is equivalent to
$F=0$, so the gauge orbits on which the action vanishes are
labelled by unitary classes of representations of $\sug(n)$ in
$\mathbb C^K$. By the standard semi--heuristic argument, these
gauge orbits are interpreted as different vacua for the
corresponding quantum theory. To specify a quantum theory, one
has to choose one and to translate the fields in
order that the zero of these translated fields corresponds to
the chosen vacuum (i.e. is the corresponding zero of the action). The
variables $A_\mu, A_k$  are thus adapted to the specific vacuum
$\varphi_0$ corresponding to the trivial representation $A_k=0$
of $\sug(n)$. If one chooses the vacuum $\varphi_\alpha$
corresponding to a representation $\Talpha R_k$ of $\sug(n)$,
(i.e. one has $[\Talpha R_k,\Talpha R_\ell]= C^m_{k\ell}
\Talpha R_n)$, one must instead use the variables $A_\mu$ and
$\Talpha B_k=A_k-\Talpha R_k$. Making this change of variable
one observes that components of $A_\mu$ become massive and that
the $\Talpha B_k$ have different masses; the whole mass
spectrum depends on $\alpha$. This is very analogous to the
Higgs mechanism. Here however the gauge invariance is not
broken, the non--invariance of the mass--terms of the $A_\mu$
is compensated by the fact that the gauge transformation of the
$\Talpha B_k$ becomes inhomogeneous (they are components of a
connection). Nevertheless, from the point of view of the
space-time intepretation this is the Higgs mechanism and the
$A_k$ are Higgs fields.\\

The above models were the first ones of classical
Yang-Mills-Higgs models based on noncommutative geometry. They
certainly admit a natural supersymmetric extension since there
is a natural extension of the derivation-based calculus to
graded matrix algebras \cite{gr}. There is also another
extension of the above calculus where $C^\infty(\mathbb
R^{s+1})\otimes M_n(\mathbb C)$ is replaced by the algebra
$\Gamma \fin (E)$ of smooth sections of the endomorphisms
bundle of a (nontrivial) smooth vector bundle $E$  (of rank
$n$) admitting a volume over a smooth ($(s+1)$-dimensional)
manifold \cite{mdv:m2}.\\

The use of the derivation-based calculus makes the above
models  quite rigid. By relaxing this i.e. by using other
differential calculi $\Omega$, other models based on
noncommutative geometry which are closer to the classical
version of the standard model have been constructed
\cite{connes:003}, \cite{colo}, \cite{rcoq}. Furthermore there
is an elegant way to combine the introduction of the (spinors)
matter fields with the differential calculus and the metric
\cite{connes:03} as well as with the reality conditions
\cite{connes:05} in noncommutative geometry, (and also with the
action principles \cite{chco}). Within this general set-up, one
can probably absorb any classical model of gauge theory.\\

A problem arises for the quantization of these classical models
based on noncommutative geometry. Namely is it possible to keep
something of the noncommutative geometrical interpretation of
these classical models at the quantum level? The best would be
to find some B.R.S. symmetry \cite{brs} ensuring that
(perturbative) quantization does not spoil the correspondence
with noncommutative geometry. Unfortunately no such symmetry
was discovered up to now. As long as no progress is obtained on
this problem, the noncommutative geometrical interpretation of
the gauge theory with Higgs field must be taken with some
circumspection in spite of its appealing features.

\section{Conclusion : Further remarks}

Concerning the noncommutative generalization of differential
geometry the point of view more or less explicit here is that the
data are encoded in an algebra $\cala$ which plays the role  of
the algebra of smooth functions. This is why although we have
described various notions in terms of an arbitrary differential
calculus $\Omega$, we have studied in some details specific
differential calculi ``naturally" associated with $\cala$ (i.e.
which do not depend on other data than $\cala$ itself) such as
the universal differential calculus $\Omega_u(\cala)$, the
generalization $\Omega_Z(\cala)$ of the K\"ahler exterior
forms, the diagonal calculus $\Omega_\diag(\cala)$ and the
derivation-based calculus. There are other possibilities, for
instance some authors consider that the data are encoded in a
graded differential algebra which plays the role of the algebra
of smooth differential forms, e.g. \cite{malt}. This latter
point of view can be taken into account here by using an
arbitrary differential calculus $\Omega$.\\

In all the above points of view, the generalization of
differential forms is provided by a graded differential
algebra. This is not always so natural. For instance it was
shown in \cite{kar} (see also \cite{kar:02}) that the subspace
$[\Omega_u(\cala),\Omega_u(\cala)]_\gr$ of graded commutators
in $\Omega_u(\cala)$ is stable by $d_u$ and that the cohomology
of the cochain complex
$\Omega_u(\cala)/[\Omega_u(\cala),\Omega_u(\cala)]$ is closely related to the cyclic homology (it is contained in the reduced cyclic homology), and is also in several respects  a noncommutative version of de Rham cohomology. This complex $\Omega_u(\cala)/[\Omega_u(\cala),
\Omega_u(\cala)]$ (which is generally not a graded algebra) is
sometimes called {\sl the noncommutative de Rham complex} \cite{ps}. It is worth noticing that, for $\cala$ noncommutative, there
is no tensor product over $\cala$ between $\cala$-modules (i.e.
no analog of the tensor product of vector bundles) and that therefore the Grothendieck group $K_0(\cala)$ (of classes of projective $\cala$-modules) has no product. Thus for $\cala$ noncommutative $K_0(\cala)\otimes\mathbb C$ is not an algebra and therefore there is no reason for a cochain complex such that its cohomology is a receptacle  for the image of the Chern character of $K_0(\cala)\otimes\mathbb C$ to be a graded algebra.\\

Also we did not describe here the approach to the differential
calculus and to the metric aspects in noncommutative geometry
based on generalized Dirac operators (spectral triples)
\cite{connes:03}, \cite{connes:04}, \cite{connes:05} as well as
the related supersymmetric approach of \cite{fgr}, see in
O.~Grandjean's lectures. In these notes we did not introduce specifically generalizations of linear connections and a fortiori not generalizations of riemannian structures.\\

Finally we did not discuss differential calculus for quantum
groups, i.e. bicovariant differential calculus \cite{slw}. In
the spirit of Section 2, let us define a {\sl graded
differential Hopf algebra} to be a graded differential algebra
$\fraca$ which is also a graded Hopf algebra with coproduct $\Delta$
such that $\Delta:\fraca\rightarrow \fraca\otimes \fraca$ is a
homomorphism of graded differential algebras (i.e. in
particular the differential $d$ of $\fraca$ satisfies the
graded co-Leibniz rule), with counit $\varepsilon$ such that
$\varepsilon\circ d=0$ and with antipode $S$ homogeneous of
degree 0 such that $S\circ d=d\circ S$. If $\fraca$ is a graded
differential Hopf algebra, then the subalgebra $\fraca^0$ of
elements of degree 0 of $\fraca$ is an ordinary Hopf algebra,
i.e. a quantum group, and $\fraca$ is a bicovariant
differential calculus over $\fraca^0$. Notice that if $G$ is a
Lie group then the graded differential algebra $\Omega(G)$ of
differential forms on $G$ is in fact a graded differential Hopf
algebra which is graded commutative, (in order to be  correct, one has to complete the tensor product in the definition of the coproduct or to use, instead of $\Omega(G)$, the graded differential subalgebra of forms generated by the representative functions on $G$).

\end{document}